\documentclass[11pt]{amsart}

\usepackage{amsmath,amssymb,amsxtra,amsthm,txfonts}
\usepackage{color,bm}
\usepackage{graphicx}
\usepackage{times}
\usepackage{tikz}
\usepackage{cases}
\usepackage{changes}

\usepackage{marginnote}
\usepackage{mathrsfs}
\usepackage{amsfonts}
\usepackage{enumerate}
\usepackage{hyperref}
\usepackage{pifont}

\hypersetup{
 pdftitle={Sharp endpoint Harnack inequality for nonlocal equations with logarithmic kernels},
 pdfauthor={Liguang Liu and Peilun Tian},
 pdfsubject={Same-ball endpoint weak Harnack inequalities},
 pdfkeywords={nonlocal equations, logarithmic kernels, weak Harnack inequality, Lorentz-Zygmund spaces},
 colorlinks = true,
 citecolor = blue,
 linkcolor = blue,
 urlcolor = blue
}

\makeatletter
\DeclareRobustCommand{\HarnackOld}{\color{red}\def\@linkcolor{red}\def\@citecolor{red}}
\DeclareRobustCommand{\HarnackNew}{\color{blue}\def\@linkcolor{blue}\def\@citecolor{blue}}
\makeatother
\makeatletter
\@namedef{subjclassname@1991}{Subject}
\@namedef{subjclassname@2000}{Subject}
\@namedef{subjclassname@2010}{Subject}
\@namedef{subjclassname@2020}{{\rm 2020} Mathematics Subject Classification}
\makeatother

\allowdisplaybreaks[2]
\AtBeginDocument{%
  \setlength{\parskip}{3pt plus 1pt minus .5pt}%
  \setlength{\abovedisplayskip}{10pt plus 2pt minus 2pt}%
  \setlength{\belowdisplayskip}{10pt plus 2pt minus 2pt}%
  \setlength{\abovedisplayshortskip}{5pt plus 2pt}%
  \setlength{\belowdisplayshortskip}{7pt plus 2pt minus 2pt}%
  \setlength{\jot}{6pt}%
}

\makeatletter
\def\ps@headings{%
  \ps@empty
  \def\@evenhead{%
    \setTrue{runhead}%
    \normalfont\scriptsize
    \thepage\hfil
    \def\thanks{\protect\thanks@warning}%
    \leftmark{}{}%
  }%
  \def\@oddhead{%
    \setTrue{runhead}%
    \normalfont\scriptsize
    \def\thanks{\protect\thanks@warning}%
    \rightmark{}{}\hfil\thepage
  }%
  \let\@mkboth\markboth
}
\makeatother
\newcommand{\rn}{\mathbb R^n}

\newtheorem{theorem}{Theorem}[section]
\newtheorem{lemma}[theorem]{Lemma}
\newtheorem{corollary}[theorem]{Corollary}
\newtheorem{proposition}[theorem]{Proposition}

\theoremstyle{definition}
\newtheorem{remark}[theorem]{Remark}
\newtheorem{definition}[theorem]{Definition}

\numberwithin{equation}{section}

\allowdisplaybreaks

\begin{document}

\title[Sharp Endpoint Weak Harnack Inequality]
{Sharp Endpoint Weak Harnack Inequality for Nonlocal Equations with Nearly Critical Kernels}

\author{Liguang Liu}
\address{School of Mathematics,
Renmin University of China,
Beijing 100872, People's Republic of China}
\email{liuliguang@ruc.edu.cn}

\author{Peilun Tian}
\address{School of Mathematics,
Renmin University of China,
Beijing 100872, People's Republic of China}
\email{peiluntian@ruc.edu.cn}

\thanks{L. Liu is supported by the National Natural Science Foundation of China (\# 12371102).}




\begin{abstract}
For $p\in (1,n)$, we establish endpoint weak Harnack inequalities
for weak supersolutions of a class of nonlocal quasilinear
equations with kernels of nearly critical order.
We consider symmetric kernels $K_{\phi,\,p}(x,y)$ comparable
to $\phi(|x-y|)/|x-y|^{n+p}$, where $\phi$ is quasi-increasing
and satisfies a mild upper scaling condition and suitable
integrability assumptions.

Writing $q_{\mathrm H}=n(p-1)/(n-p)$, we prove a weak
Harnack inequality at the critical exponent $q_{\mathrm H}$
in terms of decreasing rearrangements, with an explicit
correction determined by $\phi$ and a nonlocal tail term
accounting for the exterior values of the supersolution.
The proofs combine variational methods, a critical Sobolev
inequality in rearrangement form, and nonlocal truncation
energy estimates with linear dependence on the truncation
height.

A typical example is given by logarithmic kernels
$K_{\gamma,\,p}$, with $\gamma\in(1,\infty)$, satisfying
\[
K_{\gamma,\,p}(x,y)\approx
\frac{1}{|x-y|^{n+p}}
\left(\log\frac{e}{|x-y|\wedge1}\right)^{-\gamma}.
\]
For these kernels,  the general theorem yields endpoint weak
Harnack estimates in the Lorentz--Zygmund space
$L^{q_{\mathrm H},\infty}(\log L)^{-\beta_\gamma}$,
where $\beta_\gamma=(\gamma-1)/(p-1)$, as well as a
logarithmically weighted integral estimate with
integrability threshold $1+\beta_\gamma q_{\mathrm H}$.
We further prove that the critical exponent $q_{\mathrm H}$,
the Lorentz summability index $\infty$,
the logarithmic correction exponent $\beta_\gamma$,
and the integrability threshold $1+\beta_\gamma q_{\mathrm H}$ are all optimal.
The sharpness results are established by constructing weak
supersolutions through a nonlocal obstacle problem.
\end{abstract}

\maketitle

\vspace*{-1.3\baselineskip}
\begingroup
\setlength{\parskip}{0pt}
\linespread{1}\selectfont
\tableofcontents
\endgroup

\arraycolsep=0pt
\numberwithin{equation}{section}

\section{Introduction}

\subsection{Harnack inequalities for local operators}
The classical Harnack inequality (see, for example, \cite{GilbargTrudinger2001}) states that if $u\ge0$ is harmonic in a ball $B_{2r}\subset\rn$ of radius $2r$, then
\[
  \operatorname*{ess\,sup}_{B_r}u
  \le C_n\operatorname*{ess\,inf}_{B_r}u,
\]
where $C_n$ is a dimensional constant independent of $u$ and $r$.
The Harnack inequality holds with a constant depending on $n$ and $\Lambda_0/\lambda_0$ for nonnegative weak solutions of uniformly elliptic equations in divergence form
\[
  -\operatorname{div}(A\nabla u)=0\quad\text{in }\, B_{2r},
\]
where $A$ is a measurable matrix-valued function satisfying the ellipticity and boundedness conditions
\[
\lambda_0|\xi|^2\le  A(x)\xi\cdot\xi,\qquad |A(x)\xi|\le\Lambda_0|\xi|
\]
 for a.e. $x\in B_{2r}$ and all $\xi\in\rn$ and for some constants $0<\lambda_0\le\Lambda_0<\infty$;  this was proved by Moser \cite{Moser1961} via an iteration argument. These results form part of the classical regularity theory for elliptic equations,  whose foundations were laid by De Giorgi and Nash \cite{DeGiorgi1957, Nash1958}.
Moser's iteration  in \cite{Moser1961} also yields, for nonnegative weak supersolutions, the corresponding \emph{weak Harnack inequality}
\begin{align}\label{eq-wHI-Laplace}
  \left(\fint_{B_r}u^q\,dx\right)^{\frac 1 q}
  \le C\operatorname*{ess\,inf}_{B_r}u
\end{align}
for some $q>0$ depending only on the dimension and the ellipticity constants. In the Laplace case, when $n>2$, the estimate \eqref{eq-wHI-Laplace} holds throughout the range $0<q<n/(n-2)$; see also \cite{Trudinger1967, GilbargTrudinger2001}.

For the local quasilinear operator, the \emph{$p$-Laplacian}
\[
  (-\Delta)_p u:=-\operatorname{div}
       (|\nabla u|^{p-2}\nabla u),
\]
Trudinger \cite{Trudinger1967} established the corresponding Harnack theory, which was subsequently developed for quasiminima and related variational problems by DiBenedetto and Trudinger \cite{DiBenedettoTrudinger1984}.
 In particular, fix $p\in(1,n)$ and set
 \begin{align}\label{eq-qH}
 q_{\mathrm H}:=\frac{n(p-1)}{n-p}.
\end{align}
If $u\ge0$ is a weak supersolution of
$(-\Delta)_p u=0$ in $B_{2r}$, and
$0<q<q_{\mathrm H}$,
then
\begin{align}\label{eq-wHI-p-Laplace}
  \left(\fint_{B_r}u^q\,dx\right)^{\frac{1}{q}}
 \le C(n,p,q)\operatorname*{ess\,inf}_{B_r}u.
\end{align}
The exponent $q_{\mathrm H}$ is the critical index for this inequality; indeed, it is the integrability threshold for the $p$-Laplace fundamental solution $|x|^{(p-n)/(p-1)}$, and the estimate \eqref{eq-wHI-p-Laplace}  fails for $q\ge q_{\mathrm H}$.

\subsection{A general nonlocal quasilinear operator}
Let $p\in(1,\infty)$.
Throughout this paper, we work with  a measurable function \(\phi:\,(0,\infty)\to[0,\infty)\) that is positive on a set of positive measure and satisfies
\begin{equation}\label{eq-adm}
 \int_0^\infty\min\{1,\, t^{-p}\}\,\phi(t)\,\frac{dt}{t}<\infty.
\end{equation}
For any $r\in (0,\infty)$, define
\begin{equation}\label{eq-Phi}
\Phi(r):=\int_0^\infty
\phi(r\tau)\min\{1,\,\tau^{-p}\}\,\frac{d\tau}{\tau}.
\end{equation}
By \eqref{eq-adm}, the function $\Phi$ is well-defined
and finite for every $r>0$.
Further properties of $\Phi$ will be established in
 Lemma~\ref{lem-basic-nearly critical phi}.
Associated with \(\phi\), the nonlocal quasilinear operator has the formal expression
\begin{align}\label{eq-Lphi}
 \mathcal L_{\phi,\,p}u(x)
 :=\mathrm{p.v.}\int_{\mathbb R^n}
 |u(x)-u(y)|^{p-2}(u(x)-u(y))\,K_{\phi,\,p}(x,y)\,dy,
\end{align}
 where the kernel \(K_{\phi,\,p}:\,\mathbb R^n\times\mathbb R^n\to[0,\infty)\) is a measurable, symmetric kernel such that
\begin{align}\label{eq-ker-phi}
 \Lambda^{-1}\frac{\phi(|x-y|)}{|x-y|^{n+p}}
 \le K_{\phi,\,p}(x,y)\le
 \Lambda\frac{\phi(|x-y|)}{|x-y|^{n+p}}
\end{align}
 for almost all distinct \(x, y\in\rn\) and some  fixed constant
 \(\Lambda\in(1,\infty)\).
 The operator $\mathcal L_{\phi,\,p}$ in \eqref{eq-Lphi} arises from
 the first variation of the following Besov type energy
 $$[f]^p_{{\mathsf \Lambda}^{\phi,\,p}(\rn)}:= \int_{\rn}\int_{\rn}
 \frac{|f(x)-f(y)|^p}{|x-y|^{n+p}}\,\phi(|x-y|)\,dy\,dx.$$
Define the space
 \[
 {\mathsf \Lambda}^{\phi,\,p}(\rn)
 :=\left\{f\in L^p(\rn):\, [f]^p_{ {\mathsf \Lambda}^{\phi,\,p}(\rn)}<\infty\right\},
\]
which is equipped with the norm
\[\|f\|_{{\mathsf \Lambda}^{\phi,\,p}(\rn)}:=\|f\|_{L^p(\rn)}+[f]_{{\mathsf \Lambda}^{\phi,\,p}(\rn)}.\]
The condition \eqref{eq-adm} on \(\phi\) ensures that
$
C_c^\infty(\rn)\subset{\mathsf \Lambda}^{\phi,\,p}(\rn),
$
so that the space \({\mathsf \Lambda}^{\phi,\,p}(\rn)\) is non-trivial.
Indeed, for \(f\in C_c^\infty(\rn)\), we have
\[
 \|f(\cdot+h)-f\|_{L^p(\rn)}
 \le\min\left\{|h| \sum_{j=1}^n\|\partial_j f\|_{L^p(\rn)},\,2\|f\|_{L^p(\rn)}\right\},
\]
which, together with polar coordinates and \eqref{eq-adm}, implies
\begin{align*}
 &\int_{\mathbb R^n}\int_{\mathbb R^n}
 \frac{|f(x)-f(y)|^p}{|x-y|^{n+p}}\,\phi(|x-y|)\,dy\,dx\\
 &\quad=  \int_{\mathbb R^n}  \|f(\cdot+h)-f\|_{L^p(\rn)}^p\frac{\phi(|h|)}{|h|^{n+p}}\,dh\notag\\
 &\quad\le \omega_{n-1}\left(\sum_{j=1}^n\|\partial_j f\|_{L^p(\rn)}+2\|f\|_{L^p(\rn)}\right)^p
 \int_0^\infty\min\{t^p,\,1\}\,\phi(t)\,t^{-p}\,\frac{dt}{t} \notag\\
 &\quad< \infty. \notag
\end{align*}
Moreover, condition \eqref{eq-adm} is both sufficient and necessary for \(C_c^\infty(\rn)\subset{\mathsf \Lambda}^{\phi,\,p}(\rn)\); see \cite{LiuZhangHigherBBMMS}.

A basic example satisfying \eqref{eq-adm} is
$\phi(t)=t^{(1-s)p}$, where $s\in(0,1)$.
In this case, the general framework in
\eqref{eq-Lphi}--\eqref{eq-ker-phi}
reduces to the fractional $(s,p)$-model, with the kernel $K_{\phi,\,p}$
satisfying
\[
 \Lambda^{-1}\frac{1}{|x-y|^{n+sp}}
 \le K_{\phi,\,p}(x,y)\le
 \Lambda\frac{1}{|x-y|^{n+sp}}.
\]
Its prototype is the \emph{fractional $(s,p)$-Laplace
operator}
\[
(-\Delta)_{s,\,p}u(x)
:=\mathrm{p.v.}\int_{\mathbb R^n}
\frac{|u(x)-u(y)|^{p-2}(u(x)-u(y))}
{|x-y|^{n+sp}}\,dy.
\]
The corresponding energy space is the known \emph{fractional Sobolev space \(W^{s,\,p}(\rn)\)}, consisting of all functions \(f\in L^p(\rn)\)
such that
\[
 [f]_{W^{s,\,p}(\rn)}
   =\left(\int_{\rn}\int_{\rn}\frac{|f(x)-f(y)|^p}{|x-y|^{n+sp}}\,dy\,dx\right)^{\frac 1p}<\infty.
\]
Unlike in the local case, Harnack estimates for such
nonlocal equations involve contributions from the
exterior values of the solution, measured by the
nonlocal tail
\[
 \mathsf{Tail}_{(s,\,p)}(v;x_0,R)
 :=\left(R^{sp}\int_{\mathbb R^n\setminus B_R(x_0)}
       \frac{|v(y)|^{p-1}}{|x_0-y|^{n+sp}}\,dy
   \right)^{\frac{1}{p-1}}.
\]
For this fractional model,  Di Castro, Kuusi, and
Palatucci \cite{DiCastro-Kuusi-Palatucci}
(see also \cite{KorvenpaaKuusiPalatucci2017MathAnn} and, for $p=2$, \cite{Kassmann2007, Kassmann2009, Kassmann2011})
established Harnack and weak Harnack inequalities
involving the nonlocal tail.
More precisely, if $u$ is a weak solution
of $\mathcal L_{\phi,\,p}u=0$ in the ball $B_R(x_0)$  and
$u\ge0$ a.e. on $B_R(x_0)$, then for every $0<r\le R/2$ the \emph{Harnack inequality}
\[
 \operatorname*{ess\,sup}_{B_r(x_0)}u
 \le C\operatorname*{ess\,inf}_{B_r(x_0)}u
 +C\left(\frac rR\right)^{\frac{sp}{p-1}}
      \mathsf{Tail}_{(s,\,p)}(u_-;x_0,R)
\]
holds, where $u_-:=\max\{-u,0\}$. If $u$ is merely a weak supersolution, then the
\emph{weak Harnack inequality}
\[
 \left(\fint_{B_r(x_0)}u^q\,dx\right)^{\frac{1}{q}}
 \le C\operatorname*{ess\,inf}_{B_{2r}(x_0)}u
 +C\left(\frac rR\right)^{\frac{sp}{p-1}}
      \mathsf{Tail}_{(s,\,p)}(u_-;x_0,R)
\]
 holds for $0<r\le R/4$ whenever
\[ 0<q<\begin{cases}
\frac{n(p-1)}{n-sp}&\quad \text{if}\ \,sp<n;\\
+\infty &\quad \text{if}\ \,sp\ge n.
\end{cases}\]
The constants $C$ here depend only on $n,p,s,\Lambda$,
and additionally on $q$ in the weak Harnack inequality.
Note that, when $1<p<n$, the fractional
critical exponent $n(p-1)/(n-sp)$ converges to the
local critical exponent
$q_{\mathrm H}=n(p-1)/(n-p)$ as $s\uparrow1$.

Of particular interest in this paper is the logarithmic
kernel
\begin{equation}\label{eq-phi-gamma}
\phi_\gamma(t):=
\left(\log\frac{e}{t\wedge1}\right)^{-\gamma}
\qquad \text{for all}\ \,t\in(0,\infty),
\end{equation}
which satisfies \eqref{eq-adm} only when $\gamma>1$.
For this choice, we write
$$\mathcal L_{\gamma,\,p}:=\mathcal L_{\phi_\gamma,\,p},
\qquad
K_{\gamma,\,p}:=K_{\phi_\gamma,\,p}.$$
In this case, the energy space associated with $\mathcal L_{\gamma,\,p}$
is the well-known \emph{logarithmic Besov space}
${\mathsf\Lambda}_p^\gamma(\mathbb R^n)$, which consists of all functions $f\in L^p(\mathbb R^n)$ such that
\begin{align*}
 [f]_{{\mathsf\Lambda}_p^\gamma(\mathbb R^n)}:=\left(\int_{\mathbb R^n}\int_{\mathbb R^n}
 \frac{|f(x)-f(y)|^p}{|x-y|^{n+p}}
 \left(\log\frac{e}{|x-y|\wedge1}\right)^{-\gamma}\,dy\,dx\right)^\frac 1p<\infty.
\end{align*}
The space ${\mathsf\Lambda}_p^\gamma(\rn)$ has smoothness
slightly below order one, with a logarithmic deficit
controlled by $\gamma$.
For any $s\in(0,1)$ and $\gamma>1$, we have continuous embeddings
\[
W^{1,\,p}(\rn)\hookrightarrow{\mathsf\Lambda}_p^\gamma(\rn)
\hookrightarrow W^{s,\,p}(\rn).
\]
Function spaces with logarithmic smoothness are closely related to the work of Br\'ezis and Wainger \cite{BrezisWainger1980} on logarithmic refinements of
Sobolev embeddings.
Such logarithmic function spaces were systematically studied by Edmunds and Haroske (see, for example, \cite{EdmundsHaroske1999, Haroske2000, Edmunds-Haroske-JAT-2000}) and others (see \cite{Farkas-Leopold-2006-AnnMatPuraAppl, Caetano-Gogatishvili-Opic-JAT-2008, Caetano-Gogatishvili-Opic-JAT-2011}) in an even more general scale, and have recently attracted intensive study in \cite{Cobos-Dominguez-Triebel-JFA-2016, CobosDominguez2015, DHT,DT2023}, and in \cite{LiuXiaoWuYuan2021AdvMath, DLX2026} from the perspective of logarithmic capacity theory.

Another motivation comes from the
Bourgain--Brezis--Mironescu (BBM) formula.
The classical BBM formula \cite{BBM2001, BBM2002} states that,
for every $f\in W^{1,\,p}(\rn)$,
\[
 \lim_{s\uparrow1}(1-s)[f]_{W^{s,\,p}(\mathbb R^n)}^p
   =C_{n,p}\int_{\mathbb R^n}|\nabla f|^p\,dx,
\]
where \(C_{n,p}\) is a positive constant depending only on $n$ and $p$.
It turns out that an analogous BBM-type formula remains valid for the space \({\mathsf\Lambda}_p^\gamma(\mathbb R^n)\) (see \cite{LiuZhangHigherBBMMS} or \cite{Foghem2025NoDEA}) in the following form: for every $f\in W^{1,\,p}(\rn)$,
\[
 \lim_{\gamma\downarrow1}(\gamma-1)[f]_{{\mathsf\Lambda}_p^\gamma(\mathbb R^n)}^p
   =C_{n,p}\int_{\mathbb R^n}|\nabla f|^p\,dx.
\]
This shows that the operator $\mathcal L_{\gamma,\,p}$
is naturally connected with the local $p$-Laplacian,
much like the fractional $(s,p)$-Laplacian.

These connections motivate the study of weak Harnack inequalities for the
general operator $\mathcal L_{\phi,\,p}$ defined in
\eqref{eq-Lphi}--\eqref{eq-ker-phi}, and in particular the study of such
inequalities at the critical exponent
\[
q_{\mathrm H}=\frac{n(p-1)}{n-p}
\]
for kernels $\phi$ behaving like $\phi_\gamma$. Note that $q_{\rm H}$ is
exactly the critical weak Harnack exponent for the local $p$-Laplacian $(-\Delta)_p$ (see \eqref{eq-qH} ).

\subsection{Sharp endpoint weak   Harnack inequalities}
Motivated by the logarithmic kernel $\phi_\gamma$ in \eqref{eq-phi-gamma},
we introduce a broader class of nearly critical kernels
satisfying the following structural conditions.

\begin{definition}\label{def-ker-critical}
A function $\phi:\, (0,\infty)\to (0,\infty)$ is called
a \emph{nearly critical kernel} if it satisfies
the following conditions:
\begin{enumerate}[\rm (i)]
\item $\phi$ admits a continuous extension to
$[0,\infty)$ with $\phi(0)=0$, and satisfies the integrability condition \eqref{eq-adm}.

\item $\phi$ is \emph{quasi-increasing} on $(0,\infty)$, namely, there exists
a constant $C>0$ such that
\[
\phi(\tau)\le C\phi(t)
\qquad\text{for all }0<\tau\le t<\infty.
\]

\item $\phi$ satisfies the following \emph{mild upper scaling}
condition: for every $\varepsilon\in(0,1)$,
there exists a constant $C_\varepsilon>0$ such that
\begin{equation}\label{condition2}
\phi(t)\le C_\varepsilon
\left(\frac{t}{\tau}\right)^\varepsilon\phi(\tau)
\qquad\text{for all }0<\tau\le t<\infty.
\end{equation}
\end{enumerate}
\end{definition}

Throughout the paper, unless otherwise stated,
$\Omega$ always denotes a nonempty open set with $\Omega\Subset\rn$. Associated with the kernel \(K_{\phi,\,p}\) in \eqref{eq-ker-phi},  we shall use the following local energy
\begin{align}\label{corresponding energy}
\mathcal E_{\phi,\,p}(u;\,\Omega)
&:=\int_\Omega\int_\Omega
|u(x)-u(y)|^pK_{\phi,\,p}(x,y)\,dy\,dx\\
&\quad+2\int_\Omega\int_{\Omega^c}
|u(x)-u(y)|^pK_{\phi,\,p}(x,y)\,dy\,dx\notag\\
&=\iint_{\mathcal C_\Omega}
|u(x)-u(y)|^pK_{\phi,\,p}(x,y)\,dy\,dx, \notag
\end{align}
where
\begin{equation}\label{eq-C_Omega}
\mathcal C_\Omega
:= (\mathbb R^n\times\mathbb R^n)
\setminus
\big((\mathbb R^n\setminus\Omega)
\times(\mathbb R^n\setminus\Omega)\big).
\end{equation}
It is obvious that if $u=0$ a.e. outside of $\Omega$ then
$$
\mathcal E_{\phi,\,p}(u;\,\Omega)=\int_{\rn}\int_{\rn}
|u(x)-u(y)|^pK_{\phi,\,p}(x,y)\,dy\,dx.
$$
Set
\begin{align*}
\mathcal F^{\phi,\,p}(\Omega)
:=\left\{
f:\ \rn\to\mathbb R\ \text{measurable} :\ f|_\Omega\in L^p(\Omega),\ \mathcal E_{\phi,\,p}(f;\,\Omega)<\infty
\right\}.
\end{align*}
To study the weak Harnack inequality for \(\mathcal L_{\phi,\,p}\), we need the notion of the weak supersolution.

\begin{definition}\label{def-supersol}
Let \(p\in(1,\infty)\). A function \(u\in \mathcal F^{\phi,\,p}(\Omega)\) is called a \emph{weak supersolution} to
$$
\mathcal L_{\phi,\,p}u=0 \quad \text{in }\, \Omega
$$
if for every nonnegative \(\zeta\in \mathcal F^{\phi,\, p}(\Omega)\) satisfying \(\zeta=0\) a.e. in \(\rn \setminus \Omega\),
\begin{align*}
 \iint_{\mathcal C_\Omega}
 |u(x)-u(y)|^{p-2}(u(x)-u(y))
 (\zeta(x)-\zeta(y))K_{\phi,\, p}(x,y)\,dy\,dx\ge0.
\end{align*}

\end{definition}

Now, we define the corresponding tail spaces.

\begin{definition}\label{Tail}
Let \(p\in(1,\infty)\) and $\phi:\, (0,\infty)\to (0,\infty)$ satisfy \eqref{eq-adm}. Associated to $\phi$, define \(\Phi\) as in \eqref{eq-Phi}.
The \emph{tail space} is defined by
\begin{equation}\label{Def-Tail space}
 \mathbb L^{\phi}_{p}(\mathbb R^n):=
 \left\{u\in L^{p-1}_{\mathrm{loc}}(\mathbb R^n):\
 \int_{\mathbb R^n}
 \frac{|u(x)|^{p-1}}{(1+|x|)^{n+p}}\phi(|x|)\,dx<\infty\right\}.
\end{equation}
For \(u\in\mathbb L^{\phi}_{p}(\mathbb R^n)\), \(x_0\in\mathbb R^n\) and \(r>0\),
set
\begin{equation}\label{Def-tail}
 \mathsf{Tail}^{\phi}_{p}(u;x_0,r):=
 \left(\frac{r^p}{\Phi(r)}
 \int_{\mathbb R^n\setminus B_r(x_0)}
          \frac{|u(y)|^{p-1}}{|x_0-y|^{n+p}}\,\phi(|x_0-y|)\,dy\right)^{\frac{1}{p-1}}.
\end{equation}


\end{definition}

For a measurable function $f$ on a ball $B\subset\rn$, let $f_B^*$ denote
the \textit{nonincreasing rearrangement} of $f$, that is,
\begin{align}\label{eq-f_B*}
 f_B^*(t)=\inf\left\{\lambda\ge0:\
       \frac{|\{x\in B:|f(x)|>\lambda\}|}{|B|}\le t\right\}, \qquad t\in(0,\infty).
\end{align}
Note that $f_B^*(t)=0$ for all $t\in[1,\infty)$.
The main result of this paper is the following \emph{endpoint weak Harnack inequality}.

\begin{theorem}
\label{thm-wHI-qH-phi}
Let \(p\in(1,n)\) and \(q_{\rm H}\) be as in \eqref{eq-qH}. For a nearly
critical kernel \(\phi\), define the associated function \(\Phi\) as in
\eqref{eq-Phi}.
If  \(u \in \mathcal F^{\phi,\, p}(\Omega)\) is a  weak supersolution to
\begin{align}\label{eq-Lgammap-phi}
\mathcal L_{\phi,\,p}u=0 \quad \text{in }\, \Omega,
\end{align}
 and $u\ge 0$ a.e. on a ball \(B_{R}(x_0) \Subset \Omega\) with $x_0\in\rn$ and $R>0$,
then for any  \(0<r\le R/2\),
\begin{equation}\label{intro:endpoint-phi}
\begin{aligned}
 &\sup_{t \in (0,\infty)}t^{\frac{1}{q_{\mathrm H}}}
 \left(\frac{\Phi(rt^{1/n})}{\Phi(r)}\right)^{\frac{1}{p-1}}
                         u_{B_r(x_0)}^*(t)\\
 &\qquad\le C\left(\operatorname*{ess\,inf}_{B_r(x_0)}u
 +\left(\frac{r^p\Phi(R)}{R^p\Phi(r)}\right)^{\frac{1}{p-1}}
                         \mathsf{Tail}^{\phi}_{p}(u_-;x_0,R)\right),
\end{aligned}
\end{equation}
where $C=C(n,p,\phi,\Lambda)$ is independent of $x_0, r, R$, $\Omega$ and  $u$.
\end{theorem}

In what follows, to simplify notation, we write
\begin{align}\label{eq-Uu}
\mathcal U_{\phi,\, p}(u,x_0,r,R):=\operatorname*{ess\,inf}_{B_r(x_0)}u
 +\left(\frac{r^p\Phi(R)}{R^p\Phi(r)}\right)^{\frac{1}{p-1}}
                         \mathsf{Tail}^{\phi}_{p}(u_-;x_0,R).
\end{align}

For $0<q<\infty$, $0<s\le \infty$ and $\sigma\in\mathbb R$,
we define the \emph{Lorentz-Zygmund space $L^{q,s}(\log L)^{\sigma}(B)$} consisting of all measurable functions $f$ on $B$ such that
\begin{align}\label{eq-LqsLogL}
 \|f\|_{L^{q,s}(\log L)^{\sigma}(B)}
 :=\left(\int_0^1
 \left[t^{\frac 1 q} [\log(e/t)]^{\sigma}f_B^*(t)\right]^s\,
                                  \frac{dt}{t}\right)^{\frac 1 s}<\infty
\end{align}
with the usual modification when \(s=\infty\), that is,
 \begin{align}\label{eq-LogL-fz}
 \|f\|_{L^{q,\infty}(\log L)^{\sigma}(B)}
 :=\sup_{t\in(0,1)}
 t^{\frac 1 q} [\log(e/t)]^{\sigma}f_B^*(t)<\infty.
\end{align}

Since the logarithmic kernel $\phi_\gamma$ in
\eqref{eq-phi-gamma} is nearly critical for every
$\gamma\in(1,\infty)$, Theorem~\ref{thm-wHI-qH-phi} applies to
weak supersolutions of
$\mathcal L_{\gamma,\,p}u=0$ in $\Omega$.
The following theorem gives \textit{sharp endpoint} weak Harnack
estimates for the logarithmic model
$\mathcal L_{\gamma,\,p}$.

\begin{theorem}\label{thm-wHI-qH}
Let \(p\in(1,n)\), $\gamma\in(1,\infty)$ and \(q_{\rm H}\) be as in \eqref{eq-qH}.
Define $\phi_\gamma$ as in \eqref{eq-phi-gamma}, and
\begin{align}\label{eq-Phi-gamma}
\Phi_\gamma(r):=\int_0^1 \phi_\gamma(rt)\,\frac{dt}{t}\qquad \text{for all}\ \, r\in(0,\infty).
\end{align}
Suppose that \(u \in \mathcal F^{\phi_\gamma,\, p}(\Omega)\) is a  weak supersolution to
\begin{align}\label{eq-Lgammap}
\mathcal L_{\gamma,\,p}u=0 \quad \text{in }\, \Omega
\end{align}
and $u\ge 0$ a.e. on a ball \(B_{R}(x_0) \Subset \Omega\) with $x_0\in\rn$ and $R>0$.
Let $\mathcal U_{\gamma,\,p}$ be defined as in \eqref{eq-Uu} by taking
$\phi=\phi_\gamma$ and $\Phi=\Phi_\gamma$ therein.
Then, the following hold:

\begin{enumerate}[\rm (i)]

\item
For any  \(0<r\le R/2\), we have
\begin{equation}\label{intro:endpoint}
 \sup_{t \in (0,\infty)}t^{\frac{1}{q_{\mathrm H}}}
 \left(\frac{\Phi_\gamma(rt^{1/n})}{\Phi_\gamma(r)}\right)^{\frac{1}{p-1}}
                         u_{B_r(x_0)}^*(t)
                         \le C\,\mathcal U_{\gamma,\, p}.
\end{equation}

\item
Set
$ \beta_\gamma:=(\gamma-1)/(p-1).$
For any $0<r\le\min\{1,R/2\}$,
\begin{equation}\label{intro:LZ-functional}
 \|u\|_{L^{q_{\rm H},\infty}(\log L)^{-\beta_{\gamma}}(B_r(x_0))}
 \le C\mathcal U_{\gamma,\, p}.
\end{equation}

\item
For any
$0<r\le\min\{1,R/2\}$  and any
\begin{equation}\label{LH:integral-threshold}
 \lambda>1+\beta_\gamma q_{\mathrm H}
        =1+\frac{n(\gamma-1)}{n-p},
\end{equation}
 we have
\begin{equation}\label{LH:critical-integral}
 \left(\fint_{B_r(x_0)}
 \frac{u^{q_{\mathrm H}}}{[\log(e+u/\mathcal U_{\gamma,\, p})]^\lambda}\,dx
 \right)^{1/q_{\mathrm H}}
 \le C\, \mathcal U_{\gamma,\, p},
\end{equation}
with the convention that the left-hand side
is zero when $\mathcal U_{\gamma,\,p}=0$.
\end{enumerate}
The constants $C$ in (i) and (ii) depend only on
$n,p,\gamma,\Lambda$, while the constant in (iii)
may additionally depend on $\lambda$.
All these constants are independent of
$u$, $x_0$, $r$, $R$, and $\Omega$.
\end{theorem}

\begin{remark}
The general framework of $\mathcal L_{\phi,\,p}$ with
nearly critical kernels $\phi$ falls within the setting of
Ok and Song \cite{OkSong}. However, the endpoint weak
Harnack inequalities established in
Theorems~\ref{thm-wHI-qH-phi} and \ref{thm-wHI-qH}
are not direct consequences of \cite{OkSong}. Indeed, under the same hypotheses as in Theorem~\ref{thm-wHI-qH-phi}, applying \cite[Theorem~1.5]{OkSong} yields that if
$
0<q<
q_{\mathrm H}$,
then for every ball \(B_{R}(x_0) \Subset \Omega\) on which $u\ge0$ a.e. and every \(0<r\le R/2\), one has
\begin{equation}\label{eq-wHI-q}
 \left(\fint_{B_{r/2}(x_0)}u^q\,dx\right)^{1/q}
 \le C\,\mathcal U_{\phi,\, p}(u,x_0,r,R),
\end{equation}
where \(C=C(n,p, \phi,\Lambda,q)\) is independent of \(x_0\), \(r\), \(R\), and \(u\).
For a detailed discussion of \eqref{eq-wHI-q}, we refer the reader to Section~\ref{ss-OkSong} below.
Moreover, for $1<p<n$, we show in Section~\ref{sec:subcritical-consequence} that the endpoint
estimate \eqref{intro:endpoint-phi} implies the subcritical weak
Harnack inequality \eqref{eq-wHI-q} for every $0<q<q_{\mathrm H}$.
In particular, when $\phi=\phi_\gamma$, each  estimate
in Theorem \ref{thm-wHI-qH} also implies \eqref{eq-wHI-q}.
In view of these discussions, the main contribution of Theorems~\ref{thm-wHI-qH-phi} and \ref{thm-wHI-qH} lies in establishing the endpoint weak Harnack inequality at the critical exponent \(q=q_{\mathrm H}\), giving
a weighted rearrangement estimate at the endpoint of the subcritical range \(0<q<q_{\mathrm H}\) obtained in \cite{OkSong}.
\end{remark}

\begin{theorem}
\label{intro:sharpness}
The three weak Harnack inequalities in
Theorem~\ref{thm-wHI-qH} are all sharp in the following sense:
\begin{enumerate}[\rm (i)]

\item The exponent \(q_{\mathrm H}\) in \eqref{intro:endpoint} and \eqref{intro:LZ-functional} is optimal and cannot be replaced by any \(q>q_{\mathrm H}\).

\item The exponent \(1/(p-1)\) in \eqref{intro:endpoint} and the exponent \(\beta_{\gamma}\) in \eqref{intro:LZ-functional} are optimal, i.e., they cannot be replaced by any smaller values.

\item  In \eqref{intro:LZ-functional}, the Lorentz--Zygmund space \(L^{q_{\rm H},\infty}(\log L)^{-\beta_{\gamma}}(B_r(x_0))\) is optimal in the sense that its index \(\infty\) cannot be replaced by any finite \(s\in(0,\infty)\).

\item  The threshold \(1+\beta_\gamma q_{\mathrm H}\) for \(\lambda\) in \eqref{LH:integral-threshold} is optimal: the inequality \eqref{LH:critical-integral} fails if \(\lambda\le 1+\beta_\gamma q_{\mathrm H}\).
\end{enumerate}

\end{theorem}

The paper is organized as follows.
Section~\ref{sec2} is devoted to establishing the existence and uniqueness of solutions to the exterior Dirichlet problem and the obstacle problem for nonlocal operators \(\mathcal{L}_{\phi,p}\) with general kernels \(\phi\) satisfying only \eqref{eq-adm}.

Section~\ref{sec:logsp} develops the key estimates needed to establish Theorem~\ref{thm-wHI-qH-phi} and Theorem~\ref{thm-wHI-qH}. Indeed, in Sections~\ref{sec:kernel} and~\ref{sec-local-stab}, we establish basic estimates for the nearly critical kernels and some functional-analytic properties of the energy space, respectively. Moreover, Section~\ref{ss-f*-sobolev} establishes a critical Sobolev inequality on \(\mathbb{R}^n\) in rearrangement form.  In Section~\ref{ss-OkSong}, we apply \cite{OkSong} to obtain the weak Harnack inequality in the subcritical range, which is further used in Section~\ref{sec-linear-growth} to obtain an energy estimate for localized truncations that grows linearly in the truncation level.

The proofs of Theorems~\ref{thm-wHI-qH-phi} and~\ref{thm-wHI-qH} are completed in Section~\ref{sec-wHI}. In Section~\ref{ss-wHI-qH-phi}, we prove Theorem~\ref{thm-wHI-qH-phi} by combining the key estimates established in Sections~\ref{ss-f*-sobolev}, \ref{ss-OkSong}, and~\ref{sec-linear-growth}. In Section~\ref{ss-wHI-cor}, taking \(\phi\) to be the logarithmic kernel \(\phi_\gamma\), we prove Theorem~\ref{thm-wHI-qH} as a consequence of Theorem~\ref{thm-wHI-qH-phi}. Furthermore, in Section~\ref{sec:subcritical-consequence}, we show that the critical weak Harnack estimates in Theorems~\ref{thm-wHI-qH-phi} and~\ref{thm-wHI-qH} directly imply the subcritical weak Harnack inequality.

In Section~\ref{sec:proof-sharpness}, for the logarithmic kernel \(\phi_\gamma\), using the results from Section~\ref{sec2}, we construct a family of normalized obstacle minimizers and prove the sharpness assertions of Theorem~\ref{intro:sharpness}.

\bigskip

\noindent{\bf Notation.\,}
We use the following notation and notions throughout the whole paper.
\begin{itemize}

\item Let $\mathbb N:=\{1,2,\dots\}$ and $\mathbb Z=\{0,\pm 1,\pm2,\dots\}$.
\item For $t\in\mathbb R$ we write
$J_p(t)=|t|^{p-2}t$, with $J_p(0)=0$.

\item For any $a,b\in\mathbb R$, set $a\wedge b:=\min\{a,b\}$
and $a\vee b:=\max\{a,b\}$.

\item For a measurable function $f$, write $f_+=\max\{f,0\}$ and
$f_-=\max\{-f,0\}$.

\item For any measurable set $E\subset\rn$, let $|E|$ be its Lebesgue measure and $E^c:=\rn\setminus E$.
For a measurable set $E$ of positive finite
measure and $f\in L^1(E)$, let
$$f_E=\fint_E f\,dx=\frac{1}{|E|}\int_E f\,dx.$$

\item Let $B_r(x_0)$ denote the ball
centered at $x_0\in\rn$ with radius $r\in(0,\infty)$. We abbreviate $B_r:=B_r(x_0)$ when the
center is fixed.
Let $\mathbb{S}^{n-1} := \partial B_1(0)$ be the unit sphere in
$\rn$, and let $\omega_{n-1} := \sigma(\mathbb{S}^{n-1})$ be its surface area.
We write $\nu_n=|B_1(0)|$.

\item The letters $C$ and $c$ usually denote positive constants that are
independent of the main parameters involved. If a constant depends on
additional parameters, say $\alpha,\beta,\ldots$, we write
$C(\alpha,\beta,\ldots)$ or $C_{\alpha,\beta,\ldots}$, whose values may vary from line to line. The symbol
$u\lesssim v$ (resp., $u\gtrsim v$) means that $u\le Cv$ (resp., $u\ge cv$)
for a specified range of the variables. If $u\lesssim v\lesssim u$, then we
write $u\approx v$.

\end{itemize}

\section{The weak solution in a general framework}\label{sec2}

In this section, we study the exterior Dirichlet problem and the obstacle problem for nonlocal equations with general admissible kernels. The existence and uniqueness results are of independent interest, since they apply to the general class of kernels satisfying \eqref{eq-adm}. We also show that every obstacle minimizer is a weak supersolution. In Section~\ref{sec:proof-sharpness} below, for the logarithmic kernel $\phi_\gamma$, we use this property to construct a family of weak supersolutions that establishes the sharpness of the indices in our endpoint weak Harnack inequalities.

\subsection{Existence and uniqueness for the Dirichlet problem}\label{ss-DP}

\begin{lemma}\label{lem:bounded-support-difference}
Let \(1\leq p<\infty\) and $D:=\operatorname{diam} \Omega$.
If \(f\in L^p(\mathbb R^n)\) and \(f=0\) a.e. in \(\mathbb R^n\setminus\Omega\), then
\[
|h|\,\|f\|_{L^p(\mathbb R^n)}
\leq 2 D\,\|\Delta_h f\|_{L^p(\mathbb R^n)}
\]
for every \(h\in\mathbb R^n\) with \(0<|h|\leq D\),
where \(\Delta_h f(x):=f(x+h)-f(x)\).
\end{lemma}

\begin{proof}
Since  \(\Omega\Subset\mathbb R^n\) is a nonempty precompact open set, it follows that $D\in(0,\infty)$.
Fix \(h\in\mathbb R^n\) with \(0<|h|\leq D\), and set
\[
N:=\left\lfloor\frac{D}{|h|}\right\rfloor+1.
\]
Then
\[
D<N|h|\leq D+|h|\leq 2D.
\]
Since \(f=0\) a.e. outside \(\Omega\) and
\(\operatorname{diam}\Omega=D\), we have
\[
f(x)f(x+Nh)=0
\quad\text{for a.e. }x\in\mathbb R^n.
\]
Consequently,
\begin{align}\label{eq-Delta_hf1}
\|f\|_{L^p(\mathbb R^n)}
\leq
\|f(\cdot+Nh)-f(\cdot)\|_{L^p(\mathbb R^n)}.
\end{align}
On the other hand, the telescoping identity gives
\[
f(x+Nh)-f(x)
=\sum_{k=0}^{N-1}\Delta_hf(x+kh).
\]
By Minkowski's inequality and the translation invariance
of the \(L^p\) norm,
\begin{align}\label{eq-Delta_hf2}
\|f(\cdot+Nh)-f(\cdot)\|_{L^p(\mathbb R^n)}
&\leq \sum_{k=0}^{N-1}
\|\Delta_hf(\cdot+kh)\|_{L^p(\mathbb R^n)}\\
&=N\|\Delta_hf\|_{L^p(\mathbb R^n)}\notag\\
&\leq \frac{2D}{|h|}
\|\Delta_hf\|_{L^p(\mathbb R^n)}. \notag
\end{align}
Combining \eqref{eq-Delta_hf1} and \eqref{eq-Delta_hf2} completes the proof.
\end{proof}

\begin{lemma}\label{general:nonlocal-Poincare}
Let \(1\leq p<\infty\) and $\phi$ satisfy \eqref{eq-adm}.
Then, for any $f\in {\mathsf \Lambda}^{\phi,\, p}(\rn)$, with $f=0$ a.e. in $\rn \setminus  \Omega$,
\begin{equation}\label{general:Poincare}
  \|f\|_{L^p(\mathbb R^n)}
  \le C [f]_{{\mathsf \Lambda}^{\phi,\, p}(\rn)},
\end{equation}
where $C$ is a positive constant  depending only on $n$, $p$, $\Omega$ and $\phi$.
\end{lemma}

\begin{proof} Let $D:=\operatorname{diam} \Omega$.
By Lemma \ref{lem:bounded-support-difference}, for every
$h\in\mathbb R^n$ with $0<|h|\leq D$, we have
\[
|h|\,\|f\|_{L^p(\mathbb R^n)}
\leq 2D\,\|\Delta_hf\|_{L^p(\mathbb R^n)}.
\]
If $|h|>D$, then the supports of
$f$ and $f(\cdot+h)$ are disjoint, which implies
\[
  \|\Delta_h f\|_{L^p(\mathbb R^n)}^p
  =2\|f\|_{L^p(\mathbb R^n)}^p.
\]
From these last two formulas, together with the Tonelli theorem and polar coordinates,
we deduce
\begin{align*}
[f]_{{\mathsf \Lambda}^{\phi,\, p}(\rn)}^p
 &=\int_{\mathbb R^n}
   \|\Delta_h f\|_{L^p(\mathbb R^n)}^p\,
   \frac{\phi(|h|)}{|h|^{n+p}}\,dh\\
   &\ge (2D)^{-p}\int_{|h|\le D}
   \|f\|_{L^p(\mathbb R^n)}^p\,
   \frac{\phi(|h|)}{|h|^{n}}\,dh
   +2\int_{|h|> D}
   \|f\|_{L^p(\mathbb R^n)}^p\,
   \frac{\phi(|h|)}{|h|^{n+p}}\,dh \\
 &\ge \omega_{n-1} \|f\|_{L^p(\mathbb R^n)}^p
   \left(
     (2D)^{-p}
       \int_0^D\phi(t)\,\frac{dt}{t}
     +2\int_D^\infty
       \phi(t)\,\frac{dt}{t^{p+1}}
   \right).
\end{align*}
By \eqref{eq-adm}, we have
\[
     0< \frac{1}{(2D)^p}
       \int_0^D\phi(t)\,\frac{dt}{t}
     +2\int_D^\infty
       \phi(t)\,\frac{dt}{t^{p+1}}< \infty.
\]
Therefore, \eqref{general:Poincare} follows.
\end{proof}

Let $p\in(1,\infty)$ and $\phi$ satisfy \eqref{eq-adm}.
Associated with $\mathcal F^{\phi,\,p}(\Omega)$, we use the \textit{zero-exterior space}
\begin{equation*}
 \mathcal F_0^{\phi,\,p}(\Omega)
 :=\{v\in  \mathcal F^{\phi,\,p}(\Omega):\
                        v=0\ \text{a.e. on }\Omega^c\}.
\end{equation*}
Obviously, for any $v\in  \mathcal F_0^{\phi,\,p}(\Omega)$,
by \eqref{eq-ker-phi}, one has
\begin{align*}
\mathcal E_{\phi,\, p}(v;\,\Omega)\approx \left(\left[v\right]_{{\mathsf \Lambda}^{\phi,\,p}(\mathbb R^n)}\right)^p,
    \end{align*}
    with implicit constants depending only on $\Lambda$.
Consequently,
Lemma~\ref{general:nonlocal-Poincare} yields that for any $v\in \mathcal F_0^{\phi,\, p}(\Omega)$,
\begin{equation}\label{existence:coercivity}
\mathcal E_{\phi,\, p}(v;\,\Omega)
\ge c(n,p,\Omega,\phi,\Lambda)
                 \|v\|_{L^p(\rn)}^p.
\end{equation}

For fractional kernels, Di Castro, Kuusi, and Palatucci
\cite[Theorem~2.3]{DKP2016} proved that the unique energy minimizer
is precisely the weak solution of the exterior Dirichlet problem. We next show that the same conclusion holds for every admissible kernel $\phi$ satisfying \eqref{eq-adm}.

\begin{definition}
For any given $g \in  \mathcal F^{\phi,\,p}(\Omega)$,
define the \emph{admissible class}
\[
 \mathcal F_g^{\phi,\, p}(\Omega)
 =\{g+v:\ v\in \mathcal F_0^{\phi,\, p}(\Omega)\}.
\]
A function $u\in  \mathcal F_g^{\phi,\, p}(\Omega)$ is called a
\emph{weak solution} to the exterior Dirichlet problem
\begin{equation}\label{eq-exter-DP}
 \begin{cases}
 \mathcal L_{\phi,\, p}u=0&\text{ in }\,\Omega,\\
 u=g&\text{a.e. on }\,\mathbb R^n\setminus\Omega,
 \end{cases}
\end{equation}
if for every  $\zeta\in \mathcal F_0^{\phi,\, p}(\Omega)$,
\begin{align*}
 \iint_{\mathcal C_\Omega}
 |u(x)-u(y)|^{p-2}(u(x)-u(y))
 (\zeta(x)-\zeta(y))K_{\phi,\,p}(x,y)\,dy\,dx=0.
\end{align*}
\end{definition}

\begin{theorem}\label{thm-1}
Let $g \in  \mathcal F^{\phi,\,p}(\Omega).$
 A function $u\in\mathcal F_g^{\phi,\,p}(\Omega)$ is a weak solution
to \eqref{eq-exter-DP} if and only if it is a minimizer
of the energy \(\mathcal E_{\phi,\, p}(\cdot;\,\Omega)\), which means that
\[  \mathcal E_{\phi,\, p}(u;\,\Omega)\le  \mathcal E_{\phi,\, p}(\omega;\,\Omega)
\qquad\text{for all}\ \,\omega \in \mathcal F_g^{\phi,\, p}(\Omega).
\]
\end{theorem}

\begin{proof}
Suppose \(u\) is a minimizer of the energy \(\mathcal E_{\phi,\, p}(\cdot;\,\Omega)\)
in \(\mathcal F_g^{\phi,\, p}(\Omega)\). Then for any \(\varphi \in \mathcal F_0^{\phi,\, p}(\Omega)\),
\[
\frac{d}{dt} \mathcal E_{\phi,\, p}(u + t\varphi;\, \Omega) \bigg|_{t=0} =0.
\]
Moreover, by the definition of \(\mathcal E_{\phi,\, p}(\cdot;\,\Omega)\) in \eqref{corresponding energy},
\begin{align*}
&\frac{d}{dt} \mathcal E_{\phi,\, p}(u + t\varphi;\, \Omega) \bigg|_{t=0} \\
&\quad= \iint_{\mathcal C_{\Omega}}  K_{\phi,\,p}(x,y) \left.\left(\frac{d}{dt}  |u(x)-u(y)+t(\varphi(x)-\varphi(y))|^p\right) \right|_{t=0} \, dx \, dy\\
&\quad= p\iint_{\mathcal C_{\Omega}} |u(x)-u(y)|^{p-2}(u(x)-u(y))(\varphi(x)-\varphi(y)) K_{\phi,\,p}(x,y) \, dx \, dy.
\end{align*}
Hence, \(u\) is a weak solution to \eqref{eq-exter-DP}.

Conversely, let $u \in \mathcal F_g^{\phi,\, p}(\Omega)$ be a weak solution  to \eqref{eq-exter-DP}. Given any $v \in \mathcal F_g^{\phi,\, p}(\Omega)$, note that the Young inequality implies
\begin{align*}
|u(x)-u(y)|^{p-1}|v(x)-v(y)|\le \frac{1}{p} |v(x)-v(y)|^p
+ \frac{p-1}{p}|u(x)-u(y)|^p
\end{align*}
and, hence,
\begin{align*}
&\iint_{\mathcal C_{\Omega}} |u(x)-u(y)|^{p-2}(u(x)-u(y))(v(x)-v(y))K_{\phi,\,p}(x,y) \, dx \, dy\\
&\le\frac{1}{p} \iint_{\mathcal C_{\Omega}} |v(x)-v(y)|^pK_{\phi,\,p}(x,y) \, dx \, dy+\frac{p-1}{p} \iint_{\mathcal C_{\Omega}}|u(x)-u(y)|^pK_{\phi,\,p}(x,y)\, dx\, dy.
\end{align*}
Thus, for any \(v \in \mathcal F_g^{\phi,\, p}(\Omega)\), since \(\varphi := u - v \in \mathcal F_0^{\phi,\, p}(\Omega)\), we have
\begin{align*}
0 &= \iint_{\mathcal C_{\Omega}}|u(x)-u(y)|^{p-2}(u(x)-u(y))(\varphi(x)-\varphi(y)) K_{\phi,\,p}(x,y) \, dx \, dy&\\
&= \iint_{\mathcal C_{\Omega}} |u(x)-u(y)|^p K_{\phi,\,p}(x,y)\, dx \, dy&\\
& \quad - \iint_{\mathcal C_{\Omega}} |u(x)-u(y)|^{p-2}(u(x)-u(y))(v(x)-v(y))K_{\phi,\,p}(x,y) \, dx \, dy&\\
&\geq \frac{1}{p} \iint_{\mathcal C_{\Omega}}|u(x)-u(y)|^pK_{\phi,\,p}(x,y)\, dx\, dy - \frac{1}{p} \iint_{\mathcal C_{\Omega}} |v(x)-v(y)|^pK_{\phi,\,p}(x,y) \, dx \, dy.
\end{align*}
This proves that $u$ is a minimizer of $ \mathcal E_{\phi,\, p}(\cdot;\,\Omega)$ in \(\mathcal F_g^{\phi,\, p}(\Omega)\).
\end{proof}

Having established that every minimizer of the energy is a weak solution to the exterior Dirichlet problem \eqref{eq-exter-DP}, it is natural to ask whether such a minimizer exists and is unique. The following theorem answers this question affirmatively.

\begin{theorem}\label{existence}
Let \(g \in \mathcal F^{\phi,\,p}(\Omega)\). Then there exists a unique minimizer \(u \in \mathcal F_g^{\phi,\, p}(\Omega)\) of the energy \(\mathcal E_{\phi,\, p}(\cdot;\,\Omega)\). Consequently, the exterior Dirichlet problem \eqref{eq-exter-DP} admits a unique weak solution.
\end{theorem}

\begin{proof}
According to Theorem~\ref{thm-1}, it suffices to establish the existence and uniqueness of the minimizer. We first prove its existence by the direct method of the calculus of variations.
Let
\[
d\mu(x,y):=
\frac{\phi(|x-y|)}{|x-y|^{n+p}}\,dx\,dy.
\]
Consider the linear map
\[
\begin{aligned}
T:\ {\mathsf \Lambda}^{\phi,\, p}(\mathbb R^n)
&\to   L^p(\mathbb R^n)\times L^p(\mathbb R^{2n},\,\mu),\\
v&\longmapsto (v,\,v(x)-v(y)).
\end{aligned}
\]
By the definition of ${\mathsf \Lambda}^{\phi,\, p}(\mathbb R^n)$,
$T$ is an isometric embedding.

Next, we show that the image of $T$ is
closed. Indeed, suppose that
\[
Tv_j\to  (f,h)
\]
strongly in the product space \(L^p(\mathbb R^n)\times L^p(\mathbb R^{2n},\,\mu)\) as \(j\to\infty\), which means that both
\[
v_j\to   f\quad\text{in }\, L^p(\mathbb R^n),
\]
and
\[
v_j(x)-v_j(y)\to   h(x,y)
\quad\text{in }\, L^p(\mathbb R^{2n},\,\mu).
\]
By the almost everywhere convergence theorem for $L^p$
convergence, after passing to a common subsequence, we may
assume that
\[
\begin{cases}
v_j(x)\to   f(x),
&\text{for a.e. }x\in\mathbb R^n
\text{ with respect to Lebesgue measure},\\[3pt]
v_j(x)-v_j(y)\to   h(x,y),
&\text{for }\mu\text{-a.e. }(x,y)\in\mathbb R^{2n}.
\end{cases}
\]
Since $\mu$ is absolutely continuous with respect to
Lebesgue measure on $\mathbb R^{2n}$, the first convergence
also implies
\[
v_j(x)-v_j(y)\to   f(x)-f(y)
\quad\text{for }\mu\text{-a.e. }(x,y)\in\mathbb R^{2n}.
\]
By the uniqueness of almost everywhere limits, we obtain
\[
h(x,y)=f(x)-f(y)
\quad\text{for }\mu\text{-a.e. }(x,y)\in\mathbb R^{2n}.
\]
Consequently, $(f,h)=Tf$, and hence the image of $T$
is closed.

Since $1<p<\infty$, the product space
$L^p(\mathbb R^n)\times
L^p(\mathbb R^{2n},\,\mu)$
is reflexive. Any closed subspace of a reflexive Banach space is again reflexive. Consequently, the image \(T({\mathsf \Lambda}^{\phi,\, p}(\mathbb R^n))\) is reflexive, and hence \({\mathsf \Lambda}^{\phi,\, p}(\mathbb R^n)\) is reflexive.
Further, since $\mathcal F_0^{\phi,\, p}(\Omega)$ is a closed subspace
of ${\mathsf \Lambda}^{\phi,\, p}(\mathbb R^n)$, it is also reflexive.

Now, let
$u_j=g+v_j\in \mathcal F_g^{\phi,\, p}(\Omega)$
be a minimizing sequence, namely,
\[
\lim_{j\to\infty}\mathcal E_{\phi,\, p}(u_j;\, \Omega)
=
\inf_{u\in \mathcal F_g^{\phi,\, p}(\Omega)}
\mathcal E_{\phi,\, p}(u;\, \Omega).
\]
Since $g\in\mathcal F^{\phi,\, p}(\Omega)$,
it follows that
$\mathcal E_{\phi,\, p}(g;\,\Omega)<\infty$. By the elementary
$p$-growth inequality,
\[
\begin{aligned}
\mathcal E_{\phi,\, p}(v_j;\,\Omega)
&\le 2^{p-1}\big(
\mathcal E_{\phi,\, p}(u_j;\,\Omega)
+\mathcal E_{\phi,\, p}(g;\,\Omega)\big)\\
&\le C,
\end{aligned}
\]
where $C$ is a positive constant independent of $j$.
On the one hand, since $v_j=0$ a.e. on $\Omega^c$,
the kernel comparison \eqref{eq-ker-phi} gives
\begin{align*}
[v_j]_{{\mathsf \Lambda}^{\phi,\, p}(\mathbb R^n)}^p
&=\iint_{ \mathcal C_\Omega}
\frac{|v_j(x)-v_j(y)|^p\,\phi(|x-y|)}
{|x-y|^{n+p}}\,dy\,dx\\
&\le \Lambda\,\mathcal E_{\phi,\, p}(v_j;\, \Omega)
\le \Lambda\,C.
\end{align*}
On the other hand, it follows from  $v_j\in \mathcal F_0^{\phi,\, p}(\Omega)$ and
\eqref{existence:coercivity} that
\[
\|v_j\|_{L^p(\rn)}^p
\le C(n,p,\Omega,\phi,\Lambda)\,[v_j]_{{\mathsf \Lambda}^{\phi,\, p}(\mathbb R^n)}^p
\le C.
\]
Combining the above two estimates, we conclude that
\[
\sup_{j\in\mathbb N}\|v_j\|_{{\mathsf \Lambda}^{\phi,\, p}(\mathbb R^n)}<\infty.
\]
By reflexivity, after passing to a subsequence, there exists
$v\in\mathcal F_0^{\phi,\, p}(\Omega)$ such that
\[
v_j\rightharpoonup v
\quad\text{in }{\mathsf \Lambda}^{\phi,\, p}(\mathbb R^n).
\]
Here we have used the fact that a closed convex subset
of a Banach space is weakly closed. In particular,
$u:=g+v$ belongs to the admissible class $\mathcal F_g^{\phi,\, p}(\Omega)$.

We still need to verify that $u$ is a minimizer.
For any measurable function $f$ on $\mathbb R^n$, define
\[
\mathcal D_\phi f(x,y)
:=(f(x)-f(y))K_{\phi,\,p}(x,y)^{1/p},
\qquad (x,y)\in\mathcal C_\Omega.
\]
By the kernel comparison \eqref{eq-ker-phi},
for every $f\in\mathsf\Lambda^{\phi,\, p}(\mathbb R^n)$,
we have
\begin{align*}
\|\mathcal D_\phi f\|_{L^p(\mathcal C_\Omega)}^p
=\mathcal E_{\phi,\, p}(f;\, \Omega)\le \Lambda
[f]_{\mathsf\Lambda^{\phi,\, p}(\mathbb R^n)}^p.
\end{align*}
Thus $\mathcal D_\phi$ is a bounded linear operator
from $\mathsf\Lambda^{\phi,\, p}(\mathbb R^n)$
into $L^p(\mathcal C_\Omega)$.
Since bounded linear operators preserve weak convergence,
the weak convergence
$v_j\rightharpoonup v$ in
$\mathsf\Lambda^{\phi,\,p}(\mathbb R^n)$
implies that
\[
\mathcal D_\phi v_j\rightharpoonup
\mathcal D_\phi v
\quad\text{in }\,L^p(\mathcal C_\Omega).
\]
Since $u_j=g+v_j$, $u=g+v$, and
$\mathcal D_\phi g\in L^p(\mathcal C_\Omega)$,
we further obtain
\[
\mathcal D_\phi u_j\rightharpoonup
\mathcal D_\phi u
\quad\text{in }\, L^p(\mathcal C_\Omega).
\]
By the weak lower semicontinuity of the $L^p$ norm,
\begin{align*}
\mathcal E_{\phi,\, p}(u;\, \Omega)
=\|\mathcal D_\phi u\|_{L^p(\mathcal C_\Omega)}^p
\le\liminf_{j\to\infty}
\|\mathcal D_\phi u_j\|_{L^p(\mathcal C_\Omega)}^p
=\liminf_{j\to\infty}
\mathcal E_{\phi,\, p}(u_j;\, \Omega).
\end{align*}
Consequently,
\begin{align*}
\inf_{w\in\mathcal F_g^{\phi,\, p}(\Omega)}
\mathcal E_{\phi,\, p}(w;\, \Omega)
&\le \mathcal E_{\phi,\, p}(u;\, \Omega)\\
&\le \liminf_{j\to\infty}
\mathcal E_{\phi,\, p}(u_j;\, \Omega)\\
&=\inf_{w\in\mathcal F_g^{\phi,\, p}(\Omega)}
\mathcal E_{\phi,\, p}(w;\, \Omega).
\end{align*}
Hence $u$ is a minimizer.

It remains to prove uniqueness. Let $u_1$ and $u_2$
be two minimizers in $\mathcal F_g^{\phi,\,p}(\Omega)$,
and set $v=(u_1+u_2)/2$.
Since $v\in\mathcal F_g^{\phi,\,p}(\Omega)$,
by minimality and the strict convexity of
$t\mapsto |t|^p$, we have
\begin{align*}
  \mathcal E_{\phi,\, p}(u_1;\, \Omega)
  &\le \mathcal E_{\phi,\, p}(v;\, \Omega)\\
  &\le \frac12\mathcal E_{\phi,\, p}(u_1;\, \Omega)
       +\frac12\mathcal E_{\phi,\, p}(u_2;\, \Omega)\\
  &=\mathcal E_{\phi,\, p}(u_1;\, \Omega)<\infty.
\end{align*}
Thus equality holds in the integrated convexity inequality.
Since the pointwise convexity gap is nonnegative and
$t\mapsto |t|^p$ is strictly convex for $p>1$, it follows that
\[
  u_1(x)-u_1(y)=u_2(x)-u_2(y)
\]
for almost every
$(x,y)\in\mathcal C _\Omega$
such that $K_{\phi,\,p}(x,y)>0$.
Consequently,
\[
  \mathcal E_{\phi,\, p}(u_1-u_2;\, \Omega)=0.
\]
Moreover, $u_1-u_2\in \mathcal F_0^{\phi,\, p}(\Omega)$ because
$u_1,u_2\in \mathcal F_g^{\phi,\, p}(\Omega)$.
Applying \eqref{existence:coercivity}, we obtain
\[
  \|u_1-u_2\|_{L^p(\mathbb R^n)}^p
    \le C\mathcal E_{\phi,\, p}(u_1-u_2;\, \Omega)=0.
\]
Hence $u_1=u_2$ a.e. on $\mathbb R^n$.
This proves the uniqueness.
\end{proof}

\subsection{Existence and uniqueness for the nonlocal obstacle problem}
\label{ss:obstacle}

For fractional kernels, Korvenp\"a\"a, Kuusi, and Palatucci
\cite[Theorem~1]{KKP2016Obstacle} proved that the obstacle problem
has a unique solution and that this solution is a weak supersolution
in the whole domain. We next show that the same conclusion holds for every admissible kernel $\phi$ satisfying \eqref{eq-adm}.

Given $g,\psi\in\mathcal F^{\phi,\,p}(\Omega)$,
we define the admissible class
\[
\mathcal A_{\psi,\,g}
:=
\left\{
w\in\mathcal F_g^{\phi,\,p}(\Omega):\
w\ge\psi\ \,\text{a.e. on }\,\Omega
\right\}.
\]
Assume that $\mathcal A_{\psi,\,g}\ne\emptyset$. Then it is convex and
strongly closed in $\mathcal F_g^{\phi,\,p}(\Omega)$.
Indeed, if $w_j\in\mathcal A_{\psi,\,g}$ and
$w_j\to w$ strongly in $\mathcal F_g^{\phi,\,p}(\Omega)$,
then $w_j\to w$ in $L^p(\Omega)$; passing to a subsequence if necessary,
we may assume that $w_j\to w$ a.e. in $\Omega$, and hence
$w\ge\psi$ a.e. in $\Omega$.
Consequently, $\mathcal A_{\psi,\,g}$ is also weakly closed.
By the same direct method as in the proof of
Theorem~\ref{existence}, there exists a unique minimizer
$u\in\mathcal A_{\psi,\,g}$ of
$\mathcal E_{\phi,\,p}(\cdot;\Omega)$.

The following theorem shows that this minimizer is
a weak supersolution in the whole domain $\Omega$.

\begin{theorem}\label{thm:obstacle-supersolution}
Assume that $\mathcal A_{\psi,\,g}\ne\emptyset$ and let $u$ be the minimizer of
$\mathcal E_{\phi,\,p}(\cdot;\Omega)$ in
$\mathcal A_{\psi,\,g}$.
Then $u$ is a weak supersolution to
\[
\mathcal L_{\phi,\,p}u=0
\qquad\text{in }\,\Omega,
\]
that is, for every nonnegative
$\zeta\in\mathcal F_0^{\phi,\,p}(\Omega)$,
\begin{equation*}
\iint_{\mathcal C_\Omega}
|u(x)-u(y)|^{p-2}(u(x)-u(y))
(\zeta(x)-\zeta(y))K_{\phi,\,p}(x,y)\,dy\,dx
\ge0.
\end{equation*}
\end{theorem}

\begin{proof}
Let $\zeta\in\mathcal F_0^{\phi,\,p}(\Omega)$ be nonnegative.
For every $\tau>0$, we have
\[
u+\tau\zeta\in\mathcal A_{\psi,\,g},
\]
since $u+\tau\zeta\ge u\ge\psi$ a.e. on $ \Omega$ and
$u+\tau\zeta=g$ a.e. on $\Omega^c$.
By the minimality of $u$,
\[
\mathcal E_{\phi,\,p}(u+\tau\zeta;\Omega)
\ge
\mathcal E_{\phi,\,p}(u;\Omega).
\]
Consequently,
\begin{equation}\label{eq:obstacle-difference-quotient}
\begin{aligned}
\iint_{\mathcal C_\Omega}
\frac{
|u(x)-u(y)+\tau(\zeta(x)-\zeta(y))|^p
-|u(x)-u(y)|^p
}{p\tau}
K_{\phi,\,p}(x,y)\,dy\,dx\ge 0.
\end{aligned}
\end{equation}

We now justify passing to the limit as $\tau\downarrow0$.
For $0<\tau\le1$, the mean value theorem gives
\[
\begin{aligned}
\left|
|a+\tau b|^p-|a|^p
\right|\le
p\tau(|a|+|b|)^{p-1}|b|
\le C_p\,p\tau\big(|a|^{p-1}|b|+|b|^p\big)
\end{aligned}
\]
for all $a,b\in\mathbb R$. Taking
\[
a=u(x)-u(y),\qquad
b=\zeta(x)-\zeta(y),
\]
we find that the integrand in \eqref{eq:obstacle-difference-quotient}
is bounded, uniformly in $\tau\in(0,1]$, by
$$
C_p
\left(
|u(x)-u(y)|^{p-1}|\zeta(x)-\zeta(y)|
+|\zeta(x)-\zeta(y)|^p
\right)K_{\phi,\,p}(x,y).
$$
Moreover, applying
the H\"older inequality yields
\[
\begin{aligned}
&\iint_{\mathcal C_\Omega}
\big(|u(x)-u(y)|^{p-1}|\zeta(x)-\zeta(y)|
+|\zeta(x)-\zeta(y)|^p\big)
K_{\phi,\,p}(x,y)\,dy\,dx\\
&\quad\le
\mathcal E_{\phi,\,p}(u;\Omega)^{\frac{p-1}{p}}
\mathcal E_{\phi,\,p}(\zeta;\Omega)^{\frac1p}
+\mathcal E_{\phi,\,p}(\zeta;\Omega)
<\infty.
\end{aligned}
\]
Thus, the dominated
convergence theorem allows us to let $\tau\downarrow0$ in
\eqref{eq:obstacle-difference-quotient}, which yields
\[
\iint_{\mathcal C_\Omega}
|u(x)-u(y)|^{p-2}(u(x)-u(y))
(\zeta(x)-\zeta(y))
K_{\phi,\,p}(x,y)\,dy\,dx\ge 0.
\]
Therefore, $u$ is a weak supersolution in $\Omega$.
\end{proof}

\section{Kernel and Energy Estimates }
\label{sec:logsp}


\subsection{Basic estimates for the nearly critical kernel }
\label{sec:kernel}

For the function $\phi\colon(0,\infty)\to[0,\infty)$ satisfying \eqref{eq-adm}, let
$\Phi$ be the associated function given in \eqref{eq-Phi}. We split
$\Phi$ as
\begin{align}\label{eq-Phi0infty}
\Phi(r)
&=\int_0^1
\phi(r\tau)\,\frac{d\tau}{\tau}+\int_1^\infty
\phi(r\tau)\tau^{-p}\,\frac{d\tau}{\tau}=:\Phi_0(r)+\Phi_{\infty}(r).
\end{align}
The part $\Phi_\infty$ is closely connected to the tail estimate for the
kernel $K_{\phi,\,p}$ given in \eqref{eq-ker-phi}.

\begin{lemma}\label{1.3}
For any $r\in(0,\infty)$ and $x\in\mathbb R^n$,
\begin{equation}\label{eq-est-tail}
  \int_{\mathbb R^n\setminus B_r(x)}
     \frac{\phi(|x-y|)}{|x-y|^{n+p}}\,dy=\omega_{n-1} \frac{\Phi_\infty(r)}{r^p}.
\end{equation}
Consequently,
\begin{equation}\label{eq-est-K-tail}
 \Lambda^{-1} \omega_{n-1} r^{-p}\Phi_\infty(r) \le \int_{\mathbb R^n\setminus B_r(x)}
     K_{\phi,\,p}(x,y)\,dy\le \Lambda \omega_{n-1} \frac{\Phi_\infty(r)}{r^p}.
\end{equation}
\end{lemma}

\begin{proof}
Using polar coordinates, we directly obtain
\[
\int_{\mathbb R^n\setminus B_r(x)}
     \frac{\phi(|x-y|)}{|x-y|^{n+p}}\,dy
     =\omega_{n-1}\int_r^\infty \phi(t)\,t^{-p-1}\,dt=\omega_{n-1} r^{-p}\Phi_\infty(r),
\]
as desired.
\end{proof}

The following lemma quantifies the global cutoff behavior of the kernel  \(K_{\phi,\,p}\) that satisfies \eqref{eq-ker-phi}, showing that the truncated integral at scale
$a$
is uniformly comparable to \(\Phi(a)/a^p\).

\begin{lemma}\label{kernel:global-cutoff}
For any $a>0$,
\begin{equation}\label{kernel:cutoff-bound}
\sup_{x\in\mathbb R^n}\int_{\mathbb R^n}
\min\left\{1,\,\frac{|x-y|^p}{a^p}\right\}
\frac{\phi(|x-y|)}{|x-y|^{n+p}}\,dy
=\omega_{n-1} \frac{\Phi(a)}{a^p}.
\end{equation}
In particular,
\begin{equation}\label{kernel:cutoff-bound-K}
\Lambda^{-1} \omega_{n-1} \frac{\Phi(a)}{a^p}
\le \sup_{x\in\mathbb R^n}\int_{\mathbb R^n}
\min\left\{1,\,\frac{|x-y|^p}{a^p}\right\}
K_{\phi,\,p}(x,y)\,dy
\le \Lambda \omega_{n-1} \frac{\Phi(a)}{a^p}.
\end{equation}
\end{lemma}

\begin{proof}
By the kernel comparison \eqref{eq-ker-phi},
the two-sided estimate \eqref{kernel:cutoff-bound-K} follows directly
from \eqref{kernel:cutoff-bound}.
Further, for any \(x\in\rn\), using polar coordinates, we have
\begin{align*}
&\int_{\mathbb R^n}
 \min\left\{1,\,\frac{|x-y|^p}{a^p}\right\}
 \frac{\phi(|x-y|)}{|x-y|^{n+p}}\,dy\\
     &\quad =\omega_{n-1}\int_0^\infty \min\left\{1,\frac{t^p}{a^p}\right\} \phi(t)\,t^{-p-1}\,dt\notag\\
     &\quad=\omega_{n-1}a^{-p}\int_0^\infty \min\left\{1,t^{-p}\right\} \phi(at)\,\frac{dt}{t}\\ \notag
     &\quad=\omega_{n-1}a^{-p}\Phi(a), \notag
\end{align*}
which proves  \eqref{kernel:cutoff-bound}.
\end{proof}

\begin{remark}
 If a Lipschitz function $\eta$ satisfies $0\le\eta\le1$ and $|\nabla \eta|\le C_\eta/a$ with $C_\eta\ge1$, then
\[
  |\eta(x)-\eta(y)|
  \le C_\eta\min\left\{1,\ \frac{|x-y|}{a}\right\}.
\]
Lemma~\ref{kernel:global-cutoff} consequently gives
\begin{equation}\label{LH:global-cutoff}
  \sup_{x\in \rn}\int_{\mathbb R^n}|\eta(x)-\eta(y)|^p
                            K_{\phi,\,p}(x,y)\,dy
  \le C_\eta^p\,\Lambda\,\omega_{n-1}\, \frac{\Phi(a)}{a^p}.
\end{equation}

\end{remark}

We end this subsection with the following basic properties of
$\Phi$, $\Phi_0$, and $\Phi_\infty$ for  nearly critical kernels.

\begin{lemma}\label{lem-basic-nearly critical phi}
Let $\phi$ be a nearly critical kernel. Then, the following properties hold:
\begin{enumerate}[\rm (i)]

 \item There exist constants $0<c<C<\infty$ such that for any $r \in(0,\infty)$,
\begin{equation}\label{eq-Phiinfty-phi-Phi0}
c\Phi_\infty(r)\le \phi(r) \le C \Phi_0(r)
\end{equation}
and
\begin{equation}\label{eq-Phi0-Phi}
\Phi_0(r)\le \Phi(r)\le C \Phi_0(r).
\end{equation}

\item For any $\varepsilon\in(0, 1)$,  there exist constants $c,\,C_\varepsilon>0$ such that for $0<a\le b<\infty$,
\begin{equation}\label{eq-Phi-power-comparison}
c\le\frac{\Phi(b)}{\Phi(a)}
 \le C_\varepsilon\left(\frac ba\right)^\varepsilon.
\end{equation}
\end{enumerate}
\end{lemma}

\begin{proof}
For (i), note that \eqref{eq-Phi0-Phi} follows directly from
\eqref{eq-Phi0infty} and  \eqref{eq-Phiinfty-phi-Phi0}.
To show \eqref{eq-Phiinfty-phi-Phi0}, we fix some $\varepsilon\in(0,1)$.
For any $r\in(0,\infty)$, applying the mild upper
scaling condition~\eqref{condition2} yields
\[
\phi(r)\le C_\varepsilon
 \tau^{-\varepsilon}\phi(r\tau)
\qquad\text{for all }\, \tau\in(0,1]
\]
and
\[
\phi(\tau r)\le C_\varepsilon
\tau^\varepsilon\phi(r)
\qquad\text{for all }\,  \tau\in [1,\infty).
\]
Consequently, we obtain
\begin{align*}
\Phi_0(r)=\int_0^1
\phi(r\tau)\,\frac{d\tau}{\tau}
&\ge C_{\varepsilon}^{-1}\phi(r)
\int_0^1 \tau^{\varepsilon-1}\,d\tau=(\varepsilon\,C_{\varepsilon})^{-1}\phi(r)
\end{align*}
and
\begin{align*}
\Phi_\infty(r)
&=\int_1^\infty
\phi(r\tau)\tau^{-p}\,\frac{d\tau}{\tau}
\le  C_{\varepsilon}\phi(r)
\int_1^\infty t^{\varepsilon-p-1}\,dt
=\frac{ C_{\varepsilon}}{p-\varepsilon}\phi(r).
\end{align*}
Combining the last two estimates, we obtain \eqref{eq-Phiinfty-phi-Phi0}.

We now prove (ii). Fix $\varepsilon\in(0,1)$. By \eqref{condition2}, for all
$\tau\in(0,\infty)$ and $0<a<b<\infty$, we have
\[
\phi(b\tau)\le C_\varepsilon \left(\frac b a\right)^\varepsilon \phi(a\tau).
\]
This, together with the expression of $\Phi$ in \eqref{eq-Phi}, gives
\begin{align*}
\Phi(b)=\int_0^\infty
\phi(b\tau)\min\{1,\,\tau^{-p}\}\,\frac{d\tau}{\tau}
\le C_\varepsilon \left(\frac b a\right)^\varepsilon \int_0^\infty
\phi(a\tau)\min\{1,\,\tau^{-p}\}\,\frac{d\tau}{\tau}
= C_\varepsilon \left(\frac b a\right)^\varepsilon  \Phi(a),
\end{align*}
which proves the second inequality in \eqref{eq-Phi-power-comparison}.
Similarly, the quasi-increasing property of $\phi$ yields the first
inequality in \eqref{eq-Phi-power-comparison}.
\end{proof}

\begin{remark}\label{rem-kerK}
Let $x_0\in\rn$ and $r\in(0,\infty)$.
If $|z-x_0|\le r$ and $|y-x_0|\ge 2r$, then
\[
  \frac12|y-x_0|\le |y-z|\le 2|y-x_0|.
\]
For a nearly critical kernel $\phi$, the mild upper
scaling condition~\eqref{condition2} and quasi-increasing property  therefore imply
\begin{align}\label{eq-KKer-est}
  K_{\phi,\, p}(z,y)
  &\approx\frac{\phi(|y-z|)}{|y-z|^{n+p}}\approx\frac{\phi(|y-x_0|)}{|y-x_0|^{n+p}}
   \approx K_{\phi,\, p}(x_0,y).
   \end{align}
All constants in these comparisons are  independent of
$y,z, x_0$ and $r$.
\end{remark}

\subsection{Localization and stability of $\Lambda^{\phi,\,p}(\rn)$}\label{sec-local-stab}


This subsection is devoted to establishing the functional-analytic properties of the energy space, which  ensure that the tails are finite and that the localized test functions used below have finite energy.

The following lemma shows that functions with finite local energy possess finite tails.

\begin{lemma}\label{energy:tests}
Let \(p\in(1,\infty)\) and \(\phi\) be a nearly critical kernel.
Assume  that \(u:\,\mathbb R^n\to\mathbb R\) satisfies \(u\in L^p(\Omega)\) and
\(\mathcal E_{\phi,\,p}(u;\,\Omega)<\infty.\)
Then \(u\in L^p_{\mathrm{loc}}(\mathbb R^n)\). Moreover,
\begin{align}\label{eq-int-u}
\int_{\mathbb R^n}(1+|y|)^{-(n+p)}|u(y)|^p\,dy<\infty.
\end{align}
In particular, \(u\in\mathbb L^{\phi}_{p}(\mathbb R^n)\) and \(\mathsf{Tail}^\phi_p (u;x_0,R)<\infty\) for all \(x_0\in\rn\) and \(R>0\).
\end{lemma}

\begin{proof} The whole proof is divided into three steps.

\medskip

\textit{Step 1:\, Proof of $u\in L^p_{\mathrm{loc}}(\mathbb R^n)$.}
Since \(\Omega\Subset\rn\) is open, we may choose a small ball \(B_{2r}(x_0)\subset\Omega\), where \(r\in(0,\infty)\) and \(x_0\in\rn\). Due to the precompactness of \(\Omega\), we take an
arbitrarily large \(R\) satisfying
\[
R>R_0:=2(\operatorname{diam}(\Omega)+|x_0|+1).
\]
To prove that \(u\in L^p_{\mathrm{loc}}(\mathbb R^n)\), it suffices to show that \(u\in L^p(B_R(x_0))\). Write
$$
\int_{B_R(x_0)} |u(y)|^p\,dy
\le \int_{\Omega} |u(y)|^p\,dy
+ \int_{B_R(x_0)\setminus \Omega} |u(y)|^p\,dy.
$$
The first term on the right-hand side is finite since \(u\in L^p(\Omega)\).
For the second term, we write
\begin{align}\label{eq-u-exteri}
&\int_{B_R(x_0)\setminus \Omega} |u(y)|^p\,dy\\
&\quad= \fint_{B_r(x_0)} \int_{B_R(x_0)\setminus \Omega} |u(y)|^p\,dy\,dx\notag\\
&\quad\le 2^{p-1} \fint_{B_r(x_0)} \int_{B_R(x_0)\setminus \Omega}
\left(|u(y)-u(x)|^p+|u(x)|^p\right)\,dy\,dx\notag\\
&\quad\le 2^{p-1} \fint_{B_r(x_0)} \int_{B_R(x_0)\setminus \Omega}
|u(y)-u(x)|^p\,dy\,dx+2^{p-1}\frac{R^n}{r^n}
\int_{B_r(x_0)}|u(x)|^p\,dx.\notag
\end{align}
Again, by $u\in L^p(\Omega)$, the second integral on the right-hand side of \eqref{eq-u-exteri} is finite.
So, the proof of \(u\in L^p(B_R(x_0))\) reduces to showing that
\begin{align}\label{eq-u-ux-uy}
\int_{B_r(x_0)} \int_{B_R(x_0)\setminus \Omega}
|u(y)-u(x)|^p\,dy\,dx<\infty.
\end{align}
Note that, for any \(y\in B_R(x_0)\setminus \Omega\) and \(x\in B_r(x_0)\), we have \(y\notin B_{2r}(x_0)\) and hence \(r<|x-y|<2R\), so the quasi-increasing property of $\phi$ gives
$$
\frac{\phi(|x-y|)}{|x-y|^{n+p}}
 \ge c (2R)^{-n-p} \phi(r)=:c_0.
$$
From this and $B_{r}(x_0)\subset\Omega$, it follows that
\begin{align*}
&c_0\int_{B_r(x_0)} \int_{B_R(x_0)\setminus \Omega}
|u(y)-u(x)|^p\,dy\,dx\\
&\quad
\le \int_{B_r(x_0)} \int_{B_R(x_0)\setminus \Omega}
\frac{|u(y)-u(x)|^p}{|x-y|^{n+p}}\,\phi(|x-y|)\,dy\,dx\\
&\quad
 \le\Lambda\mathcal E_{\phi,\,p}(u;\,\Omega)<\infty.
\end{align*}
This proves \eqref{eq-u-ux-uy}.
Thus, we obtain \(u\in L^p(B_R(x_0))\) for all \(R>R_0\). In particular, \(u\in L_{\mathrm{loc}}^p(\rn)\).

\medskip

\textit{Step 2:\, Proof of \eqref{eq-int-u}.}
Let $R_0$ be as in \textit{Step 1}.
Since we have proved \(u\in L_{\mathrm{loc}}^p(\rn)\), it suffices to prove that for \(R>R_0\),
\begin{align}\label{eq2-u-exteri}
 \int_{\mathbb R^n\setminus B_R(x_0)} (1+|y|)^{-(n+p)}|u(y)|^p\,dy<\infty.
\end{align}
To this end, following the same lines as in \eqref{eq-u-exteri}, we obtain
\begin{align*}
&\int_{\mathbb R^n\setminus B_R(x_0)} (1+|y|)^{-(n+p)}|u(y)|^p\,dy\\
&\quad\le 2^{p-1} \fint_{B_r(x_0)} \int_{\mathbb R^n\setminus B_R(x_0)}
(1+|y|)^{-(n+p)} |u(y)-u(x)|^p\,dy\,dx\\
&\qquad+2^{p-1}\fint_{B_r(x_0)}
\int_{\mathbb R^n\setminus B_R(x_0)} (1+|y|)^{-(n+p)} |u(x)|^p\,dy\,dx\\
&\quad=:{\rm I}+{\rm II}.\notag
\end{align*}
For $x\in B_r(x_0)$ and $y\notin B_R(x_0)$, we have $|x-y|\ge R-r$ and, hence,
by the quasi-increasing property of $\phi$,
$$C\phi(|x-y|)\ge \phi(R-r).$$ Moreover,
the triangle inequality gives $$|x-y|\le |y|+|x_0|+|x-x_0|\le |y|+|x_0|+r\le(1+|x_0|+r)(1+|y|).$$
Consequently,
\begin{align*}
{\rm I}
&\le C \fint_{B_r(x_0)} \int_{\mathbb R^n\setminus B_R(x_0)}
|y-x|^{-(n+p)} |u(y)-u(x)|^p \frac{\phi(|x-y|)}{\phi(R-r)}\,dy\,dx\\
&\le \frac{C}{r^n\phi(R-r)}
\int_{\Omega} \int_{\mathbb R^n\setminus \Omega}
|y-x|^{-(n+p)} |u(y)-u(x)|^p \phi(|x-y|)\,dy\,dx\\
&\le  \frac{C\,\Lambda }{r^n\phi(R-r)} \mathcal E_{\phi,\,p}(u;\,\Omega)<\infty.
\end{align*}
For ${\rm II}$, it is easy to see that
$$
{\rm II}\le \frac{2^{p-1}}{|B_r(x_0)|}\|u\|_{L^p(\Omega)}^p
\left(\int_{\rn} (1+|y|)^{-n-p}\,dy\right)<\infty.
$$
Combining the estimates of  ${\rm I}$ and  ${\rm II}$ yields the finiteness in \eqref{eq2-u-exteri}. This proves \eqref{eq-int-u}.

\medskip

\textit{Step 3:\, Proof of \(u\in\mathbb L^{\phi}_p(\mathbb R^n)\) and \(\mathsf{Tail}^\phi_p (u;x_0,R)<\infty\).} Fix some $\varepsilon\in(0,1)$.
Then, by the quasi-increasing property of $\phi$ and the mild upper scaling condition \eqref{condition2}, we obtain
\begin{align}\label{eq-phit}
\phi(t)\le C_\varepsilon\,(1+t)^\varepsilon
\end{align}for all \(t>0\) and for some $C_\varepsilon>0$.  Consequently,
$$\int_{\mathbb R^n}
     \frac{|u(y)|^{p-1}}{(1+|y|)^{n+p}}\phi(|y|)\,dy
     \le C_\varepsilon \int_{\mathbb R^n}
     \frac{|u(y)|^{p-1}}{(1+|y|)^{n+p-\varepsilon}}\,dy.$$
From the H\"older inequality and \eqref{eq-int-u}, it follows that
\begin{align}\label{eq-p-1-poweru}
  \int_{\mathbb R^n}
     \frac{|u(y)|^{p-1}}{(1+|y|)^{n+p-\varepsilon}}\,dy
  &\le
   \left(\int_{\mathbb R^n}\frac{|u(y)|^p}{(1+|y|)^{n+p}}\,dy
   \right)^{\frac{p-1}{p}}
   \left(\int_{\mathbb R^n}\frac1{(1+|y|)^{n+p(1-\varepsilon)}}\,dy\right)^{\frac{1}{p}}<\infty.
\end{align}
Combining the last two estimates yields \(u\in \mathbb L^{\phi}_p(\mathbb R^n)\).

Now, fix \(a\in\rn\) and \(s>0\).
It follows from \eqref{Def-tail} and \eqref{eq-phit} that
\begin{align*}
\left(\mathsf{Tail}^\phi_p (u;a,s)\right)^{p-1}
 \le  C_{n,\phi,s,p,\varepsilon}\int_{\rn\setminus B_s(a)}
\frac{|u(y)|^{p-1}}{|a-y|^{n+p-\varepsilon}}\,dy.
\end{align*}
The value \(|a-y|\) in the denominator can be replaced by \(1+|y|\), up to a constant depending on \(a\) and \(s\). To see this,
we observe that when \(|y-a|\ge s\), the triangle inequality gives
\[
1+|y|
\le 1+|y-a|+|a|\le C_{s,a} (1+|y-a|).
\]
By this and \eqref{eq-p-1-poweru}, we conclude that
\begin{align*}
\left(\mathsf{Tail}^\phi_p (u;a,s)\right)^{p-1}
\le  C_{n,\phi,s,p,\varepsilon,a}\int_{\rn\setminus B_s(a)}
\frac{|u(y)|^{p-1}}{(1+|y|)^{n+p-\varepsilon}}\,dy<\infty.
\end{align*}
This completes the proof of Lemma \ref{energy:tests}.
\end{proof}

The next lemma shows that the energy class $\mathcal F^{\phi,\,p}(\Omega)$ is stable under localized
Lipschitz truncation or composition. Typical choices of Lipschitz
composition include  $F(t)=(t_++d)^{-b}$ for $b,d>0$ and $F(t)=t\wedge k$; see
Proposition~\ref{LH:linear-truncation} below.

\begin{lemma}\label{lem-stability}
Let \(p\in(1,\infty)\) and \(\phi\) satisfy \eqref{eq-adm}.
 If $u\in\mathcal F^{\phi,\,p}(\Omega)$, \(\eta:\,\rn\to\mathbb R\) is Lipschitz, \(\operatorname{supp}\eta\subset\Omega\), and \(F:\,\mathbb R\to\mathbb R\) is Lipschitz, then \(\eta F(u)\in \Lambda^{\phi,\,p}(\rn)\).
\end{lemma}
\begin{proof}
Since \(F:\,\mathbb R\to\mathbb R\) is Lipschitz, there exists a constant \(L>0\) such that for all \(s,t\in\mathbb R\),
\[
|F(t)-F(s)|\le L|t-s|.
\]
Taking \(s=0\), this gives
\[
|F(t)|\le |F(0)|+L|t|.
\]
With \(t=u(x)\), we obtain
\begin{align}\label{eq-Fu-bdd}
|F(u(x))|\le (|F(0)|+L)\,(1+|u(x)|).
\end{align}
Since \(u\in L^p(\Omega)\), it follows from \eqref{eq-Fu-bdd} that \(F\circ u\in L^p(\Omega)\) as well. Invoking the fact that \(\eta\) is compactly supported, we obtain \(\eta F(u)\in L^p(\rn)\).

Let \(\phi\) satisfy \eqref{eq-adm}.
Choose an open set \(V\Subset\Omega\) such that
\(\operatorname{supp}\eta\subset V\).
Then
\begin{align}\label{eq-etaFu-norm}
[\eta F(u)]_{{\mathsf \Lambda}^{\phi,\,p}(\rn)}^p = \iint_{\mathcal C_V}
       \frac{|\eta(x)F(u(x))-\eta(y)F(u(y))|^p}{|x-y|^{n+p}}\,\phi(|x-y|)\,
                                   dy\,dx.
\end{align}
Upon setting $\delta:=\operatorname{dist}(\operatorname{supp}\eta,\,V^c)>0$,  we have
\begin{align}\label{eq-etaFu-Vc}
  &\int_V\int_{V^c}
       \frac{|\eta(x)F(u(x))-\eta(y)F(u(y))|^p}{|x-y|^{n+p}}\,\phi(|x-y|)\,
                                   dy\,dx\\
  &\quad =\int_{\operatorname{supp}\eta}|\eta(x)F(u(x))|^p
          \left(\int_{V^c}\frac {\phi(|x-y|)}{|x-y|^{n+p}}\,dy\right)dx\notag\\
          &\quad \le \int_{\operatorname{supp}\eta}|\eta(x)F(u(x))|^p
          \left(\int_{\rn\setminus B_\delta(x)}\frac {\phi(|x-y|)}{|x-y|^{n+p}}\,dy\right)dx\notag\\
  &\quad =\omega_{n-1} \frac{\Phi_\infty(\delta)}{\delta^p}
                                  \|\eta F(u)\|_{L^p(\rn)}^p<\infty, \notag
\end{align}
where the last equality follows from \eqref{eq-est-tail}.

For \(x,y\in V\), we write
\begin{align*}
  \eta(x)F(u(x))-\eta(y)F(u(y))
  &=\eta(x)[F(u(x))-F(u(y))]+F(u(y))[\eta(x)-\eta(y)].
\end{align*}
Since \(\eta\) is a compactly supported Lipschitz function, it is bounded and, hence,
\[
  |\eta(x)-\eta(y)|\le C_\eta\min\{1,\,|x-y|\}.
\]
By this and \eqref{eq-Fu-bdd}, together with the Lipschitz property of \(F\), we obtain
\begin{align*}
  |\eta(x)F(u(x))-\eta(y)F(u(y))|
  &\le C\left(|u(x)-u(y)|+ (1+|u(y)|)\min\{1,\,|x-y|\}\right),
\end{align*}
where \(C=C(F,\eta,n)>0\). Consequently,
\begin{align*}
  &|\eta(x)F(u(x))-\eta(y)F(u(y))|^p\\
  &\quad \le C 2^{p-1}
  \left(|u(x)-u(y)|^p+(1+|u(y)|)^p\min\{1,\,|x-y|^p\}\right).
\end{align*}
Thus, applying \eqref{kernel:cutoff-bound} at the scale \(a=1\), we obtain
\begin{align}\label{eq-etaFu-V}
  &\int_V\int_V
        \frac{|\eta(x)F(u(x))-\eta(y)F(u(y))|^p}{|x-y|^{n+p}}\,\phi(|x-y|)\,\,dy\,dx\\
        &\quad \le C 2^{p-1}
        \int_V\int_V
        \frac{|u(x)-u(y)|^p}{|x-y|^{n+p}}\,\phi(|x-y|)\,\,dy\,dx\notag\\
        &\qquad
        +C  2^{p-1}\,\int_V\int_V
       (1+|u(y)|)^p \frac{ \min\{1,\,|x-y|^p\}}{|x-y|^{n+p}}\,\phi(|x-y|)\,\,dy\,dx\notag\\
  &\quad\le C  2^{p-1}\,
  \left(\mathcal E_{\phi,\,p}(u;\, \Omega)
  + \omega_{n-1}\Phi(1)\int_V(1+|u(y)|)^p\,dy\right)<\infty. \notag
\end{align}
Using the definition of \(\mathcal C_V\) in \eqref{eq-C_Omega}, we deduce from \eqref{eq-etaFu-Vc} and \eqref{eq-etaFu-V} that the seminorm in \eqref{eq-etaFu-norm} is finite. Consequently, \(\eta F(u)\in {\mathsf \Lambda}^{\phi,\,p}(\rn)\).
\end{proof}

\subsection{A critical Sobolev inequality in rearrangement form}\label{ss-f*-sobolev}
Let \(f^*\) denote the \emph{nonincreasing rearrangement} of a Lebesgue measurable function \(f\) on \(\rn\), defined by
\[
 f^*(t):=\inf\{\lambda\ge0:
 |\{x\in\mathbb R^n:\ |f(x)|>\lambda\}|\le t\},\qquad t>0.
\]
Comparing this with the notation \(f_B^*\) in \eqref{eq-f_B*}, we have
\[
f_B^*(t)=(f{\mathbf 1}_B)^*(t|B|)
\]
for any ball \(B\subset\rn\) and \(t>0\).

The following proposition presents a global Sobolev
inequality for nearly critical kernels. In contrast
to the classical Sobolev inequality for $W^{1,\,p}(\mathbb R^n)$,
the presence of the weight $\Phi$
yields an upper bound for $f^*$
at all scales.

\begin{proposition}\label{logsp:fullscale}
Let \(1<p<n\) and \(\phi\) be a nearly critical kernel. Set \(p^*:=np/(n-p)\).
Then there exists a constant \(C=C(n,p,\phi)>0\) such that for any \(f\in {\mathsf \Lambda}^{\phi,\,p}(\mathbb R^n)\),
\begin{equation}\label{logsp:fullscale-est}
\sup_{t>0} t^{1/p^*}\Phi(t^{1/n})^{1/p}f^*(t)
\le C[f]_{{\mathsf \Lambda}^{\phi,\,p}(\mathbb R^n)},
\end{equation}
where \(\Phi\) is the function defined in \eqref{eq-Phi}.
\end{proposition}

\begin{proof}
For any \(s>0\), let \(B_s:=B(0,s)\).
Fix a function \(0\le \eta\in C_c^\infty(B_1)\) satisfying
\(\int_{\mathbb R^n}\eta\,dx=1\). For any \(x,\, h\in\rn\) and \(s>0\), set
\[
  \eta_s(x):=s^{-n}\eta(x/s),\qquad f_s:=\eta_s*f,\qquad
  \Delta_hf(x):=f(x+h)-f(x).
\]

\medskip

\textit{Step 1:\, Estimates for $f_s$ and $\nabla f_s$.}
Since \(\eta_s\) is supported in \(B_s\) and has integral one, it follows that
\[
  f(x)-f_s(x)
  =\int_{|h|<s}\eta_s(h)[f(x)-f(x-h)]\,dh.
\]
Moreover, by the H\"older inequality and \(\int_{|h|<s} \eta_s(h)\,dh=1\), we have
\begin{align*}
  |f(x)-f_s(x)|^p
&\le \left(\int_{|h|<s}\eta_s(h)|f(x)-f(x-h)|^p\,dh\right)
\left(\int_{|h|<s}\eta_s(h)\,dh\right)^{p-1}\\
&=\int_{|h|<s}\eta_s(h)|f(x)-f(x-h)|^p\,dh.
\end{align*}
This, together with the Fubini theorem, gives
\begin{align}\label{eq-f-fs}
  \|f-f_s\|_{L^p(\rn)}^p
  &\le\int_{\mathbb R^n}\int_{|h|<s}
         \eta_s(h)|f(x)-f(x-h)|^p\,dh\,dx\\
  &=\int_{|h|<s}\eta_s(h)
         \|f(\cdot)-f(\cdot-h)\|_{L^p}^p\,dh\notag\\
  &\le \|\eta\|_{L^\infty(\rn)}\,s^{-n}\int_{|h|<s}\|\Delta_hf\|_{L^p}^p\,dh,  \notag
\end{align}
where translation invariance gives
\(\|f(\cdot)-f(\cdot-h)\|_{L^p(\rn)}=\|\Delta_hf\|_{L^p(\rn)}\).

In a way similar to the estimate of \eqref{eq-f-fs}, we now treat
\(\|\nabla f_s\|_{L^p(\rn)}^p\). Using \(\int_{\mathbb R^n}\nabla\eta_s(h)\,dh=0\), we write
\[
  \nabla f_s(x)=\int_{|h|<s}\nabla\eta_s(h)
                                      [f(x-h)-f(x)]\,dh.
\]
By H\"older's inequality and the fact \(\|\nabla\eta_s\|_{L^1(\rn)}=s^{-1}\|\nabla\eta\|_{L^1(\rn)}\), we have
\begin{align*}
  |\nabla f_s(x)|^p
  &\le\left(\int_{|h|<s}|\nabla\eta_s(h)|\,dh\right)^{p-1}
       \left(\int_{|h|<s}|\nabla\eta_s(h)|
                                  |f(x-h)-f(x)|^p\,dh\right)\\
  &\le  \|\nabla\eta\|_{L^1(\rn)}^{p-1}\|\nabla\eta\|_{L^\infty(\rn)} \,s^{-n-p} \int_{|h|<s}|f(x-h)-f(x)|^p\,dh.
\end{align*}
Integrating with respect to \(x\) and using the change of variables \(h \mapsto -h\), we obtain
\begin{equation}\label{eq-nabla-fs}
  \|\nabla f_s\|_{L^p(\rn)}^p
  \le C s^{-n-p}\int_{|h|<s}\|\Delta_hf\|_{L^p(\rn)}^p\,dh,
\end{equation}
where the constant \(C\) depends only on \(n,p\) and \(\eta\).

\medskip

\textit{Step 2:\,  A rearrangement estimate at arbitrary scales.}
For each \(s>0\), by the Young inequality, we have
\[
  \|f_s\|_{L^p(\rn)}\le\|f\|_{L^p(\rn)},\qquad
  \|\nabla f_s\|_{L^p(\rn)}\le s^{-1}\|\nabla\eta\|_{L^1(\rn)}\,\|f\|_{L^p(\rn)},
\]
which implies that \(f_s\in W^{1,p}(\mathbb R^n)\).
Then, the classical Sobolev inequality (see \cite{AdamsFournier2003}) gives
\[
  \|f_s\|_{L^{p^*}(\rn)}\le C_{n,p}\|\nabla f_s\|_{L^p(\rn)},
\]
where \(p^*=np/(n-p)\).

With this fact, the properties of the rearrangement yield that for any \(t>0\) and \(s>0\),
\begin{align*}
  f^*(t)
  &\le (f-f_s)^*(t/2)+f_s^*(t/2)\\
  &\le (2/t)^{1/p}\|f-f_s\|_{L^p(\rn)}
           +(2/t)^{1/p^*}\|f_s\|_{L^{p^*}(\rn)} \notag\\
           &\le (2/t)^{1/p}\|f-f_s\|_{L^p(\rn)}
           +C_{n,\,p}(2/t)^{1/p^*}\|\nabla f_s\|_{L^p(\rn)}.
\end{align*}
Consequently, by \((a+b)^p\le 2^{p-1}(a^p+b^p)\) for all $a,b>0$,
\begin{align*}
 [f^*(t)]^p
  \le C \left(t^{-1}\|f-f_s\|_{L^p(\rn)}^p
           +t^{-p/p^*}\|\nabla f_s\|_{L^p(\rn)}^p\right).
\end{align*}
Inserting \eqref{eq-f-fs} and \eqref{eq-nabla-fs} into the right-hand side of this last estimate, we obtain
\begin{align}\label{eq-f*t-pth}
 [f^*(t)]^p
  &\le C\left(t^{-1} s^{-n}
  + t^{-p/p^*}   s^{-n-p} \right)\,\left(\int_{|h|<s}\|\Delta_hf\|_{L^p(\rn)}^p\,dh\right),
\end{align}
where \(C\) is a positive constant depending only on \(n\) and \(p\).

\medskip

\textit{Step 3:\,  Integration with respect to the mollification parameter.}
Multiplying both sides of \eqref{eq-f*t-pth} by
\(t^{p/p^*} \phi(s)/s\) and integrating with respect
to \(s\) over \((0,t^{1/n})\), by the definition of $\Phi_0$ in \eqref{eq-Phi0infty}, we obtain
\begin{align}\label{eq-Phi0-reag}
  & t^{p/p^*}\Phi_0(t^{1/n})[f^*(t)]^p\\
  &\quad =t^{p/p^*}[f^*(t)]^p\int_0^{t^{1/n}} \phi(s)\,\frac{ds}{s} \notag\\
  &\quad \le C
\int_0^{t^{1/n}}
  \left(t^{p/p^*-1} s^{-n}
  +  s^{-n-p} \right)\,\left(\int_{|h|<s}\|\Delta_hf\|_{L^p(\rn)}^p\,dh\right) \phi(s)\,\frac{ds}{s} \notag\\
  &\quad = C \int_{\rn}
  \left(\int_{|h|<s\le t^{1/n}} \left(t^{p/p^*-1} s^{-n}
  +  s^{-n-p} \right)\, \phi(s)\,\frac{ds}{s}\right)\,\|\Delta_hf\|_{L^p(\rn)}^p\,dh. \notag
\end{align}
We now estimate the inner integral. Since
\(p/p^*=1-p/n\), for \(|h|<s\le t^{1/n}\) we have
$$t^{p/p^*-1} s^{-n}
  +  s^{-n-p} =t^{-p/n}s^{-n}
  +  s^{-n-p} \le 2s^{-n-p}\le 2|h|^{-n} s^{-p}.$$
Further, by the estimate $\Phi_\infty(r)\le C \phi(r)$ (see \eqref{eq-Phiinfty-phi-Phi0}), we obtain
\begin{align*}
\int_{|h|<s\le t^{1/n}} \left(t^{p/p^*-1} s^{-n}
  +  s^{-n-p} \right)\, \phi(s)\,\frac{ds}{s}
&\le 2|h|^{-n}
\int_{|h|}^\infty \phi(s)\,s^{-p-1}\,ds\\
&=2|h|^{-n-p}\, \Phi_\infty(|h|)\\
&\le C |h|^{-n-p}\, \phi(|h|).
\end{align*}
Inserting this estimate into \eqref{eq-Phi0-reag},
and using the equivalence $\Phi_0\approx \Phi$ from \eqref{eq-Phi0-Phi}, we conclude that
\begin{align*}
t^{p/p^*}\Phi(t^{1/n})[f^*(t)]^p\le
C\int_{\mathbb R^n}
\|\Delta_hf\|_{L^p(\mathbb R^n)}^p
\frac{\phi(|h|)}{|h|^{n+p}}\,dh.
\end{align*}
Taking the \(p\)-th root and then the supremum over
\(t>0\), we obtain
\eqref{logsp:fullscale-est}.
\end{proof}

\subsection{Subcritical weak Harnack inequality via Ok and Song}\label{ss-OkSong}


Let \(p\in(1,\infty)\).
For the operator \(\mathcal L_{\phi,\,p}\) in \eqref{eq-Lphi} associated  with a kernel \(K_{\phi,\,p}\) satisfying \eqref{eq-ker-phi}, Ok and Song \cite{OkSong} established various Harnack inequalities, provided that \(\phi:\,[0,\infty)\to[0,\infty)\) satisfies the following conditions (see \cite[(1.5)--(1.7)]{OkSong}):
\begin{enumerate}[\rm (P1)]
    \item \(\phi(0)=0\), \(\phi(t)>0\) for \(t>0\);

    \item $\phi$ satisfies the Dini condition
$\int_0^1\phi(t)\,\frac{dt}{t}<\infty;$

 \item $\phi$ satisfies the two scaling conditions
\begin{equation}\label{condition2-Ok}
 \frac{\phi(t)}{t^{(1-s)p}}
 \le L\frac{\phi(\tau)}{\tau^{(1-s)p}},
 \qquad
 \frac{\phi(t)}{t^{(1-\widetilde s)p}}
 \ge L^{-1}\frac{\phi(\tau)}{\tau^{(1-\widetilde s)p}},
 \end{equation}
 for all \(0<\tau\le t<\infty\), and for some constants \(0<s<\widetilde s<\infty\) and \(L\ge1\).
\end{enumerate}
From (P2) and the first condition in \eqref{condition2-Ok} (which implies \(\phi(t)\le Ct^{(1-s)p}\) whenever \(t>1\)), it follows easily that
\[
 \int_0^\infty\min\{1, t^{-p}\}\,\phi(t)\,\frac{dt}{t}<\infty,
\]
that is, such \(\phi\) satisfies the integrability condition \eqref{eq-adm}.

For a nearly critical kernel $\phi$, applying \cite[Theorem~1.5]{OkSong} directly yields the following
\textit{subcritical weak Harnack inequality}.

\begin{proposition}
\label{prop-sub-wHI}
Let $p\in(1,\infty)$ and $\phi$ be a nearly critical kernel.  If  \(u \in \mathcal F^{\phi,\, p}(\Omega)\) is a  weak supersolution to
\[
\mathcal L_{\phi,\,p}u=0 \quad \text{in }\, \Omega
\]
and $u\ge 0$ a.e. in a ball \(B_{R}(x_0) \Subset \Omega\) with $x_0\in\rn$ and $R>0$,
then for  any
\begin{equation}\label{q range}
0<q<\begin{cases}
q_{\rm H}=\dfrac{n(p-1)}{n-p}\quad & \text{if}\ \, 1<p<n;\\
\infty \quad& \text{if}\ \, p\ge n,
\end{cases}
\end{equation}
and any  \(0<r\le R/2\),
\begin{equation}\label{regularity:weak-Harnack}
 \left(\fint_{B_{r/2}(x_0)}u^q\,dx\right)^{1/q}
 \le C\left(\operatorname*{ess\,inf}_{B_r(x_0)}u
 +\left(\frac{r^p\Phi (R)}{R^p\Phi (r)}\right)^{1/(p-1)}
     \mathsf{Tail}^{\phi}_{p} (u_-;x_0,R)\right),
\end{equation}
where \(\Phi\) is the function defined in \eqref{eq-Phi}, and $C=C(n,p,\phi,\Lambda,q)$ is independent of $x_0, r, R$, $\Omega$ and  $u$.
\end{proposition}

\begin{proof}
If $\phi$ is a nearly critical kernel, then it  satisfies properties (P1) and (P2), and, moreover, condition \eqref{condition2-Ok} holds for arbitrary \(s\in(0,1)\) and for \(\widetilde s=1\).
By Lemma~\ref{energy:tests}, the required tail is finite.  Thus, applying \cite[Theorem~1.5]{OkSong} to the kernel \(K_{\phi,\,p}\) yields
that
\eqref{regularity:weak-Harnack} holds whenever
$$
0<q<\begin{cases}
\dfrac{n(p-1)}{n-sp}\quad & \text{if}\ \, sp<n;\\
\infty \quad& \text{if}\ \, sp\ge n.
\end{cases}
$$
We note that, the primitive used in \cite{OkSong} is $\Phi_0$, and \eqref{eq-Phi0-Phi} gives $\Phi_0\approx\Phi$.
Due to the arbitrariness of $s\in(0,1)$, we obtain that \eqref{regularity:weak-Harnack} holds for $q$ satisfying \eqref{q range}.
\end{proof}

\begin{remark}

In general, the range of the exponent $s$ in \eqref{condition2-Ok} can not exceed $1$.
For example, when \(\gamma>1\), the logarithmic kernel
 \(\phi_\gamma(t)=(\log\frac{e}{t\wedge1})^{-\gamma}\)  is nearly critical;
however,
 the first inequality in \eqref{condition2-Ok} fails for \(s\ge 1\).
Consequently, when \(1<p<n\), the weak Harnack inequality of Ok and Song
\cite[Theorem~1.5]{OkSong} yields estimates for \(0<q<q_{\mathrm H}\)
but does not yield the endpoint estimate at \(q=q_{\mathrm H}\).
  A natural question is what happens at the critical index \(q=q_{\mathrm H}\), namely, what is the appropriate version of the weak Harnack inequality in this case. This is precisely what Theorem \ref{thm-wHI-qH-phi} addresses.
\end{remark}

\subsection{Full-space truncation energy with linear level growth}\label{sec-linear-growth}

We start with the following Caccioppoli-type estimate, which is taken from \cite[Lemma~5.4]{OkSong}; see also \cite[Lemma~5.1]{DiCastro-Kuusi-Palatucci}.

\begin{lemma}\label{Caccioppoli}
Let $p\in(1,\infty)$ and $\phi$ be a nearly critical kernel.
Suppose that  \(u \in \mathcal F^{\phi,\, p}(\Omega)\) is a  weak supersolution to \eqref{eq-Lgammap-phi} and $u\ge 0$ a.e. in a ball \(B_{R}(x_0) \Subset \Omega\) with center $x_0\in\rn$ and radius $R>0$. Set
$$v:=u+d,$$
where $d>0$. If $a\in(0,p-1)$, then for any $r\in(0, R/2)$ and  nonnegative $\chi\in C_c^\infty(B_{r}(x_0))$,
\begin{align}\label{caccioppoli}
  &\int_{B_{r}(x_0)}\int_{B_{r}(x_0)}
       |v(x)^{a/p}\chi(x)-v(y)^{a/p}\chi(y)|^p
                                      K_{\phi,\, p}(x,y)\,dy\,dx\\\notag
  &\quad\le C\int_{B_{r}(x_0)}\int_{B_{r}(x_0)}
       \max\{v(x)^a,v(y)^a\}|\chi(x)-\chi(y)|^p
                                      K_{\phi,\, p}(x,y)\,dy\,dx\\\notag
  &\qquad+C\left[\sup_{x\in\operatorname{supp}\chi}
                  \int_{B_{r}(x_0)^c}K_{\phi,\, p}(x,y)\,dy
       +d^{1-p}\frac{\Phi(R)}{R^p}
                       [\mathsf{Tail}^{\phi}_{p}(u_-;x_0,R)]^{p-1}\right]
                  \int_{B_{r}}v^a\chi^p\,dx,\notag
\end{align}
where \(\Phi\) is defined as in \eqref{eq-Phi}, and $C=C(a,n,p,\phi,\Lambda)$ is independent of $u,x_0,r,R,d,\Omega$.

\end{lemma}

\begin{proof}
The estimate \eqref{caccioppoli} is a direct application of \cite[Lemma~5.4]{OkSong}.
Its proof follows essentially from \cite[Lemma~5.1]{DiCastro-Kuusi-Palatucci}, by applying the definition of a supersolution (see Definition~\ref{def-supersol}) to the test function
\[
 \zeta:=\chi^p(u_++d)^{a-p+1}.
\]
Indeed, \(\zeta\in \mathcal F^{\phi,\,p}(\Omega)\). To see this, we note that \(F(t):=(t_++d)^{a-p+1}\) is Lipschitz, with
\[
 |F(t)-F(s)|\le(p-1-a)d^{a-p}|t-s| \qquad \text{for all } t,s\in\mathbb R,
\]
so that Lemma~\ref{lem-stability} can be applied.
\end{proof}

The main goal of this subsection is to establish a full-space energy estimate
for the localized Lipschitz truncation
$w_k:=\eta \min\{(u+d)_+,\,k\}$, with the energy growing linearly
in the truncation height $k$.
This estimate will play a key role in establishing  Theorem~\ref{thm-wHI-qH}.

\begin{proposition}\label{LH:linear-truncation}
Let $p\in(1,\infty)$ and $\phi$ be a nearly critical kernel.
 Suppose that  \(u \in \mathcal F^{\phi,\,p}(\Omega)\) is a  weak supersolution to \eqref{eq-Lgammap-phi} and $u\ge 0$ a.e. in \(B_{R}(x_0) \Subset \Omega\) with $x_0\in\rn$ and $R>0$.
Let $0<r\le 10^{-10n}R$ and
\begin{align}\label{eq-d>tail}
d>0, \qquad
  d\ge\left(\frac{r^p\Phi(R)}{R^p\Phi(r)}\right)^{1/(p-1)}
                       \mathsf{Tail}^{\phi}_{p}(u_-;x_0,R),
\end{align}
where \(\Phi\) is the function defined in \eqref{eq-Phi}.
Choose $\eta\in C_c^\infty(B_{3r/2}(x_0))$ such that
$0\le\eta\le1$, $\eta=1$ on $B_r(x_0)$, and $|\nabla\eta|\le 4/r$.
For $k\in(0,\infty)$, set
\[
v=u+d\qquad\text{and}\qquad  w_k=\eta\min\{v_+,\,k\}.
\]
Then $w_k\in\Lambda^{\phi,\,p}(\mathbb R^n)$ and
\begin{equation}\label{LH:localized-energy}
\begin{aligned}
  \int_{\mathbb R^n}\int_{\mathbb R^n}
       |w_k(x)-w_k(y)|^pK_{\phi,\, p}(x,y)\,dy\,dx
 \le Ck r^{n-p}\Phi(r)
        \left(\operatorname*{ess\,inf}_{B_{2r}(x_0)}v\right)^{p-1},
\end{aligned}
\end{equation}
where the constant $C=C(n,p,\phi,\Lambda)>0$ is independent of $u,r,R,d,k$ and $\Omega$.
\end{proposition}

The proof of Proposition \ref{LH:linear-truncation} relies on the following Lemma~\ref{v plus d} and Lemma~\ref{LH:flux-tail}.
Specifically, Lemma~\ref{v plus d} provides a weak Harnack inequality for the shifted supersolution \(v=u+d\).

\begin{lemma}\label{v plus d}
Let all the assumptions be as in Proposition \ref{LH:linear-truncation}.
Then, for any $q$ satisfying \eqref{q range},
and any $\kappa\ge 1$ with $\kappa r< R/2$,
\begin{equation}\label{LH:truncation-moments}
  \left(\fint_{B_{\kappa r/2}(x_0)}(u+d)^{q}\,dx\right)^{1/q}
  \le C\operatorname*{ess\,inf}_{B_{\kappa r}(x_0)}(u+d),
\end{equation}
where $C=C(n,p,q,\phi,\Lambda,\kappa)$ is independent of $u,x_0,r,R,d,\Omega$.

\end{lemma}

\begin{proof}
Apply Proposition~\ref{prop-sub-wHI} to the function $u$, with $\kappa r$ in place of $r$.
Moreover, by \eqref{eq-Phi-power-comparison} and the choice of $d$ in \eqref{eq-d>tail}, we have
\begin{align*}
  &\left(\frac{(\kappa r)^p\Phi(R)}{R^p\Phi(\kappa r)}\right)^{1/(p-1)}
                                \mathsf{Tail}^{\phi}_{p}(u_-;x_0,R)\\
  &\qquad=\kappa^{p/(p-1)}
       \left(\frac{\Phi(r)}{\Phi(\kappa r)}\right)^{1/(p-1)}
       \left(\frac{r^p\Phi(R)}{R^p\Phi(r)}\right)^{1/(p-1)}
                                \mathsf{Tail}^{\phi}_{p}(u_-;x_0,R)\\
  &\qquad\le C_{\kappa, p}\,d.
\end{align*}
It follows that
\begin{align*}
  \left(\fint_{B_{\kappa r/2}(x_0)}u^{q}\,dx\right)^{1/q}
  \le C\left(\operatorname*{ess\,inf}_{B_{\kappa r}(x_0)}u+\kappa^{p/(p-1)}d\right).
\end{align*}
For $A\ge0$, since $d > 0$, we have
\[
  (A+d)^{q}\le\max\{1,\,2^{q-1}\}(A^{q}+d^{q}).
\]
Thus,
\begin{align*}
  \fint_{B_{\kappa r/2}(x_0)}(u+d)^{q}\,dx
  &\le \max\{1,\,2^{q-1}\}\left(\fint_{B_{\kappa r/2}(x_0)}u^{q}\,dx+d^{q}\right)\\
  &\le C\left(\operatorname*{ess\,inf}_{B_{\kappa r}(x_0)}(u+d)\right)^{q},
\end{align*}
which proves \eqref{LH:truncation-moments}.
\end{proof}

\begin{lemma}\label{LH:flux-tail}
Let all the assumptions be as in Proposition \ref{LH:linear-truncation}. Then for every $\zeta\in C_c^\infty(B_{3r}(x_0))$ with
$0\le\zeta\le1$ and $|\nabla\zeta|\le4/r$,
\begin{equation}\label{LH:cutoff-flux}
\begin{aligned}
  &\int_{B_{4r}(x_0)}\int_{B_{4r}(x_0)}
       |v(x)-v(y)|^{p-1}|\zeta(x)-\zeta(y)|
                         K_{\phi,\, p}(x,y)\,dy\,dx\\
  &\qquad\le Cr^{n-p}\Phi(r)
              \left(\operatorname*{ess\,inf}_{B_{4r}(x_0)}v\right)^{p-1},
\end{aligned}
\end{equation}
where $C=C(n,p,\phi,\Lambda)$ is independent of $u,x_0,r,R,d,\Omega$.
\end{lemma}

\begin{proof}
 Fix $x_0\in\rn$. We simply write $B_t=B_t(x_0)$ for all $t> 0$.
Since $u\ge0$ on
$B_R$, we have $v\ge d>0$ on $B_R$ and hence, for $0<r\le 10^{-10n}R$, $$\operatorname*{ess\,inf}_{B_{4r}}v\ge d>0.$$

Let $a\in(0,p-1)$ and, when $p<n$, assume additionally that
$a>p-\frac{n}{n-p}.$
 If $A>B>0$,  then
\[
  A^{a/p}-B^{a/p}
  =\frac ap\int_B^A t^{a/p-1}\,dt
  \ge\frac ap A^{a/p-1}(A-B).
\]
By symmetry, for all $A,B>0$, we have
\[
  |A^{a/p}-B^{a/p}|
    \ge\frac ap \left(\max\{A, B\}\right)^{a/p-1}|A-B|
\]
and, hence,
\begin{align*}
  |A-B|^{p-1}
  &\le\left(\frac pa\right)^{p-1}
       |A^{a/p}-B^{a/p}|^{p-1}
       \max\{A,\,B\}^{\frac{(p-a)(p-1)}{p}}\\
  &\le C_a|A^{a/p}-B^{a/p}|^{p-1}(A+B)^{\frac{(p-a)(p-1)}{p}}.
\end{align*}
Applying this last inequality with $A=v(x)$ and $B=v(y)$, we obtain
\begin{align*}
  |v(x)-v(y)|^{p-1}\le C_a|v(x)^{a/p}-v(y)^{a/p}|^{p-1}(v(x)+v(y))^{\frac{(p-a)(p-1)}{p}}.
\end{align*}
This, combined with the H\"older inequality, yields
\begin{align*}
  &\int_{B_{4r}}\int_{B_{4r}}
       |v(x)-v(y)|^{p-1}|\zeta(x)-\zeta(y)|
                                      K_{\phi,\, p}(x,y)\,dy\,dx\\
    &\quad\le C_a
       \int_{B_{4r}}\int_{B_{4r}}|v(x)^{a/p}-v(y)^{a/p}|^{p-1}(v(x)+v(y))^{\frac{(p-a)(p-1)}{p}}|\zeta(x)-\zeta(y)|
       K_{\phi,\, p}(x,y)\,dy\,dx
       \notag\\
  &\quad\le C_a
       \left(\int_{B_{4r}}\int_{B_{4r}}
       |v(x)^{a/p}-v(y)^{a/p}|^pK_{\phi,\, p}(x,y)\,dy\,dx
       \right)^{\frac{p-1}{p}}\notag\\
  &\qquad\times
       \left(\int_{B_{4r}}\int_{B_{4r}}
       (v(x)+v(y))^{(p-a)(p-1)}|\zeta(x)-\zeta(y)|^p
                                     K_{\phi,\, p}(x,y)\,dy\,dx
       \right)^{\frac{1}{p}}. \notag
\end{align*}
In view of this last estimate, the proof of \eqref{LH:cutoff-flux} can be reduced to verifying that
\begin{align}\label{eq-I1-est}
{\rm I}_1&:=\int_{B_{4r}}\int_{B_{4r}}
       |v(x)^{a/p}-v(y)^{a/p}|^pK_{\phi,\, p}(x,y)\,dy\,dx\\
       &\le  Cr^{n-p}\Phi(r)
              \left(\operatorname*{ess\,inf}_{B_{4r}}v\right)^a \notag
\end{align}
and
\begin{align}\label{eq-I2-est}
{\rm I}_2&:=\int_{B_{4r}}\int_{B_{4r}}
       (v(x)+v(y))^{(p-a)(p-1)}|\zeta(x)-\zeta(y)|^p
                                     K_{\phi,\, p}(x,y)\,dy\,dx\\
                                     &
                                     \le Cr^{n-p}\Phi(r)
              \left(\operatorname*{ess\,inf}_{B_{4r}}v\right)^{(p-a)(p-1)}. \notag
\end{align}

Let us first prove \eqref{eq-I2-est}.  For all $s,t>0$, the
following elementary inequality holds:
\[
    (s+t)^{(p-a)(p-1)}\le C_{p,a}\left(s^{(p-a)(p-1)}+t^{(p-a)(p-1)}\right).
\]
Applying this inequality with $s=v(x)$ and $t=v(y)$,
and using the symmetry of the integrand together with
\eqref{LH:global-cutoff}, we obtain
\begin{align*}
 \mathrm{I}_2
  &\le C_{p,a}\int_{B_{4r}}v(y)^{(p-a)(p-1)}
      \left(\int_{B_{4r}}|\zeta(x)-\zeta(y)|^p
                                      K_{\phi,\, p}(x,y)\,dx\right)dy\\
  &\le C\frac{\Phi(r)}{r^p}\int_{B_{4r}}v^{(p-a)(p-1)}\,dx\\
  &\le Cr^{n-p}\Phi(r)\left(\operatorname*{ess\,inf}_{B_{4r}}v\right)^{(p-a)(p-1)}.
\end{align*}
For the last inequality, we have applied the subcritical weak Harnack
estimate \eqref{LH:truncation-moments} with $\kappa=8$, using $B_{4r}\subset B_{8r}$. It is
automatically applicable when $p\ge n$, whereas, when $p<n$,
our assumption on $a$ ensures that
\[
(p-a)(p-1)<q_{\mathrm H}.
\]
This completes the proof of \eqref{eq-I2-est}.

It remains to prove \eqref{eq-I1-est}. Take a function  $\chi\in C_c^\infty(B_{6r})$ such that
\[
  0\le\chi\le1,\qquad \chi=1\text{ on }B_{4r},\qquad
  |\nabla \chi|\le C/r.
\]
Since $\chi=1$ on $B_{4r}$, we write
$$
{\rm I}_1=\int_{B_{4r}}\int_{B_{4r}}
       |v(x)^{a/p}\chi(x)-v(y)^{a/p}\chi(y)|^p
                                      K_{\phi,\, p}(x,y)\,dy\,dx.
$$
Applying Lemma~\ref{Caccioppoli} on the integral domain $B_{8r}\times B_{8r}$ yields
\begin{align}\label{eq-I1-caccioppoli}
  {\rm I}_1&\le \int_{B_{8r}}\int_{B_{8r}}
       |v(x)^{a/p}\chi(x)-v(y)^{a/p}\chi(y)|^p
                                      K_{\phi,\, p}(x,y)\,dy\,dx\\
  &\le C\int_{B_{8r}}\int_{B_{8r}}
       \max\{v(x)^a,v(y)^a\}|\chi(x)-\chi(y)|^p
                                      K_{\phi,\, p}(x,y)\,dy\,dx \notag\\
  &\quad+C\left[\sup_{x\in\operatorname{supp}\chi}
                  \int_{B_{8r}^c}K_{\phi,\, p}(x,y)\,dy
       +d^{1-p}\frac{\Phi(R)}{R^p}
                       [\mathsf{Tail}^{\phi}_{p}(u_-;x_0,R)]^{p-1}\right]
                  \int_{B_{8r}}v^a\,dx. \notag
\end{align}
 By symmetry and \eqref{LH:global-cutoff}, we have
\begin{align*}
  &\int_{B_{8r}}\int_{B_{8r}}
       \max\{v(x)^a,v(y)^a\}|\chi(x)-\chi(y)|^p
                                      K_{\phi,\, p}(x,y)\,dy\,dx\\
  &\quad\le2\int_{B_{8r}}v(y)^a
       \left(\int_{B_{8r}}|\chi(x)-\chi(y)|^p
                                      K_{\phi,\, p}(x,y)\,dx\right)dy\\
  &\quad\le C\frac{\Phi(r)}{r^p}\int_{B_{8r}}v^a\,dx.
\end{align*}
Since $\operatorname{supp}\chi\subset B_{6r}$, it follows from Lemma~\ref{1.3}
and \eqref{eq-Phiinfty-phi-Phi0} that
\[
  \sup_{x\in\operatorname{supp}\chi}
      \int_{B_{8r}^c}K_{\phi,\, p}(x,y)\,dy
  \le\sup_{x\in\rn}\int_{|x-y|>r}K_{\phi,\, p}(x,y)\,dy
  \le C\frac{\Phi(r)}{r^p}.
\]
Also, the assumption on $d$ in \eqref{eq-d>tail} is equivalent to
\[
  d^{1-p}\frac{\Phi(R)}{R^p}
                 [\mathsf{Tail}^{\phi}_{p}(u_-;x_0,R)]^{p-1}
  \le\frac{\Phi(r)}{r^p}.
\]
Combining the preceding three estimates with
\eqref{eq-I1-caccioppoli}, we obtain
\begin{align*}
\mathrm I_1
\le C\frac{\Phi(r)}{r^p}
\int_{B_{8r}}v^a\,dx
\le Cr^{n-p}\Phi(r)
\fint_{B_{8r}}v^a\,dx.
\end{align*}
Applying
\eqref{LH:truncation-moments} with $q=a$ and
$\kappa=16$, we obtain \eqref{eq-I1-est}.

Finally, combining the estimates in \eqref{eq-I1-est}
and \eqref{eq-I2-est} yields \eqref{LH:cutoff-flux}.
\end{proof}

\begin{proof}[Proof of Proposition \ref{LH:linear-truncation}]
For any  $k>0$ and $t\in\mathbb R$, let
$T_k(t):=\min\{t_+,\,k\}$, which is a Lipschitz
function. So, applying Lemma~\ref{lem-stability} gives
$w_k=\eta T_k(v)\in\Lambda^{\phi,\,p}(\mathbb R^n).$
It remains to prove \eqref{LH:localized-energy},
which we divide into the following four steps.

\medskip

\textit{Step 1:\, Splitting the integral in
\eqref{LH:localized-energy} into three terms.}
For simplicity, we suppress the dependence on $x_0$ and
write $B_t:=B_t(x_0)$ for all $t>0$.
For any $x,y\in\mathbb R^n$, we have
\begin{align*}
  w_k(x)-w_k(y)
  =\eta(x)[T_k(v(x))-T_k(v(y))]
                +T_k(v(y))[\eta(x)-\eta(y)],
\end{align*}
Since $0\le\eta\le1$, this gives
\[  |w_k(x)-w_k(y)|^p
\le2^{p-1}|T_k(v(x))-T_k(v(y))|^p
+2^{p-1}T_k(v(y))^p|\eta(x)-\eta(y)|^p.\]
Further, using
$\operatorname{supp} w_k\subset B_{3r/2} \subset B_{2r}$ and the symmetry of $K_{\phi,\,p}$, we write the left-hand side of \eqref{LH:localized-energy} as follows:
\begin{align*}
  &\int_{\rn}\int_{\rn}
           |w_k(x)-w_k(y)|^pK_{\phi,\, p}(x,y)\,dy\,dx\\
           &\quad= \int_{B_{2r}}\int_{B_{2r}}
           |w_k(x)-w_k(y)|^pK_{\phi,\, p}(x,y)\,dy\,dx\\
           &\qquad+2\int_{B_{2r}}\int_{B_{2r}^c}
           |w_k(x)-w_k(y)|^pK_{\phi,\, p}(x,y)\,dy\,dx\\
 &\quad\le 2^{p-1} \int_{B_{2r}}\int_{B_{2r}}
           |T_k(v(x))-T_k(v(y))|^pK_{\phi,\, p}(x,y)\,dy\,dx\\
           &\qquad+2^{p-1}\int_{B_{2r}}\int_{B_{2r}}
           T_k(v(y))^p|\eta(x)-\eta(y)|^pK_{\phi,\, p}(x,y)\,dy\,dx\\
           &\qquad+2\int_{B_{2r}}\int_{B_{2r}^c}
           |w_k(x)-w_k(y)|^pK_{\phi,\, p}(x,y)\,dy\,dx\\
            &\quad=:\rm I+\rm II+\rm III.
\end{align*}

\medskip

\textit{Step 2:\, Estimates of the terms $\rm II$ and $\rm III$.}
Applying  \eqref{LH:global-cutoff} gives
\begin{align*}
  \rm II&=2^{p-1}\int_{B_{2r}}\int_{B_{2r}}
           T_k(v(y))^p|\eta(x)-\eta(y)|^pK_{\phi,\, p}(x,y)\,dy\,dx
\\
  &\le C\frac{\Phi(r)}{r^p}\int_{B_{2r}}T_k(v(y))^p\,dy.
\end{align*}
Moreover, since $\operatorname{supp}w_k\subset B_{3r/2}$
and $0\le w_k\le T_k(v)$, we deduce from \eqref{eq-est-K-tail}
that
\begin{align*}
\mathrm{III}
&=2\int_{B_{3r/2}}|w_k(x)|^p
\left(\int_{B_{2r}^c}K_{\phi,\,p}(x,y)\,dy\right)dx\\
&\le 2\int_{B_{3r/2}}|w_k(x)|^p
\left(\int_{|x-y|>r/2}K_{\phi,\,p}(x,y)\,dy\right)dx\\
&\le C\frac{\Phi(r)}{r^p}
\int_{B_{2r}}T_k(v(x))^p\,dx,
\end{align*}
where in the last step we have used
$\Phi_\infty(r/2)
\le C\Phi(r)$
by \eqref{eq-Phi-power-comparison} and \eqref{eq-Phiinfty-phi-Phi0}.

Since $u\ge0$ on $B_R$, we have
\[
\operatorname*{ess\,inf}_{B_{2r}}v\ge d>0.
\]
In particular, $v\ge0$ on $B_{2r}$, and hence
\(T_k(v)^p=\min\{v^p,k^p\}\le kv^{p-1}.\)
Combining this with \eqref{LH:truncation-moments}, we obtain
\[
\int_{B_{2r}}T_k(v)^p\,dx
\le k\int_{B_{2r}}v^{p-1}\,dx
\le Ckr^n
\left(\operatorname*{ess\,inf}_{B_{2r}}v\right)^{p-1}.
\]
Consequently, the preceding estimates of $\mathrm{II}$ and
$\mathrm{III}$ yield
\[
\mathrm{II}+\mathrm{III}
\le Ckr^{n-p}\Phi(r)
\left(\operatorname*{ess\,inf}_{B_{2r}}v\right)^{p-1}.
\]
This gives the desired bound for $\mathrm{II}+\mathrm{III}$.

\medskip

\textit{Step 3:\, A reduction for estimating $\rm I$.}
Choose $\xi\in C_c^\infty(B_{3r})$ such that
$0\le\xi\le1$, $\xi=1$  on $B_{2r}$, and $|\nabla\xi|\le 2/r$.
Define
\[
\psi_k:=\xi^p\bigl(k-T_k(v)\bigr).
\]
Since the map $t\mapsto k-T_k(t)$ is Lipschitz
and takes values in $[0,k]$, it follows from Lemma~\ref{lem-stability}
 that
\[
\psi_k\in\Lambda^{\phi,\,p}(\mathbb R^n),
\qquad
0\le\psi_k\le k\xi^p,
\qquad
\operatorname{supp}\psi_k\subset B_{3r} \subset\Omega.
\]
In particular, $\psi_k\in\mathcal F^{\phi,\,p}(\Omega)$.

Let $J_p(t):=|t|^{p-2}t$ for $t\in\mathbb R$.
Since $v=u+d\in\mathcal F^{\phi,\,p}(\Omega)$
is a weak supersolution to \eqref{eq-Lgammap-phi},
testing the supersolution inequality with $\psi_k$
yields
\begin{equation}\label{LH:descending-test}
0\le\int_{\mathbb R^n}\int_{\mathbb R^n}
J_p(v(x)-v(y))(\psi_k(x)-\psi_k(y))
K_{\phi,\,p}(x,y)\,dy\,dx.
\end{equation}
Moreover, by the H\"older inequality and
$\operatorname{supp}\psi_k\subset\Omega$,
\begin{align*}
&\int_{\mathbb R^n}\int_{\mathbb R^n}
|v(x)-v(y)|^{p-1}|\psi_k(x)-\psi_k(y)|
K_{\phi,\,p}(x,y)\,dy\,dx\\
&\quad =\iint_{\mathcal C_\Omega}
|v(x)-v(y)|^{p-1}|\psi_k(x)-\psi_k(y)|
K_{\phi,\,p}(x,y)\,dy\,dx\\
&\quad \le
\mathcal E_{\phi,\,p}(u;\Omega)^{\frac{p-1}{p}}
\mathcal E_{\phi,\,p}(\psi_k;\Omega)^{\frac1p}
<\infty.
\end{align*}
Thus the integral in \eqref{LH:descending-test}
is absolutely convergent.

The support condition $\operatorname{supp}\psi_k\subset
B_{3r}\subset B_{4r}$ and the symmetry of
$K_{\phi,\,p}$ allow us to rewrite the integral on
the right-hand side of \eqref{LH:descending-test} as
\begin{align}\label{eq-two-terms}
 &\int_{B_{4r}}\int_{B_{4r}}
      J_p(v(x)-v(y))(\psi_k(x)-\psi_k(y))
                                   K_{\phi,\, p}(x,y)\,dy\,dx\\
  &\quad+2\int_{B_{4r}}\int_{B_{4r}^c}
      J_p(v(x)-v(y))\psi_k(x)K_{\phi,\, p}(x,y)\,dy\,dx. \notag
\end{align}
In view of \eqref{LH:descending-test}, the expression
in \eqref{eq-two-terms} is nonnegative.

Next, we relate the first term in \eqref{eq-two-terms}
to the term $\rm I$.
For
all $A,B\in\mathbb R$, we have
\[
  (A-B)(T_k(A)-T_k(B))\ge0,\qquad
  |T_k(A)-T_k(B)|\le |A-B|,
\]
which implies
\begin{align}\label{LH:truncation-monotonicity}
  J_p(A-B)(T_k(A)-T_k(B))
  &=|A-B|^{p-1}|T_k(A)-T_k(B)|\\
  &\ge |T_k(A)-T_k(B)|^p. \notag
\end{align}
For any $x,y\in B_{4r}$, expanding
$\psi_k=\xi^p(k-T_k(v))$ and rearranging the terms, we obtain
\begin{align*}
  \psi_k(x)-\psi_k(y)
  &=-\frac{\xi(x)^p+\xi(y)^p}{2}
                          [T_k(v(x))-T_k(v(y))]\\
  &\quad+\frac{2k-T_k(v(x))-T_k(v(y))}{2}
                          [\xi(x)^p-\xi(y)^p].
\end{align*}
Combining this identity with $0\le T_k\le k$
and \eqref{LH:truncation-monotonicity} yields
\begin{align*}
  J_p(v(x)-v(y))(\psi_k(x)-\psi_k(y))
  &\le-\frac{\xi(x)^p+\xi(y)^p}{2}
                       |T_k(v(x))-T_k(v(y))|^p\\
  &\qquad+k|v(x)-v(y)|^{p-1}|\xi(x)^p-\xi(y)^p|.
\end{align*}
Note that $0\le\xi\le1$ and $\frac{\xi(x)^p+\xi(y)^p}{2}=1$ for $(x,y)\in B_{2r}\times B_{2r}$.
Integrating over $B_{4r}\times B_{4r}$ and dropping
the nonpositive contribution outside
$B_{2r}\times B_{2r}$, we obtain
\begin{align}\label{eq-two-terms-1st}
  &\int_{B_{4r}}\int_{B_{4r}}
       J_p(v(x)-v(y))(\psi_k(x)-\psi_k(y))
                                      K_{\phi,\, p}(x,y)\,dy\,dx\\
  &\quad\le-\int_{B_{2r}}\int_{B_{2r}}
       |T_k(v(x))-T_k(v(y))|^pK_{\phi,\, p}(x,y)\,dy\,dx\notag\\
  &\qquad+k\int_{B_{4r}}\int_{B_{4r}}
       |v(x)-v(y)|^{p-1}|\xi(x)^p-\xi(y)^p|
                                      K_{\phi,\, p}(x,y)\,dy\,dx. \notag
\end{align}
In particular, \eqref{eq-two-terms-1st} provides an upper bound
for the first term in \eqref{eq-two-terms}.

Combining \eqref{eq-two-terms-1st} with
\eqref{eq-two-terms}, and using the nonnegativity
of the expression in \eqref{eq-two-terms}, we deduce that
\begin{align*}
2^{1-p}\cdot {\rm I}&= \int_{B_{2r}}\int_{B_{2r}}
       |T_k(v(x))-T_k(v(y))|^pK_{\phi,\, p}(x,y)\,dy\,dx\\
       &\le k\int_{B_{4r}}\int_{B_{4r}}
       |v(x)-v(y)|^{p-1}|\xi(x)^p-\xi(y)^p|
                                      K_{\phi,\, p}(x,y)\,dy\,dx\\
                                      &\quad+2\int_{B_{4r}}\int_{B_{4r}^c}
      J_p(v(x)-v(y))\psi_k(x)K_{\phi,\, p}(x,y)\,dy\,dx\\
      &=:{\rm I}_1+{\rm I}_2.
\end{align*}
Thus, it remains to estimate $\mathrm I_1$ and
$\mathrm I_2$ to obtain the desired bound for $\mathrm I$.

\medskip

\textit{Step 4:\, Completion of the estimate for $\rm I$
via estimates of $\rm I_1$ and $\rm I_2$.}
Applying
Lemma~\ref{LH:flux-tail} with $\zeta=\xi$, and using $|\xi(x)^p-\xi(y)^p|\le p|\xi(x)-\xi(y)|$, gives
\begin{align*}
\mathrm I_1
&\le Ckr^{n-p}\Phi(r)
\left(\operatorname*{ess\,inf}_{B_{4r}}v\right)^{p-1}
\le Ckr^{n-p}\Phi(r)
\left(\operatorname*{ess\,inf}_{B_{2r}}v\right)^{p-1},
\end{align*}
  as desired.

Now, we turn to the estimate of ${\rm I}_2$.  By $\operatorname{supp}\psi_k\subset B_{3r}$, we have
\begin{align*}
{\rm I}_2
=2\int_{B_{3r}}\int_{B_{4r}^c}
       J_p(v(x)-v(y))\psi_k(x)K_{\phi,\, p}(x,y)\,dy\,dx.
       \end{align*}
For any $A\ge0$ and $B\in\mathbb R$, using the monotonicity of $J_p$ gives
\[
  J_p(A-B)\le J_p(A+B_-)
             =(A+B_-)^{p-1}
             \le C_p(A^{p-1}+B_-^{p-1}).
\]
Since $v(x)\ge d>0$  for a.e. $x\in B_{4r}\subset B_R$,  we apply this inequality with $A=v(x)$ and $B=v(y)$. With this and the fact $\psi_k\le k\xi^p\le k$, we obtain
\begin{align*}
{\rm I}_2
&\le Ck\int_{B_{3r}}v(x)^{p-1}
               \left(\int_{B_{4r}^c}K_{\phi,\, p}(x,y)\,dy\right)\,dx\\
  &\quad+Ck\int_{B_{3r}}               \left(\int_{B_{4r}^c}v_-(y)^{p-1}
                                     K_{\phi,\, p}(x,y)\,dy\right)\,dx. \notag
\end{align*}

For any $x\in B_{3r}$ and $y\notin B_{4r}$, we have $|x-y|>r$, so \eqref{eq-est-K-tail} and \eqref{eq-Phiinfty-phi-Phi0} give
\[
  \int_{B_{4r}^c}K_{\phi,\, p}(x,y)\,dy
  \le\int_{|x-y|>r}K_{\phi,\, p}(x,y)\,dy
    \le C\frac{\Phi(r)}{r^p}.
\]
This, combined with the weak Harnack estimate in Lemma~\ref{v plus d} (with $q=p-1$ and $\kappa=8$), further implies
\begin{align}\label{eq-I1-two-term1}
\int_{B_{3r}}v(x)^{p-1}
               \left(\int_{B_{4r}^c}K_{\phi,\, p}(x,y)\,dy\right)\,dx
&\le C\frac{\Phi(r)}{r^p} \int_{B_{3r}}v(x)^{p-1} \,dx\\
 & \le Cr^{n-p}\Phi(r)\left(\operatorname*{ess\,inf}_{B_{8r}}v\right)^{p-1}\notag\\
 &\le Cr^{n-p}\Phi(r)\left(\operatorname*{ess\,inf}_{B_{2r}}v\right)^{p-1}.\notag
\end{align}

By $|x_0-y|/4\le |x-y|\le 7|x_0-y|/4$ and the kernel comparison in Remark~\ref{rem-kerK}, for every $x\in B_{3r}$ and
$y\in B_{4r}^c$, we have
\[
K_{\phi,\,p}(x,y)\le C K_{\phi,\,p}(x_0,y).
\]
Moreover, since $v_-=0$ on $B_R$, $v_-\le u_-$
on $\mathbb R^n$, and $4r<R$, it follows that
\begin{align*}
\int_{B_{4r}^c}v_-(y)^{p-1}
K_{\phi,\,p}(x,y)\,dy
&=\int_{B_R^c}v_-(y)^{p-1}
K_{\phi,\,p}(x,y)\,dy\\
&\le C\int_{B_R^c}u_-(y)^{p-1}
K_{\phi,\,p}(x_0,y)\,dy.
\end{align*}
By the definition of the tail and the choice of $d$
in \eqref{eq-d>tail}, we further have
\begin{align*}
\int_{B_R^c}u_-(y)^{p-1}
K_{\phi,\,p}(x_0,y)\,dy
&\le\Lambda\,\frac{\Phi(R)}{R^p}
[\mathsf{Tail}^{\phi}_{p}(u_-;x_0,R)]^{p-1}\\
&\le\Lambda\,\frac{\Phi(r)}{r^p}d^{p-1}\\
&\le\Lambda\,\frac{\Phi(r)}{r^p}
\left(\operatorname*{ess\,inf}_{B_{2r}}v\right)^{p-1}.
\end{align*}
Consequently, we obtain
\begin{align}\label{eq-I1-two-term2}
\int_{B_{3r}}
\left(\int_{B_{4r}^c}v_-(y)^{p-1}
K_{\phi,\,p}(x,y)\,dy\right)\,dx
&\le Cr^n
\int_{B_R^c}u_-(y)^{p-1}
K_{\phi,\,p}(x_0,y)\,dy\\
&\le Cr^{n-p}\Phi(r)
\left(\operatorname*{ess\,inf}_{B_{2r}}v\right)^{p-1}. \notag
\end{align}

By \eqref{eq-I1-two-term1} and \eqref{eq-I1-two-term2}, we obtain
\begin{align*}
  \mathrm I_2
  \le Ck r^{n-p}\Phi(r)
  \left(\operatorname*{ess\,inf}_{B_{2r}}v\right)^{p-1}.
\end{align*}
Combined with the estimate for $\mathrm I_1$, this yields
\begin{align*}
  \mathrm I
  \le Ck r^{n-p}\Phi(r)
  \left(\operatorname*{ess\,inf}_{B_{2r}}v\right)^{p-1}.
\end{align*}

Finally, the estimates for $\mathrm I$ obtained in
\textit{Step 4}, together with those for $\mathrm{II}$ and
$\mathrm{III}$ in \textit{Step 2}, establish
\eqref{LH:localized-energy}.
\end{proof}

\section{Endpoint weak Harnack inequality }\label{sec-wHI}

%

The main aim of this section is to prove Theorems~\ref{thm-wHI-qH-phi} and \ref{thm-wHI-qH}, and to show that they imply the weak Harnack inequality \eqref{eq-wHI-q} in the subcritical regime \(0<q<q_{\mathrm H}\) when \(p\in(1,n)\).

\subsection{Proof of Theorem~\ref{thm-wHI-qH-phi}}\label{ss-wHI-qH-phi}

Indeed, Propositions~\ref{logsp:fullscale} and \ref{LH:linear-truncation}
are the two key ingredients in the proof of the
weak Harnack inequality at the endpoint exponent $q_{\rm H}$.

\begin{proof}[Proof of Theorem~\ref{thm-wHI-qH-phi}]
Fix a ball $B_R(x_0)\Subset\Omega$.  For any $0<r\le R/2$,  set
\[
  D:=\left(\frac{r^p\Phi(R)}{R^p\Phi(r)}\right)^{1/(p-1)}
                         \mathsf{Tail}^{\phi}_{p}(u_-;x_0,R).
\]
Note that the definition of the tail $\mathsf{Tail}^{\phi}_{p}$ gives
\[
  D^{p-1}=\frac{r^p}{\Phi(r)}
       \int_{\mathbb R^n\setminus B_R(x_0)}
           u_-(y)^{p-1}\frac{\phi(|x_0-y|)}{|x_0-y|^{n+p}}\,dy.
\]
Moreover, by \eqref{eq-Uu},  we have
\[\mathcal U_{\phi,\,p}:=\mathcal U_{\phi,\,p}(u,x_0,r,R)=\operatorname*{ess\,inf}_{B_r(x_0)}u+D.\]
To prove
\eqref{intro:endpoint-phi}, since $u_{B_r}^*(t)=0$ for $t\ge1$, it suffices
to prove
\begin{equation}\label{LH:endpoint-main}
  \sup_{0<t<1}t^{1/q_{\mathrm H}}
       \left(\frac{\Phi(rt^{1/n})}{\Phi(r)}\right)^{1/(p-1)}
                                u_{B_r(x_0)}^*(t)
  \le C\,\mathcal U_{\phi,\,p},
\end{equation}
where  $C=C(n,p,\phi,\Lambda)$ is independent of $x_0, r, R$, $\Omega$ and  $u$.

\medskip
\textit{Step 1:\, Establishing \eqref{LH:endpoint-main}
for small balls $B_r(x_0)$ with $0<r\le 10^{-10n}R$.}
For any $\epsilon>0$, let
\[
  d:=D+\epsilon,
  \qquad v:=u+d.
\]
Take $\eta\in C_c^\infty(B_{3r/2}(x_0))$ such that
$0\le\eta\le1$, $\eta=1$ on $B_r(x_0)$, and $|\nabla\eta|\le 4/r$.
For any $k>0$, let $$w_k:=\eta\min\{v_+,\,k\}$$ be the localized
truncation in Proposition~\ref{LH:linear-truncation}, and define
\[
  E_k:=\{x\in B_r(x_0):\,v(x)>k\},\qquad
  m(k):=\frac{|E_k|}{|B_r(x_0)|}.
\]

Consider first the case $|E_k|>0$.
Since $w_k=k$ on $E_k$ and $0\le w_k\le k$ on
$\mathbb R^n$, we have
\[
  w_k^*(|E_k|/2)=k.
\]
Indeed, for every $0\le\lambda<k$, the superlevel set
$\{x\in B_r(x_0):\, w_k(x)>\lambda\}$ contains $E_k$, and hence has measure
greater than $|E_k|/2$.
Applying Proposition~\ref{logsp:fullscale} to $f=w_k$ at $t=|E_k|/2$,
and using  the localized energy estimate \eqref{LH:localized-energy} in
Proposition \ref{LH:linear-truncation},
we obtain
\begin{align}\label{eq-vip}
  &k^p\left(\frac{|E_k|}{2}\right)^{(n-p)/n}
                  \Phi\left(\left(\frac{|E_k|}{2}\right)^{1/n}\right)\\
  &\quad\le C\int_{\mathbb R^n}\int_{\mathbb R^n}
         \frac{|w_k(x)-w_k(y)|^p\phi(|x-y|)}{|x-y|^{n+p}}\,dy\,dx\notag\\
  &\quad\le Ck r^{n-p}\Phi(r)\left(\operatorname*{ess\,inf}_{B_{r}(x_0)}v\right)^{p-1}. \notag
\end{align}
Recall that $|E_k|=\nu_n r^n m(k)$. By this and the fixed-ratio estimate in  \eqref{eq-Phi-power-comparison}, we have
\[
\Phi\left(\left(\frac{|E_k|}{2}\right)^{1/n}\right)
= \Phi\left(\left(\frac{\nu_n r^n m(k)}{2}\right)^{1/n}\right)
\ge C\Phi\left(r[m(k)]^{1/n}\right).
\]
Moreover,
\begin{align*}
 \operatorname*{ess\,inf}_{B_{r}(x_0)}v
  =\operatorname*{ess\,inf}_{B_{r}(x_0)}u+d = \mathcal U_{\phi,\,p}+\epsilon.
\end{align*}
Using these facts and $q_{\mathrm H}=\frac{n(p-1)}{n-p}$,
we derive from  \eqref{eq-vip} that
\begin{equation}\label{LH:endpoint-distribution}
  k\,m(k)^{1/q_{\mathrm H}}
       \left(\frac{\Phi(r[m(k)]^{1/n})}{\Phi(r)}\right)^{1/(p-1)}
  \le C\left(\mathcal U_{\phi,\,p}+\epsilon\right).
\end{equation}
If $|E_k|=0$, then $m(k)=0$ and, hence, the inequality \eqref{LH:endpoint-distribution}
 is trivial.

For any $0<t<1$ and $0<k<v_{B_r(x_0)}^*(t)$, the definition of
the normalized rearrangement in \eqref{eq-f_B*} gives
$m(k)>t$. Since $\Phi$ is increasing, we have
\begin{align*}
 k\,m(k)^{1/q_{\mathrm H}}
       \left(\frac{\Phi(r[m(k)]^{1/n})}{\Phi(r)}\right)^{1/(p-1)}
       \ge k\,t^{1/q_{\mathrm H}}
       \left(\frac{\Phi(rt^{1/n})}{\Phi(r)}\right)^{1/(p-1)}.
\end{align*}
Combining this with \eqref{LH:endpoint-distribution}
and letting $k\uparrow v_{B_r}^*(t)$, we find that
\[
  t^{1/q_{\mathrm H}}
       \left(\frac{\Phi(rt^{1/n})}{\Phi(r)}\right)^{1/(p-1)}
       v_{B_r(x_0)}^*(t)\le C\left(\mathcal U_{\phi,\,p}+\epsilon\right).
\]
From $0\le u\le v$ on $B_r$, it follows that
\[
  \sup_{0<t<1}t^{1/q_{\mathrm H}}
       \left(\frac{\Phi(rt^{1/n})}{\Phi(r)}\right)^{1/(p-1)}
       u_{B_r(x_0)}^*(t)\le C(\mathcal U_{\phi,\,p}+\epsilon).
\]
Letting $\epsilon\downarrow0$ proves \eqref{LH:endpoint-main}
when $r\le 10^{-10n}R$.

\medskip

\textit{Step 2:\, Establishing \eqref{LH:endpoint-main}
on shifted small balls.}
Let $0<r\le R/2$ and $\rho:=10^{-11n}r$.
Clearly, $\rho<10^{-10n}R/2$.
For any $z\in B_r(x_0)$, we have
\[
B_\rho(z)\Subset B_{R/2}(z)\Subset B_R(x_0).
\]
We may therefore apply \emph{Step 1} to $B_\rho(z)$,
taking $B_{R/2}(z)$ as the corresponding larger ball,
to obtain
\begin{align}\label{eq-wHI-Brhoz-half}
\sup_{0<t<1}t^{1/q_{\mathrm H}}
\left(\frac{\Phi(\rho t^{1/n})}
{\Phi(\rho)}\right)^{1/(p-1)}
u_{B_\rho(z)}^*(t)\le C\,
\mathcal U_{\phi,\,p}(u,z,\rho,R/2).
\end{align}

To estimate the tail part contained in $\mathcal U_{\phi,\,p}(u,z,\rho,R/2)$,
by \eqref{Def-tail}, we write
\begin{align*}
\left(\frac{\rho^p\Phi(R/2)}{(R/2)^p\Phi(\rho)}\right)
                        \left( \mathsf{Tail}^{\phi}_{p}(u_-; z,R/2)\right)^{p-1}
&\le\Lambda\,\frac{\rho^p}{\Phi(\rho)}
       \int_{\mathbb R^n\setminus B_{R/2}(z)}
                  u_-(y)^{p-1}K_{\phi,\, p}(z,y)\,dy.
\end{align*}
Since $u_-=0$ a.e. on $B_R(x_0)$,
the kernel comparison \eqref{eq-KKer-est} gives
\begin{align*}
  \int_{\mathbb R^n\setminus B_{R/2}(z)}
                  u_-(y)^{p-1}K_{\phi,\, p}(z,y)\,dy
  &=\int_{\mathbb R^n\setminus B_R(x_0)}
                  u_-(y)^{p-1}K_{\phi,\, p}(z,y)\,dy\\
  &\le C\int_{\mathbb R^n\setminus B_R(x_0)}
                  u_-(y)^{p-1}K_{\phi,\, p}(x_0,y)\,dy.
\end{align*}
Combining this with $\rho=10^{-11n}r$ and
\eqref{eq-Phi-power-comparison}, we obtain
\begin{align*}
 \left(\frac{\rho^p\Phi(R/2)}{(R/2)^p\Phi(\rho)}\right)
                        \left( \mathsf{Tail}^{\phi}_{p}(u_-; z,R/2)\right)^{p-1}
&\le C\frac{r^p}{\Phi(r)}
       \int_{\mathbb R^n\setminus B_{R}(x_0)}
                  u_-(y)^{p-1}K_{\phi,\, p}(x_0,y)\,dy\\
                  &\le C\Lambda D^{p-1}.
\end{align*}
Taking the $(p-1)$-st root and recalling
the definition of $\mathcal U_{\phi,\,p}$
in \eqref{eq-Uu}, we conclude that
\begin{align*}
\mathcal U_{\phi,\,p}(u,z,\rho,R/2)
\le C \left(\operatorname*{ess\,inf}_{B_\rho(z)}u+D\right).
\end{align*}
From this and \eqref{eq-wHI-Brhoz-half}, it follows that
\begin{align}\label{eq-wHI-Brhoz-D}
\sup_{0<t<1}t^{1/q_{\mathrm H}}
\left(\frac{\Phi(\rho t^{1/n})}
{\Phi(\rho)}\right)^{1/(p-1)}
u_{B_\rho(z)}^*(t)\le C\, \left(\operatorname*{ess\,inf}_{B_\rho(z)}u+D\right),
\end{align}
whenever $z\in B_r(x_0)$ and $\rho=10^{-11n}r$.

Next, we estimate $\operatorname*{ess\,inf}_{B_\rho(z)}u$.
Suppose $z,z'\in B_r(x_0)$ and $|z-z'|\le\rho/2$. Then
\[
  B_{3\rho/4}((z+z')/2)\subset B_\rho(z)\cap B_\rho(z').
\]
Consequently, for every $0\le\lambda<\operatorname*{ess\,inf}_{B_\rho(z)}u$,
\begin{align*}
  \frac{|\{x\in B_{\rho}(z'):\, u(x)>\lambda\}|}{|B_\rho(z')|}
  &\ge \frac{|B_\rho(z)\cap B_\rho(z')|}{|B_\rho(z')|}\\
  &  = \frac{|B_\rho(z)\cap B_\rho(z')|}{|B_\rho(z)|}>(3/4)^n=:\theta.
\end{align*}
This, along with the definition of rearrangement in \eqref{eq-f_B*}, gives
\[
  \operatorname*{ess\,inf}_{B_\rho(z)}u
  \le u_{B_\rho(z')}^*(\theta).
\]
Applying \eqref{eq-wHI-Brhoz-D} to the ball $B_\rho(z')$
at the level $\theta=(3/4)^n$, we obtain
\begin{align*}
u_{B_\rho(z')}^*(\theta)
&\le C \theta^{-1/q_{\mathrm H}}
\left(\frac{\Phi(\rho \theta^{1/n})}
{\Phi(\rho)}\right)^{-1/(p-1)}
 \left(\operatorname*{ess\,inf}_{B_\rho(z')}u+D\right)\\
 &\le C \,\left(\operatorname*{ess\,inf}_{B_\rho(z')}u+D\right),
\end{align*}
where the second inequality follows from \eqref{eq-Phi-power-comparison}.
Combining the last two estimates yields
\begin{equation}\label{LH:adjacent-infima}
  \operatorname*{ess\,inf}_{B_\rho(z)}u
  \le u_{B_\rho(z')}^*(\theta)
  \le  C \,\left(\operatorname*{ess\,inf}_{B_\rho(z')}u+D\right).
\end{equation}

For any $z,y\in B_r(x_0)$, we consider a chain of balls
$\{B_\rho(z_j)\}_{j=0}^N$ whose centers lie on the line
segment joining $z$ and $y$, with
$z_0=z$, $z_N=y$,
$|z_j-z_{j-1}|\le\rho/2$ for all $j=1,2,\dots,N$.
Moreover,
\[
N\le 1+\frac{2|z-y|}{\rho}
\le 1+\frac{4r}{10^{-11n}r}
<10^{12n}.
\]
Iterating \eqref{LH:adjacent-infima} along the chain of
balls, we obtain
\begin{align*}
  \operatorname*{ess\,inf}_{B_\rho(z)}u
  &\le C^N\operatorname*{ess\,inf}_{B_\rho(y)}u
             +D\sum_{j=0}^N C^j
             \le C
             \left(\operatorname*{ess\,inf}_{B_\rho(y)}u+D\right),
\end{align*}
where we have used the uniform bound
$N<10^{12n}$.
By the Lebesgue differentiation theorem, for a.e.
$y\in B_r(x_0)$, we have
\[
\operatorname*{ess\,inf}_{B_\rho(y)}u\le u(y).
\]
Consequently,
\begin{align*}
  \operatorname*{ess\,inf}_{B_\rho(z)}u
  \le    C  \left(u(y)+D\right).
\end{align*}
Taking the essential infimum over $y\in B_r(x_0)$ gives
\begin{align}\label{LH:shifted-inf}
  \operatorname*{ess\,inf}_{B_\rho(z)}u
  \le C\left(\operatorname*{ess\,inf}_{B_r(x_0)}u+D\right)=C\,\mathcal U_{\phi,\,p}(u; x_0, r, R).
\end{align}

Combining \eqref{eq-wHI-Brhoz-D} with
\eqref{LH:shifted-inf}, we conclude that, for
$\rho=10^{-11n}r$ and every $z\in B_r(x_0)$,
\begin{align}\label{eq-wHI-Brhoz}
\sup_{0<t<1}t^{1/q_{\mathrm H}}
\left(\frac{\Phi(\rho t^{1/n})}
{\Phi(\rho)}\right)^{1/(p-1)}
u_{B_\rho(z)}^*(t)\le C\,
\mathcal U_{\phi,\,p}(u,x_0,r,R).
\end{align}

\medskip

\textit{Step 3:\, Passing from local estimates on the covering
balls $B_\rho(z_j)$ to an estimate on $B_r(x_0)$.}
Still let $0<r\le R/2$ and $\rho=10^{-11n}r$.
Choose a maximal set
$\{z_1,\ldots,z_L\}\subset B_r(x_0)$ such that
\[
|z_i-z_j|\ge\frac{\rho}{2}
\qquad\text{whenever }i\ne j.
\]
Then, the balls
$\{B_{\rho/4}(z_j)\}_{j=1}^L$ are mutually disjoint, every $B_{\rho/4}(z_j)\subset B_{2r}(x_0)$,
and
$$B_r(x_0)\subset \bigcup_{j=1}^L B_\rho(z_j).$$
These facts, together with $r=10^{11n}\rho$, yield
\[
L\nu_n (10^{-11n}2^{-2} r)^n
 =\sum_{j=1}^L |B_{\rho/4}(z_j)|
=\left|\cup_{j=1}^L B_{\rho/4}(z_j) \right|
\le |B_{2r}(x_0)|
=\nu_n (2r)^n,\]
which implies
$L\le (2^{3}10^{11 n})^n$.

Let $0<t<1$ and
$\lambda>\max_{1\le j\le L}u_{B_\rho(z_j)}^*(9^{-n} t).$
Then, by the definition of the normalized rearrangement in \eqref{eq-f_B*}, for
each $j=1,\ldots,L$, we have
$$|\{x\in B_\rho(z_j):\,u(x)>\lambda\}|\le 9^{-n} t |B_\rho(z_j)|.$$
Consequently, using the covering
$B_r(x_0)\subset\bigcup_{j=1}^L B_\rho(z_j)$, we obtain
\begin{align*}
  |\{x\in B_r(x_0):\, u(x)>\lambda\}|
    &\le\sum_{j=1}^L |\{x\in B_\rho(z_j):\,u(x)>\lambda\}|\\
  &\le  9^{-n} t\sum_{j=1}^L |B_{\rho}(z_j)|\\
   &< t|B_r(x_0)|,
\end{align*}
where the last step is due to
$$9^{-n} \sum_{j=1}^L |B_{\rho}(z_j)|
= 9^{-n} 4^n \sum_{j=1}^L |B_{\rho/4}(z_j)|
\le 9^{-n} 4^n
 |B_{2r}(x_0)|<|B_r(x_0)|.$$
It therefore follows that
$u_{B_r(x_0)}^*(t)\le \lambda.$
Letting $\lambda$  decrease to $\max_{1\le j\le L} u_{B_\rho(z_j)}^*(9^{-n}t)$ yields
\[
  u_{B_r(x_0)}^*(t)\le\max_{1\le j\le L}u_{B_\rho(z_j)}^*(9^{-n}t).
\]
Meanwhile, by \eqref{eq-Phi-power-comparison} and $\rho=10^{-11n}r$, we have
\[\frac{\Phi(rt^{1/n})}{\Phi(r)}\le C\frac{\Phi(\rho(9^{-n}t)^{1/n})}{\Phi(\rho)}.\]
Combining the last two estimates with \eqref{eq-wHI-Brhoz} yields that for all $t\in(0,1)$ and
 $0<r\le R/2$,
\begin{align*}
  &t^{1/q_{\mathrm H}}
       \left(\frac{\Phi(rt^{1/n})}{\Phi(r)}\right)^{1/(p-1)}
                                       u_{B_r(x_0)}^*(t)\\
  &\quad\le C\max_{1\le j\le L}
       (9^{-n}t)^{1/q_{\mathrm H}}
       \left(\frac{\Phi(\rho (9^{-n}t)^{1/n})}{\Phi(\rho)}\right)^{1/(p-1)}
                              u_{B_\rho(z_j)}^*(9^{-n}t)\\
  &\quad\le C\mathcal U_{\phi,\,p}(u,x_0, r, R).
\end{align*}
Taking the supremum over
$t\in(0,1)$ proves \eqref{LH:endpoint-main}.
\end{proof}

\subsection{Proof of Theorem~\ref{thm-wHI-qH}}\label{ss-wHI-cor}


For the function $\Phi_\gamma$ defined in \eqref{eq-Phi-gamma},
we collect the following three basic properties that will be used
in this subsection and in Section~\ref{sec:proof-sharpness} below.

\begin{lemma}\label{lem-Phi-gamma}
Let $\gamma\in(1,\infty)$ and $\phi_\gamma$ be as in \eqref{eq-phi-gamma}.

\begin{enumerate}[\rm (i)]
    \item   For all $ r\in(0,1]$,
\begin{align}\label{eq-Phi-value}
 \Phi_\gamma(r)
 = \frac{[\log(e/r)]^{1-\gamma}}{\gamma-1}.
\end{align}

\item For all $r\in(0,\infty)$,
\begin{align}\label{eq-Phigamma0}
\Phi_\gamma(r)=\int_0^1 \phi_\gamma(rt)\,\frac{dt}{t}
\approx \int_0^\infty
\phi_\gamma(r\tau)\min\{1,\,\tau^{-p}\}\,\frac{d\tau}{\tau},
\end{align}
with implicit constants depending only on $\gamma$ and $p$.

\item For any $\kappa\in(1,\infty)$, there exist positive
constants $c$ and $C$, depending only on $\gamma$ and $\kappa$, such that,
for all $r\in(0,\infty)$,
\begin{align}\label{eq-doubling-Phi-gamma}
c\,\Phi_\gamma(r)\le \Phi_\gamma(\kappa r)\le C\,\Phi_\gamma(r).
\end{align}
\end{enumerate}
\end{lemma}

\begin{proof}
For $0<r\le1$, since $\gamma>1$, a direct calculation gives
\[
\Phi_\gamma(r)
=\int_0^r \phi_\gamma(\tau)\,\frac{d\tau}{\tau}
=\int_0^r \left(\log\frac e{\tau}\right)^{-\gamma}\,\frac{d\tau}{\tau}
=\frac{[\log(e/r)]^{1-\gamma}}{\gamma-1},
\]
which proves \eqref{eq-Phi-value}.

Observe that the integral on the right-hand side of \eqref{eq-Phigamma0}
is exactly the function $\Phi$ in \eqref{eq-Phi}, but now with $\phi$
replaced by $\phi_\gamma$. Combining this observation with the definition of $\Phi_0$ in \eqref{eq-Phi0infty},
and the equivalence $\Phi\approx \Phi_0$ in \eqref{eq-Phi0-Phi}, we obtain the equivalence in \eqref{eq-Phigamma0}.

Finally, the doubling property \eqref{eq-doubling-Phi-gamma} follows directly from \eqref{eq-Phi-power-comparison} and  \eqref{eq-Phi0-Phi}.
\end{proof}

\begin{proof}[Proof of Theorem~\ref{thm-wHI-qH}]
Note that $\phi_\gamma$ in \eqref{eq-phi-gamma} is a nearly critical kernel.
By \eqref{eq-Phigamma0}, the function $\Phi_\gamma$ is comparable to
the function $\Phi$ defined in \eqref{eq-Phi}, with $\phi$ therein replaced
by $\phi_\gamma$.
Thus, taking $\phi=\phi_\gamma$ in Theorem~\ref{thm-wHI-qH-phi} and applying
\eqref{eq-Phigamma0} directly yields Theorem~\ref{thm-wHI-qH}(i).

Now, we show Theorem~\ref{thm-wHI-qH}(ii). By \eqref{eq-Phi-value} and $\beta_\gamma=(\gamma-1)/(p-1)$, for all $r,t\in(0,1]$, we have
\begin{align*}
\left(\frac{\Phi_\gamma(rt^{1/n})}{\Phi_\gamma(r)}\right)^{\frac{1}{p-1}}
&=\left(\frac{\log(e/(rt^{1/n}))}{\log(e/r)}\right)^{-\beta_\gamma}
=\left(1+\frac{\log(1/t)}{n\log(e/r)}\right)^{-\beta_\gamma}
  \ge\left[\log(e/t)\right]^{-\beta_\gamma}.
\end{align*}
From this and \eqref{eq-LogL-fz}, it follows that for all $r\in(0,1]$,
\begin{align}\label{eq-u-LlogL}
    \|u\|_{L^{q_{\rm H},\infty}(\log L)^{-\beta_{\gamma}}(B_r(x_0))}
& = \sup_{t\in(0,1)}
 t^{\frac 1 {q_{\rm H}}}\left[\log(e/t)\right]^{-\beta_\gamma}u_{B_r(x_0)}^*(t)\\
 &\le \sup_{t \in (0, 1)}t^{\frac{1}{q_{\mathrm H}}}
 \left(\frac{\Phi_\gamma(rt^{1/n})}{\Phi_\gamma(r)}\right)^{\frac{1}{p-1}}
                         u_{B_r(x_0)}^*(t). \notag
\end{align}
The supremum on the right-hand side of \eqref{eq-u-LlogL}
is precisely the quantity on the left-hand side of
\eqref{intro:endpoint}.
Combining \eqref{intro:endpoint} with \eqref{eq-u-LlogL}
immediately yields \eqref{intro:LZ-functional} for
$0<r\le\min\{1,R/2\}$.

It remains to prove Theorem~\ref{thm-wHI-qH}(iii).
Let $0<r\le\min\{1,R/2\}$.
By \eqref{intro:LZ-functional}, we have
\begin{align}\label{eq-LlogL-wHI}
\sup_{t\in(0,1)}
 t^{\frac 1 {q_{\rm H}}}\left[\log(e/t)\right]^{-\beta_\gamma}u_{B_r(x_0)}^*(t)
 \le C_0 \,\mathcal U_{\gamma,\,p},
\end{align}
where $C_0>0$ is independent of $u$,  $x_0$, $r$, $R$ and $\Omega$.

If $\mathcal U_{\gamma,\, p}=0$, then \eqref{eq-LlogL-wHI} forces
 $u_{B_r(x_0)}^*\equiv 0$, and hence
$u=0$ a.e. on ${B_r(x_0)}$. In this case, both sides of \eqref{LH:critical-integral}
vanish, and hence the inequality holds trivially.

Consider now the case $\mathcal U_{\gamma,\, p}>0$.  For any $s\in[0,\infty)$, set
\[
  G(s):=\frac{s^{q_{\mathrm H}}}{[\log(e+s)]^\lambda}.
\]
Moreover, for any $t\in(0,1)$, set
$$ g(t):=C_0t^{-\frac 1 {q_{\rm H}}}\left[\log(e/t)\right]^{\beta_\gamma},$$
where $C_0$ is the constant appearing on the right-hand
side of \eqref{eq-LlogL-wHI}.
Assume without loss of generality that
$C_0\ge e^{1/q_{\mathrm H}}$.
For any $t\in(0,1)$, by this assumption and $g(t)\ge C_0t^{-1 /{q_{\rm H}}}$, we have
$$
\log \left(e+g(t)\right)
\ge \log(g(t))
\ge \log\left(C_0t^{-\frac 1 {q_{\rm H}}}\right)
\ge \frac1{q_{\rm H}}
\log(e/t)$$
and, hence,
\begin{align}\label{eq-Gg}
G(g(t))
= (C_0)^{q_{\rm H}}t^{-1}\frac{\left[\log(e/t)\right]^{\beta_\gamma q_{\mathrm H}}}{\left[\log \left(e+g(t)\right)\right]^{\lambda}}
\le (C_0)^{q_{\rm H}} (q_{\rm H})^\lambda t^{-1} \left[\log(e/t)\right]^{\beta_\gamma q_{\mathrm H}-\lambda}.
\end{align}

Set
$f:= {u}/{\mathcal U_{\gamma,\, p}}$. It follows from \eqref{eq-LlogL-wHI} that for all $t\in(0,1)$,
\[
  f_{B_r(x_0)}^*(t)
  \le g(t).
\]
A direct calculation shows that, for any $0<a\le b<\infty$,
\begin{align}\label{eq-G-quasi-increasing}
\frac{G(a)}{G(b)}= \left(\frac ab\right)^{q_{\rm H}}
\left(\frac{\log(e+b)}{\log(e+a)}\right)^\lambda
&
\le \left(\frac ab\right)^{q_{\rm H}}
\left(\log\left(e+\frac b a\right)\right)^\lambda\\
&
\le \sup_{0<\tau\le 1} \tau^{q_{\rm H}} \left(\log\left(e+\frac 1\tau \right)\right)^\lambda
=:C_{q_{\rm H},\,\lambda}. \notag
\end{align}
Applying this estimate with
$a=f_{B_r(x_0)}^*(t)$ and $b=g(t)$, we obtain that for any $t\in(0,1)$,
\begin{align}\label{eq-Gf*-Gg}
G(f_{B_r(x_0)}^*(t)) \le C_{q_{\rm H},\,\lambda}\,G(g(t)).
\end{align}
Combining \eqref{eq-Gf*-Gg} with \eqref{eq-Gg}, we obtain
\begin{align*}
  \fint_{B_r(x_0)}\frac{f^{q_{\mathrm H}}}{[\log(e+f)]^\lambda}\,dx
  &=\int_0^1G(f_{B_r(x_0)}^*(t))\,dt\\
  &\le C_{q_{\rm H},\,\lambda}\,\int_0^1G(g(t))\,dt\\
  &\le C\int_0^1
           [\log(e/t)]^{\beta_\gamma q_{\mathrm H}-\lambda}
                                                       \frac{dt}{t}\\
  &=C\int_1^\infty \tau^{\beta_\gamma q_{\mathrm H}-\lambda}\,d\tau.
\end{align*}
The last integral is finite when
$\lambda>1+\beta_\gamma q_{\mathrm H}$,
which is exactly the condition in
\eqref{LH:integral-threshold}.
Finally, recalling that
$f=u/\mathcal U_{\gamma,\,p}$, we obtain
\eqref{LH:critical-integral}. This establishes Theorem~\ref{thm-wHI-qH}(iii).
\end{proof}

\subsection{Critical implies subcritical for the weak Harnack inequality}
\label{sec:subcritical-consequence}

In this subsection, we show that the
subcritical weak Harnack inequality \eqref{eq-wHI-q}, originally established
by Ok--Song \cite[Theorem~1.5]{OkSong} (see also \eqref{regularity:weak-Harnack}), is a direct consequence of each of the endpoint weak Harnack
inequalities in Theorems~\ref{thm-wHI-qH-phi} and~\ref{thm-wHI-qH}.

\begin{corollary}\label{cor-wHI-qH=>q}
Let $p\in(1,n)$, $q\in(0,q_{\mathrm H})$ and $\phi$ be a nearly
critical kernel. Then, there exists a constant
$C=C(n,p,\phi,q)>0$ such that for any ball $B\subset\rn$ with radius $r>0$
and any measurable
function $f$ on $B$,
\begin{equation}\label{eq-q<qH-lorentz}
\left(\fint_{B}|f(x)|^q\,dx\right)^{1/q}
\le C \sup_{0<t<1}t^{1/q_{\mathrm H}}
     \left(\frac{\Phi(rt^{1/n})}{\Phi(r)}\right)^{1/(p-1)}
                          f_B^*(t).
\end{equation}
In particular, under the hypothesis of Theorem \ref{thm-wHI-qH-phi},
it is a direct consequence of \eqref{eq-q<qH-lorentz} and \eqref{intro:endpoint-phi} that the
supersolution $u$ to \eqref{eq-Lgammap-phi} satisfies
 \eqref{eq-wHI-q}.
\end{corollary}

\begin{proof}
It suffices to prove \eqref{eq-q<qH-lorentz}.
Denote by $\mathcal C_f$ the supremum on the right-hand side of \eqref{eq-q<qH-lorentz}. Then, for any $0<t<1$,
\[
f_{B}^*(t)
\le
\mathcal C_f\,t^{-1/q_{\mathrm H}}
\left(
\frac{\Phi(r)}
{\Phi(rt^{1/n})}
\right)^{1/(p-1)}.
\]
Thus, by the equimeasurability of $|f|$ and its nonincreasing
rearrangement $f_B^*$, we have
\begin{align*}
\left(\fint_{B}|f|^q\,dx\right)^{1/q}
&=
\left(\int_0^1[f_{B}^*(t)]^q\,dt\right)^{1/q}
\notag\\
&\le
\mathcal C_f\,
\left(
\int_0^1t^{-q/q_{\mathrm H}}
\left(
\frac{\Phi(r)}
{\Phi(rt^{1/n})}
\right)^{q/(p-1)}\,
dt
\right)^{1/q}.
\end{align*}
 For all $t\in(0,1)$ and $r\in(0,\infty)$, by
\eqref{eq-Phi-power-comparison}, we have
\[
\frac{\Phi(r)}
{\Phi(rt^{1/n})}
 \le C_\varepsilon \left(\frac {1}{t^{1/n}}\right)^\varepsilon,
\]
where we choose $\varepsilon$ satisfying
\[
0<\varepsilon<
\min\left\{1,\;n(p-1)\left(\frac1q-\frac1{q_{\mathrm H}}\right)\right\}.
\]
From this and  $q<q_{\mathrm H}$, it follows that
\[
\begin{aligned}
&\int_0^1
t^{-q/q_{\mathrm H}}
\left(\frac{\Phi(r)}
{\Phi(rt^{1/n})}\right)^{q/(p-1)}\,dt
\le
C_\varepsilon\int_0^1
t^{-q/q_{\mathrm H}}\left(\frac {1}{t^{1/n}}\right)^{
\frac{\varepsilon q}{p-1}}
\,dt
<\infty.
\end{aligned}
\]
Consequently,
\[
\left(\fint_{B}|f|^q\,dx\right)^{1/q}
\le C\mathcal C_f.
\]
This completes the proof of \eqref{eq-q<qH-lorentz}.
\end{proof}

 By Corollary~\ref{cor-wHI-qH=>q}, when $\phi=\phi_\gamma$, the endpoint weak Harnack inequality \eqref{intro:endpoint}
in Theorem~\ref{thm-wHI-qH}(i) directly implies the subcritical
weak Harnack inequality \eqref{eq-wHI-q}.

\begin{corollary}
Let $p\in(1,n)$, $\gamma\in(1,\infty)$, $q\in(0,q_{\mathrm H})$ and $\beta_\gamma:=(\gamma-1)/(p-1)$. Then, there exists a constant
$C=C(n,p,\gamma,q)>0$ such that for any ball $B\subset\rn$ with radius $r>0$
and any measurable
function $f$ on $B$,
\begin{equation}\label{eq:weak-LZ-to-Lq}
\left(\fint_B |f|^q\,dx\right)^{1/q}
\le C\,
\|f\|_{L^{q_{\mathrm H},\infty}(\log L)^{-\beta_\gamma}(B)}.
\end{equation}
In particular, under the hypothesis of Theorem \ref{thm-wHI-qH}(ii),
it is a direct consequence of \eqref{eq:weak-LZ-to-Lq} and \eqref{intro:LZ-functional} that  the
supersolution $u$ to \eqref{eq-Lgammap} satisfies  \eqref{eq-wHI-q}, with $\mathcal U_{\phi,\,p}$ replaced by $\mathcal U_{\gamma,\,p}$ therein.
\end{corollary}

\begin{proof}
We only need to prove \eqref{eq:weak-LZ-to-Lq}.
By \eqref{eq-LogL-fz},
for every $t\in(0,1)$,
\[
f_B^*(t)
\le
\|f\|_{L^{q_{\mathrm H},\infty}
(\log L)^{-\beta_\gamma}(B)}
t^{-1/q_{\mathrm H}}
[\log(e/t)]^{\beta_\gamma}.
\]
Again, using the equimeasurability of $|f|$ and its nonincreasing
rearrangement $f_B^*$, we obtain
\begin{align*}
\left(\fint_B|f|^q\,dx\right)^{1/q}
&=
\left(\int_0^1[f_B^*(t)]^q\,dt\right)^{1/q}\\
&\le
\|f\|_{L^{q_{\mathrm H},\infty}
(\log L)^{-\beta_\gamma}(B)}\,
\left(
\int_0^1
t^{-q/q_{\mathrm H}}
[\log(e/t)]^{\beta_\gamma q}\,dt
\right)^{1/q}.
\end{align*}
Since $q<q_{\mathrm H}$, the last integral is finite,
which proves \eqref{eq:weak-LZ-to-Lq}.
\end{proof}

\begin{corollary}
Let $0<q<q_0<\infty$ and $\lambda>0$. Then, there exists a constant
$C=C(q, q_0,\lambda)>0$ such that for any ball $B\subset\rn$, any measurable
function $f\ge0$ on $B$ and any $\mathcal U>0$,
\begin{equation}\label{eq:weak-Log-to-Lq}
\fint_B\left(\frac{f}{\mathcal U}\right)^q\,dx
\le
1+C\fint_B
\frac{(f/\mathcal U)^{q_0}}{[\log(e+f/\mathcal U)]^\lambda}\,dx.
\end{equation}
In particular, under the hypothesis of Theorem \ref{thm-wHI-qH}(iii),
it is a direct consequence of \eqref{eq:weak-Log-to-Lq} and \eqref{LH:critical-integral} that the
supersolution $u$ to \eqref{eq-Lgammap} satisfies  \eqref{eq-wHI-q}, with $\mathcal U_{\phi,\,p}$ replaced by $\mathcal U_{\gamma,\,p}$ therein.
\end{corollary}

\begin{proof}
It suffices to prove \eqref{eq:weak-Log-to-Lq}. Since $q<q_0$, we have
\[
C_0:=\sup_{t\ge1}
t^{q-q_0}[\log(e+t)]^\lambda<\infty.
\]
Thus, for $t\ge1$,
\[
t^q\le C_0\frac{t^{q_0}}{[\log(e+t)]^\lambda},
\]
whereas $t^q\le1$ for $0\le t\le1$.
Consequently, for every $t\ge0$,
\begin{equation}\label{preceding pointwise inequality}
t^q\le
1+C_0\frac{t^{q_0}}{[\log(e+t)]^\lambda}.
\end{equation}
Applying \eqref{preceding pointwise inequality}
with $t=f/\mathcal U$ and integrating over $B$, we obtain
\[
\fint_B\left(\frac{f}{\mathcal U}\right)^q\,dx
\le
1+C_0\fint_B
\frac{(f/\mathcal U)^{q_0}}{[\log(e+f/\mathcal U)]^\lambda}\,dx.
\]
This proves \eqref{eq:weak-Log-to-Lq}.
\end{proof}


\section{Sharpness}\label{sec:proof-sharpness}

This section aims to establish the four sharpness results stated in Theorem~\ref{intro:sharpness}.
Throughout the whole section, we simply write $B_r:=B_r(0)$ for all $r\in(0,\infty)$, and let
\[
  \Omega_0=B_{1/2} ,\qquad 0<\rho<2^{-12}.
\]
 Fix $p\in(1,n)$ and $\gamma\in(1,\infty)$, and $\phi_\gamma$ be the logarithmic kernel  defined in \eqref{eq-phi-gamma}.
The idea is to consider
the obstacle problem developed in Section~\ref{ss:obstacle}
for the operator $\mathcal L_{\gamma,\,p}:=\mathcal L_{\phi_\gamma,\,p}$,  with respect to the
admissible class
\begin{align}\label{eq-Arho}
\mathcal A_\rho
:=
\left\{
w\in\mathcal F^{\phi_\gamma,\,p}(\Omega_0):\
w=0\ \,\text{a.e. on }\,\Omega_0^c,\ \,
w\ge1\ \,\text{a.e. on }\,B_\rho
\right\}.
\end{align}
Note that $\mathcal A_\rho\ne\emptyset$, since it contains smooth cutoff functions.
The class $\mathcal A_\rho$ is convex and weakly closed in $\mathcal F_0^{\phi_\gamma,\,p}(\Omega_0)$, the direct method and the uniqueness argument of Section~\ref{ss:obstacle} give a unique energy minimizer $w_\rho$. Moreover, by Theorem~\ref{thm:obstacle-supersolution}, $w_\rho$ is a weak supersolution to $\mathcal L_{\gamma,\,p}u=0$ in $\Omega_0$.
%

To simplify the notation, in what follows we write  $\mathcal F^{\gamma,\,p}:=\mathcal F^{\phi_\gamma,\,p}$
and $\mathcal E_{\gamma,\,p}:=\mathcal E_{\phi_\gamma,\,p}$.
For each $\rho$, define the normalized function
\begin{align}\label{eq-urho}
  u_\rho:=\mathcal E_{\gamma,\,p}(w_\rho;\,\Omega_0)^{-1/(p-1)}w_\rho,
\end{align}
which is also a weak supersolution in $\Omega_0$.
For the family
$\{u_\rho\}_{\rho}$, we will show that the right-hand sides of \eqref{intro:endpoint},
\eqref{intro:LZ-functional} and
\eqref{LH:critical-integral} are uniformly  bounded in $\rho$, while any strict improvement of the
involved indices forces the left-hand side of the corresponding inequality to
diverge to infinity as $\rho\downarrow0$. This will demonstrate that all the
indices in Theorem~\ref{intro:sharpness} are optimal.

\subsection[Four technical lemmas on $w_\rho$]{Four technical  lemmas on $w_\rho$}
\label{sec:sharpness-construction}
Starting from the obstacle minimizer $w_\rho$, we
establish four  estimates needed for the
sharpness arguments.

\begin{lemma}
The minimizer  $w_\rho$ of
$\mathcal A_\rho$ satisfies
    \begin{equation}\label{LH:obstacle-bounds}
0\le w_\rho\le1,\qquad
w_\rho=1\ \, \text{a.e. on }\, B_\rho ,\qquad
w_\rho=0\ \, \text{a.e. on }\, \Omega_0^c.
\end{equation}
\end{lemma}

\begin{proof}
For any $t\in\mathbb R$, let $P(t):=\min\{t_+,1\}$.
By \eqref{eq-Arho} and $P(0)=0$, we have
$$P(w_\rho)=1\ \, \text{a.e. on }\, B_\rho ,\qquad
P(w_\rho)=0\ \, \text{a.e. on }\, \Omega_0^c.$$
By the $1$-Lipschitz property of $P$,
\begin{align*}
\mathcal E_{\gamma,\,p}(P(w_\rho);\Omega_0)
&=\iint_{\mathcal C_{\Omega_0}}
|P(w_\rho(x))-P(w_\rho(y))|^p
K_{\gamma,\,p}(x,y)\,dy\,dx\\
&\le
\mathcal E_{\gamma,\,p}(w_\rho;\,\Omega_0).
\end{align*}
This shows that
$P(w_\rho)\in \mathcal F^{\gamma,\,p}(\Omega_0)$,
 and that $P(w_\rho)$ is a minimizer of $\mathcal A_\rho$.
By uniqueness, we obtain $P(w_\rho)=w_\rho$ a.e. on $\mathbb R^n$,
which yields \eqref{LH:obstacle-bounds}.
\end{proof}

\begin{lemma}
There exists a constant $C=C(n,p,\gamma,\Lambda)>1$ such that for any $\rho\in(0, 2^{-12})$,
\begin{equation}\label{LH:obstacle-energy}
C^{-1}\rho^{n-p}\Phi_\gamma(\rho)\le   \mathcal E_{\gamma,\,p}(w_\rho;\,\Omega_0)
                        \le C\rho^{n-p}\Phi_\gamma(\rho).
\end{equation}
\end{lemma}

\begin{proof}
Note that $w_\rho^*({|B_\rho |}/{2})=1$ by  \eqref{LH:obstacle-bounds}.
Applying Proposition~\ref{logsp:fullscale} to $w_\rho$ at $t=|B_\rho |/2$ and using
\eqref{eq-doubling-Phi-gamma}, we obtain
\begin{equation*}
  \rho^{n-p}\Phi_\gamma(\rho)\le C\, \mathcal E_{\gamma,\,p}(w_\rho;\,\Omega_0),
\end{equation*}
which gives the lower bound in \eqref{LH:obstacle-energy}.

For the upper bound, we fix $\eta_0\in C_c^\infty(B_2 )$ with $0\le\eta_0\le1$
and $\eta_0=1$ on $B_1 $, and set $\eta_\rho(x)=\eta_0(x/\rho)$.
For any $h\in\rn$, the triangle inequality and translation
invariance yield
\[
\|\eta_\rho(\cdot+h)-\eta_\rho(\cdot)\|_{L^p(\rn)}
\le 2\|\eta_\rho\|_{L^p(\rn)}
\le C\rho^{n/p},
\]
and, moreover,  the fundamental theorem of calculus
and the Minkowski inequality give
\[
\|\eta_\rho(\cdot+h)-\eta_\rho(\cdot)\|_{L^p(\rn)}
=\left\|\int_0^1 h\cdot\nabla \eta_\rho(\cdot +th)\,dt\right\|_{L^p(\rn)}
\le |h|\,\|\nabla\eta_\rho\|_{L^p(\rn)}
\le C|h|\rho^{n/p-1}.
\]
Combining these estimates, we obtain
\[
  \|\eta_\rho(\cdot+h)-\eta_\rho(\cdot)\|_{L^p(\rn)}^p
  \le C\min\{\rho^{n-p}|h|^p,\,\rho^n\}.
\]
 By the admissibility of $\eta_\rho$, Lemma~\ref{kernel:global-cutoff}, and \eqref{eq-Phigamma0}, we deduce
\begin{align*}
  \mathcal E_{\gamma,\,p}(w_\rho;\,\Omega_0)
  &\le\mathcal E_{\gamma,\,p}(\eta_\rho;\,\Omega_0)\\
  &\le\Lambda\int_{\rn} \frac{\|\eta_\rho(\cdot+h)-\eta_\rho(\cdot)\|_{L^p(\rn)}^p}{|h|^{n+p}}\, \phi_\gamma(|h|)\,dh\\
  &\le C \rho^n \int_{\rn} \frac{\min\{\rho^{-p}|h|^p,\,1\}}{|h|^{n+p}}\,\phi_\gamma(|h|)\,dh\\
  &\le C\rho^{n-p}\Phi_\gamma(\rho).
\end{align*}
This proves the upper bound in \eqref{LH:obstacle-energy}.
\end{proof}

We now derive a linear truncation energy estimate for the obstacle minimizer
$w_\rho$, which will be used in the proof of Lemma
\ref{eq-rearragement-wrho} below.

\begin{lemma}\label{lem-wrho-LTEE}
For any $k\in(0,\infty)$, define
$T_k(t):=\min\{t_+,\,k\}$.
Then
\begin{equation}\label{LH:obstacle-truncation}
\mathcal E_{\gamma,\,p}(T_k(w_\rho);\,\Omega_0)\le
\min\{1,\,k\}\,\mathcal E_{\gamma,\,p}(w_\rho;\,\Omega_0).
\end{equation}
\end{lemma}

\begin{proof}
For $k\ge1$, it follows from \eqref{LH:obstacle-bounds} that
$T_k(w_\rho)=w_\rho$, and hence \eqref{LH:obstacle-truncation} holds
automatically. Henceforth, we assume $k\in(0,1)$.

Let $\zeta\in \mathcal F^{\gamma,\,p}(\Omega_0)$ satisfy
$\zeta=0$ a.e. on $B_\rho \cup \Omega_0^c$.
Then, for every $\tau\in\mathbb R$, we have
\[
w_\rho+\tau\zeta\in\mathcal A_\rho,
\]
and, moreover, the minimality of $w_\rho$ gives
\[
\mathcal E_{\gamma,\,p}(w_\rho+\tau\zeta;\Omega_0)
\ge
\mathcal E_{\gamma,\,p}(w_\rho;\,\Omega_0).
\]
Thus $\tau=0$ is a minimum of
\[
\tau\longmapsto
\mathcal E_{\gamma,\,p}(w_\rho+\tau\zeta;\Omega_0).
\]
By the first-variation argument in the proof of
Theorem~\ref{thm-1}, we obtain
\begin{equation}\label{LH:obstacle-zero-variation}
\iint_{\mathcal C_{\Omega_0}}
J_p(w_\rho(x)-w_\rho(y))
(\zeta(x)-\zeta(y))K_{\gamma,\,p}(x,y)\,dy\,dx=0,
\end{equation}
where $J_p(t):=|t|^{p-2}t$ for all $t\in\mathbb R$.

In particular, for $0<k<1$, define
$\zeta:=T_k(w_\rho)-kw_\rho$.
Since $T_k$ is Lipschitz and $T_k(0)=0$, we have
\[
\zeta\in \mathcal F^{\gamma,\,p}(\Omega_0),
\qquad
\zeta=0\quad\text{a.e. on }\Omega_0^c.
\]
Furthermore, by \eqref{LH:obstacle-bounds} and  $k\in(0,1)$, we have
\[
T_k(w_\rho)-kw_\rho=k-k=0
\quad\text{a.e. on }B_\rho .
\]
We may therefore take $\zeta=T_k(w_\rho)-kw_\rho$ in
\eqref{LH:obstacle-zero-variation}. It follows that
\begin{align*}
\iint_{\mathcal C_{\Omega_0}}
J_p(w_\rho(x)-w_\rho(y))
\big(T_k(w_\rho(x))-T_k(w_\rho(y))\big)
K_{\gamma,\,p}(x,y)\,dy\,dx=k\,\mathcal E_{\gamma,\,p}(w_\rho;\,\Omega_0).
\end{align*}
Finally, since $T_k$ is increasing and $1$-Lipschitz,
for all $A,B\in\mathbb R$,
\[
\begin{aligned}
J_p(A-B)(T_k(A)-T_k(B))
&=|A-B|^{p-1}|T_k(A)-T_k(B)|\\
&\ge |T_k(A)-T_k(B)|^p.
\end{aligned}
\]
Combining the preceding estimates gives
$$\iint_{\mathcal C_{\Omega_0}}
|T_k(w_\rho(x))-T_k(w_\rho(y))|^p
K_{\gamma,\,p}(x,y)\,dy\,dx\le
k\mathcal E_{\gamma,\,p}(w_\rho;\,\Omega_0).$$
This proves that \eqref{LH:obstacle-truncation}
holds for every $k\in(0,\infty)$.
\end{proof}

Applying Lemma \ref{lem-wrho-LTEE}, we establish the following rearrangement-energy estimate for the obstacle minimizer $w_\rho$.

\begin{lemma}\label{eq-rearragement-wrho}
There exists a positive constant $C=C(n,p,\gamma,\Lambda)$ such that, for any $\rho\in(0, 2^{-12})$,
\begin{equation}\label{LH:obstacle-rearrangement-estimate}
\sup_{0<t<1}
t^{1/q_{\mathrm H}}
\left[\Phi_\gamma(t^{1/n})\right]^{1/(p-1)}
(w_\rho)_{B_1}^*(t)
\le
C\,\mathcal E_{\gamma,\,p}(w_\rho;\,\Omega_0)^{1/(p-1)}.
\end{equation}
\end{lemma}

\begin{proof}
For any $k\in(0,\infty)$, let
$$E_k:=\{x\in\mathbb R^n:w_\rho(x)>k\}.$$
Suppose $|E_k|>0$. Then, the definition of the
nonincreasing rearrangement gives
\[
(T_k(w_\rho))^*(|E_k|/2)=k.
\]
Applying Proposition~\ref{logsp:fullscale} to $T_k(w_\rho)$
at the rearrangement level $t=|E_k|/2$,  we obtain
\[
\begin{aligned}
k^p\left(\frac{|E_k|}{2}\right)^{p/p^*}
\Phi_\gamma\left(\left(\frac{|E_k|}{2}\right)^{1/n}\right)\le
C\,\mathcal E_{\gamma,\,p}(T_k(w_\rho);\,\Omega_0).
\end{aligned}
\]
Combining this with \eqref{eq-doubling-Phi-gamma} and
\eqref{LH:obstacle-truncation}, we further obtain
\begin{align*}
k^p|E_k|^{(n-p)/n}
\Phi_\gamma(|E_k|^{1/n})
&\le C\,\mathcal E_{\gamma,\,p}(T_k(w_\rho);\,\Omega_0)\\
&\le C\,k\,\mathcal E_{\gamma,\,p}(w_\rho;\,\Omega_0).
\end{align*}
Recalling that $q_{\mathrm H}=\frac{n(p-1)}{n-p}$,
we can rewrite this last inequality as
\begin{equation}\label{LH:obstacle-distribution-qH}
k|E_k|^{1/q_{\mathrm H}}
\left[\Phi_\gamma(|E_k|^{1/n})\right]^{1/(p-1)}\le
C\,\mathcal E_{\gamma,\,p}(w_\rho;\,\Omega_0)^{1/(p-1)}.
\end{equation}
Clearly, if $|E_k|=0$, then \eqref{LH:obstacle-distribution-qH}  holds trivially with $\Phi_{\gamma}(0)=0$.

We now pass from the distribution function estimate
\eqref{LH:obstacle-distribution-qH} to the nonincreasing rearrangement
estimate \eqref{LH:obstacle-rearrangement-estimate}.
Since $w_\rho=0$ a.e. outside $B_{1/2}$, for any
$0<t<1$ and
\[
0<k<(w_\rho)_{B_1}^*(t),
\]
the definition of $(w_\rho)_{B_1}^*$ gives
\[
\frac{|\{x\in B_1:\ w_\rho(x)>k\}|}{|B_1|}>t.
\]
In particular,
\[
|E_k|\ge |\{x\in B_1:\ w_\rho(x)>k\}|>|B_1|t=\nu_n t.
\]
Since $\Phi_\gamma$ is increasing, it follows that
\[
|E_k|^{1/q_{\mathrm H}}
\left[\Phi_\gamma(|E_k|^{1/n})\right]^{1/(p-1)}\ge
 (\nu_nt)^{1/q_{\mathrm H}}
\left[\Phi_\gamma\left((\nu_nt)^{1/n}\right)\right]^{1/(p-1)}.
\]
Substituting this into the left-hand side of \eqref{LH:obstacle-distribution-qH} and using \eqref{eq-doubling-Phi-gamma},
we obtain
\[
k t^{1/q_{\mathrm H}}
\left[\Phi_\gamma(t^{1/n})\right]^{1/(p-1)}\le
C\,\mathcal E_{\gamma,\,p}(w_\rho;\,\Omega_0)^{1/(p-1)}.
\]
Letting $k\uparrow(w_\rho)_{B_1}^*(t)$, we conclude that for every $t\in(0,1)$,
\[
t^{1/q_{\mathrm H}}
\left[\Phi_\gamma(t^{1/n})\right]^{1/(p-1)}
(w_\rho)_{B_1}^*(t)\le
C\,\mathcal E_{\gamma,\,p}(w_\rho;\,\Omega_0)^{1/(p-1)}.
\]
Taking the supremum over $t\in(0,1)$ gives \eqref{LH:obstacle-rearrangement-estimate}.
\end{proof}

\subsection{Sharpness of the critical indices}
As a consequence of the preceding four lemmas in
Section~\ref{sec:sharpness-construction}, we derive the following properties
of the normalized functions $u_\rho$ in \eqref{eq-urho},
 which provide all the properties needed for the forthcoming
sharpness arguments.

\begin{lemma}\label{lem-urho}
For any $\rho\in(0, 2^{-12})$, the function $u_\rho$
defined in \eqref{eq-urho} has the following properties:
\begin{enumerate}[\rm (i)]
    \item
 $u_\rho$ is a weak supersolution  to $\mathcal L_{\gamma,\,p}u=0$ in $\Omega_0$.

   \item  $u_\rho\ge 0$ a.e. on $\rn$.  In particular, for all $x_0\in\mathbb R^n$
    and $R>0$,
    \[
    \mathsf{Tail}^{\gamma}_p
    ((u_\rho)_-;x_0,R)=0.
    \]

\item For a.e. $x\in B_\rho $,
\begin{align*}
u_\rho(x)\approx
\rho^{-n/q_{\mathrm H}}
\left[\log(e/\rho)\right]^{\beta_\gamma},
\end{align*}
where the implicit constants are independent of $\rho$ and $x$.

\item There exist positive constants $c$ and $C$, independent of $\rho$, such that
\begin{equation}\label{LH:obstacle-infimum}
  0<c\le \operatorname*{ess\,inf}_{B_{1/8} }u_\rho\le C<\infty.
\end{equation}

\item There exists a positive constant $c$, independent of $\rho$, such that for any $t\in[(32\rho)^n,1/4]$,
\begin{equation}\label{LH:obstacle-rearrangement-profile}
  (u_\rho)_{B_{1/8} }^*(t)
  \ge c\,t^{-1/q_{\mathrm H}}\left[\log(e/t)\right]^{\beta_\gamma}.
\end{equation}
\end{enumerate}
\end{lemma}

\begin{proof}
By \eqref{LH:obstacle-bounds} and
    \eqref{LH:obstacle-energy}, together with the definition of $u_\rho$ in
    \eqref{eq-urho}, we directly obtain (i), (ii) and (iii).

Now, we show  (iv). By \eqref{eq-urho}, \eqref{LH:obstacle-rearrangement-estimate} and $\beta_\gamma=(\gamma-1)/(p-1)$, we have
\begin{align*}
  (u_\rho)_{B_1}^*(t)&=\mathcal E_{\gamma,\,p}(w_\rho;\,\Omega_0)^{-1/(p-1)}(w_\rho)_{B_1}^*(t)\\
 &\le Ct^{-1/q_{\mathrm H}}
       \left[\Phi_\gamma(t^{1/n})\right]^{-1/(p-1)}\notag\\
       &       \le Ct^{-1/q_{\mathrm H}}\left[\log(e/t)\right]^{\beta_\gamma}, \notag
  \end{align*}
  where the last step follows from \eqref{eq-Phi-value}.
Consequently, for any $q\in(0,q_{\mathrm H})$,
by the equimeasurability of $u_\rho$ and $(u_\rho)_{B_1}^*$, we obtain
\begin{align*}
\fint_{B_1}u_\rho^q\,dx
=\int_0^1[(u_\rho)_{B_1}^*(t)]^q\,dt\le C^q\int_0^1
t^{-q/q_{\mathrm H}}[\log(e/t)]^{\beta_\gamma q}\,dt,
\end{align*}
where the last integral converges
since $q<q_{\mathrm H}$.
 Thus, the estimate
\begin{equation}\label{LH:obstacle-infimum-upper}
\operatorname*{ess\,inf}_{B_{1/8} }u_\rho
\le
\left(\fint_{B_{1/8}}u_\rho^q\,dx\right)^{1/q}
\le \left(8^{n} \fint_{B_1}u_\rho^q\,dx\right)^{1/q}\le C
\end{equation}
holds uniformly in $\rho$. This proves the upper bound in \eqref{LH:obstacle-infimum}.

Next, we show the lower bound in \eqref{LH:obstacle-infimum}.
Let $4\rho\le s\le 1/8$.  Since $w_\rho=1$ a.e. on $B_\rho $ (see \eqref{LH:obstacle-bounds}), it follows that for
\[
  t_s:=\frac{|B_\rho |}{2|B_s |}=\frac12(\rho/s)^n
\]
we have
$(w_\rho)^*_{B_s }(t_s)\ge 1$ and, hence,
\[(u_\rho)^*_{B_s }(t_s)
= \mathcal E_{\gamma,\,p}(w_\rho;\,\Omega_0)^{-1/(p-1)}\, (w_\rho)^*_{B_s }(t_s)
\ge\mathcal E_{\gamma,\,p}(w_\rho;\,\Omega_0)^{-1/(p-1)}.\]
We apply the already proved
Theorem~\ref{thm-wHI-qH} to $u_\rho$ on $B_s $ with respect to the
outer ball $B_{1/4} $.
Since the negative tail vanishes, we have
\begin{align}\label{LH:obstacle-profile}
  C\operatorname*{ess\,inf}_{B_s}u_\rho
  &\ge t_s^{1/q_{\mathrm H}}
       \left(\frac{\Phi_\gamma(st_s^{1/n})}{\Phi_\gamma(s)}\right)^{1/(p-1)}
                            (u_\rho)^*_{B_s }(t_s)\\
&\ge t_s^{1/q_{\mathrm H}}
       \left(\frac{\Phi_\gamma(st_s^{1/n})}{\Phi_\gamma(s)}\right)^{1/(p-1)}
                            \mathcal E_{\gamma,\,p}(w_\rho;\,\Omega_0)^{-1/(p-1)}\notag\\
  &=c\left[\left(\frac\rho s\right)^{n-p}
              \frac{\Phi_\gamma(2^{-1/n}\rho)}{\Phi_\gamma(s)\,\mathcal E_{\gamma,\,p}(w_\rho;\,\Omega_0)}
      \right]^{1/(p-1)}\notag\\
  &\ge c[s^{n-p}\Phi_\gamma(s)]^{-1/(p-1)}, \notag
\end{align}
where the last inequality follows from \eqref{LH:obstacle-energy}  and
$\Phi_\gamma(2^{-1/n}\rho)\ge c\Phi_\gamma(\rho)$.  In particular, taking $s=\frac{1}{8}$ in
\eqref{LH:obstacle-profile} yields
\begin{align}\label{LH:obstacle-infimum-lower}
\operatorname*{ess\,inf}_{B_{1/8} }u_\rho\ge c>0.
\end{align}
A combination of \eqref{LH:obstacle-infimum-upper} and \eqref{LH:obstacle-infimum-lower}
gives \eqref{LH:obstacle-infimum}.

Finally, we prove (v). Recall that \eqref{LH:obstacle-profile}
holds whenever $4\rho\le s\le 1/8$.
For any $t\in[(32\rho)^n,1/4]$, set
\[
s_t:=\frac18(2t)^{1/n}.
\]
Then $4\rho\le s_t\le1/8$ and
$s_t^n=2t(1/8)^n$.
By the definition of the normalized rearrangement
$(u_\rho)_{B_{1/8}}^*$, we have
\begin{align*}
(u_\rho)_{B_{1/8}}^*(t)
&\ge \operatorname*{ess\,inf}_{B_{s_t}}u_\rho\\
&\ge c[s_t^{n-p}\Phi_\gamma(s_t)]^{-1/(p-1)}\\
&\ge ct^{-1/q_{\mathrm H}}
[\log(e/t)]^{\beta_\gamma},
\end{align*}
where the second inequality follows from
\eqref{LH:obstacle-profile}, while the third follows
from \eqref{eq-Phi-value}.
This proves \eqref{LH:obstacle-rearrangement-profile}.
\end{proof}


Applying Lemma \ref{lem-urho}, we now establish the sharpness of the critical indices in Theorem~\ref{intro:sharpness}.

\begin{proof}[Proof of Theorem~\ref{intro:sharpness}(i)--(iii)]
We examine the weak Harnack inequalities
\eqref{intro:endpoint} and \eqref{intro:LZ-functional}
for the weak supersolutions $u_\rho$ in the ball
$B_r(x_0)=B_{1/8}$.
By parts (ii) and (iv) of Lemma~\ref{lem-urho}, we have
\begin{align}\label{eq-U-urho}
0<c<\mathcal U_{\gamma,\,p}(u_\rho;0,1/8,1/4)
=\operatorname*{ess\,inf}_{B_{1/8}}u_\rho
\le C<\infty.
\end{align}
Thus, the right-hand sides of
\eqref{intro:endpoint} and \eqref{intro:LZ-functional}
are uniformly bounded in $\rho$ for
$B_r(x_0)=B_{1/8}$ and $R=1/4$.

By following the proof of \eqref{eq-u-LlogL}, we find that for all  $q\ge q_{\rm H}$ and $\sigma\le 1$,
\begin{align}\label{eq-u-LlogL-q-sigma}
    \|u_\rho\|_{L^{q,\infty}(\log L)^{-\sigma\beta_{\gamma}}(B_{1/8})}
& = \sup_{t\in(0,1)}
 t^{\frac 1 q}\left[\log(e/t)\right]^{-\sigma\beta_\gamma}(u_\rho)_{B_{1/8}}^*(t)\\
 &\le \sup_{t \in (0, 1)}t^{\frac 1 q}
 \left(\frac{\Phi_\gamma(8^{-1}t^{1/n})}{\Phi_\gamma(8^{-1})}\right)^{\frac{\sigma}{p-1}}
                         (u_\rho)_{B_{1/8}}^*(t).  \notag
\end{align}
For any $\rho\in(0, 2^{-12})$, set
$t_\rho:  =\frac12(64\rho)^n.$
Using the fact $(32\rho)^n\le t_\rho\le1/4$ and \eqref{LH:obstacle-rearrangement-profile},
we obtain
$$(u_\rho)_{B_{1/8} }^*(t_\rho)\ge c\rho^{-\frac{n}{q_{\mathrm H}}}
[\log(e/\rho)]^{\beta_\gamma},$$
which implies
\begin{align}\label{LH:obstacle-target-lower}
 \|u_\rho\|_{L^{q,\infty}(\log L)^{-\sigma\beta_{\gamma}}(B_{1/8})}
 &\ge t_\rho^{\frac1 q}[\log(e/t_\rho)]^{-\sigma \beta_{\gamma}}
                       (u_\rho)_{B_{1/8} }^*(t_\rho)\\
  &\ge c\rho^{\frac n q-\frac n {q_{\mathrm H}}}
                       [\log(e/\rho)]^{(1-\sigma)\beta_\gamma}. \notag
\end{align}
If $q>q_{\mathrm H}$ or $\sigma<1$, then the
right-hand side of \eqref{LH:obstacle-target-lower}, and hence each term
in \eqref{eq-u-LlogL-q-sigma}, tends to $\infty$ as $\rho\to0_+$.

Thus, if any one of the following replacements occurs:
\begin{itemize}
    \item the exponent $q_{\mathrm H}$ in \eqref{intro:endpoint} or \eqref{intro:LZ-functional} is replaced by some $q>q_{\mathrm H}$;
    \item the exponent $1/(p-1)$ in \eqref{intro:endpoint} is replaced by $\sigma/(p-1)$;
    \item the exponent $\beta_{\gamma}$ in \eqref{intro:LZ-functional} is replaced by $\sigma\beta_{\gamma}$,
\end{itemize}
then the left-hand sides of \eqref{intro:endpoint} and
\eqref{intro:LZ-functional} for $u_\rho$ become the quantities appearing
in \eqref{eq-u-LlogL-q-sigma}, which tend to $\infty$ as $\rho\to0_+$.
 This proves the sharpness results stated in Theorem~\ref{intro:sharpness}(i)--(ii).

Further, by \eqref{eq-LqsLogL} and the lower rearrangement bound \eqref{LH:obstacle-rearrangement-profile}, for any $s\in(0,\infty)$,
\begin{align*}
  \|u_\rho\|_{L^{q_{\rm H},s}(\log L)^{-\beta_{\gamma}}(B_{1/8})}^s&= \int_0^1
    \left(t^{1/q_{\mathrm H}}[\log(e/t)]^{-\beta_\gamma}
                         (u_\rho)_{B_{1/8} }^*(t)\right)^s\,
                         \frac{dt}{t}\\
                         &\ge c\int_{(32\rho)^n}^{1/4}\left(t^{1/q_{\mathrm H}}[\log(e/t)]^{-\beta_\gamma}
                         (u_\rho)_{B_{1/8} }^*(t)\right)^s\,
                         \frac{dt}{t}\\
  &\ge c\int_{(32\rho)^n}^{1/4}\,\frac{dt}{t}\\
  &   =c\log\frac1{4(32\rho)^n},
\end{align*}
which tends to $\infty$ as $\rho\to0_+$.
Combining this with \eqref{eq-U-urho}, we obtain that \eqref{intro:LZ-functional}
fails for the family $\{u_\rho\}_{\rho}$ if the space
$L^{q_{\mathrm H},\infty}(\log L)^{-\beta_\gamma}
(B_{1/8})$ is replaced by
$L^{q_{\mathrm H},s}(\log L)^{-\beta_\gamma}
(B_{1/8})$.
This completes the proof of
Theorem~\ref{intro:sharpness}(iii).
\end{proof}

\begin{proof}[Proof of Theorem~\ref{intro:sharpness}(iv).]
Consider now the weak Harnack inequality
\eqref{LH:critical-integral} for the weak supersolutions
$u_\rho$ on $B_r(x_0)=B_{1/8}$. As before, we have
\eqref{eq-U-urho}. Set
\[
f_\rho:=
\frac{u_\rho}{\mathcal U_{\gamma,\,p}(u_\rho;0,1/8,1/4)}.
\]
By \eqref{eq-U-urho} and
\eqref{LH:obstacle-rearrangement-profile}, there exists
$c_0>0$, independent of $\rho$, such that for all $t\in [(32\rho)^n, 1/4]$,
\[
(f_\rho)_{B_{1/8}}^*(t)
\ge c_0t^{-1/q_{\mathrm H}}
[\log(e/t)]^{\beta_\gamma}
=:g(t).
\]

Fix $\lambda\le1+\beta_\gamma q_{\mathrm H}$ and define
\[
G(s):=\frac{s^{q_{\mathrm H}}}{[\log(e+s)]^\lambda}
\qquad \text{for all}\ \, s\ge0.
\]
By the quasi-increasing property of $G$
(by \eqref{eq-G-quasi-increasing} for $\lambda>0$ and monotonicity for $\lambda\le0$), we obtain
\begin{align*}
G((f_\rho)_{B_{1/8}}^*(t))
&\ge c\,G(g(t))\\
&=c\frac{c_0^{q_{\mathrm H}}t^{-1}
[\log(e/t)]^{\beta_\gamma q_{\mathrm H}}}
{[\log(e+g(t))]^\lambda}\\
&\ge c\,t^{-1}
[\log(e/t)]^{\beta_\gamma q_{\mathrm H}-\lambda},
\end{align*}
where in the last step we used
$\log(e+g(t))\approx\log(e/t)$.
Then,
\begin{align*}
  \fint_{B_{1/8} }
       \frac{f_\rho^{q_{\mathrm H}}}
            {[\log(e+f_\rho)]^\lambda}\,dx
  &=\int_0^1G((f_\rho)_{B_{1/8} }^*(t))\,dt\\
   &\ge c\int_{(32\rho)^n}^{1/4}
             [\log(e/t)]^{\beta_\gamma q_{\mathrm H}-\lambda}\,
                                               \frac{dt}{t}\\
  &=c\int_{\log(4e)}^{\log(e/(32\rho)^n)}
                       \tau^{\beta_\gamma q_{\mathrm H}-\lambda}\,d\tau.
\end{align*}
The last integral diverges as $\rho\to0_+$ whenever
$\lambda\le1+\beta_\gamma q_{\mathrm H}$.
Consequently, the estimate \eqref{LH:critical-integral}
fails for the family $\{u_\rho\}_{\rho}$ if its logarithmic exponent threshold
$1+\beta_\gamma q_{\mathrm H}$ is replaced by any
$\lambda\le1+\beta_\gamma q_{\mathrm H}$.
This completes the proof of
Theorem~\ref{intro:sharpness}(iv).
\end{proof}

\medskip

\noindent{\bf Statements and Declarations.}
The authors declare no competing financial interests. Data availability is not applicable to this paper.


\end{document}